\documentclass[reqno]{amsart}
\usepackage[a4paper]{geometry}

\usepackage[T2A]{fontenc}
\usepackage[utf8]{inputenc}
\usepackage[ukrainian,russian,english]{babel}

\usepackage{comment}

\usepackage{subcaption}
\usepackage{amssymb}
\usepackage[mathscr]{eucal}
\usepackage{tikz}
\usetikzlibrary{positioning,calc,decorations.markings}
\usetikzlibrary{calc}
\usepackage{mathtools}
\usepackage{enumerate}
\usepackage{accents}
\usepackage{dsfont}
\usepackage{multicol}
\usepackage{color}
\usepackage[colorlinks=true]{hyperref}
\hypersetup{urlcolor=blue,citecolor=red,linkcolor=blue}
\usepackage[initials]{amsrefs}
\definecolor{darkgreen}{rgb}{0.13, 0.55, 0.13}
\numberwithin{equation}{section}
\theoremstyle{plain}
\newtheorem{theorem}{Theorem}[section]
\newtheorem{proposition}[theorem]{Proposition}
\newtheorem{corollary}[theorem]{Corollary}
\newtheorem{lemma}[theorem]{Lemma}
\theoremstyle{definition}
\newtheorem*{rh-pb*}{Basic RH problem}
\newtheorem*{rh-I*}{RH problem normalized as $I$ at $\infty$}
\newtheorem*{sol-rh-pb*}{Soliton RH problem}
\newtheorem*{data*}{Data of this RH problem associated with $\BS{u_0(x)}$}
\theoremstyle{remark}
\newtheorem{remark}[theorem]{Remark}
\newtheorem*{notations*}{Notations}
\providecommand{\BS}[1]{\boldsymbol{#1}}  
\providecommand{\D}[1]{\mathbb{#1}}
\newcommand{\dd}{\mathrm{d}}
\newcommand{\eul}{\mathrm{e}}
\newcommand{\ii}{\mathrm{i}}
\newlength{\dhatheight}
\newcommand{\doublehat}[1]{%
    \settoheight{\dhatheight}{\ensuremath{\hat{#1}}}%
    \addtolength{\dhatheight}{-0.3ex}%
    \hat{\vphantom{\rule{1pt}{\dhatheight}}%
    \smash{\hat{#1}}}}
\renewcommand{\Im}{\operatorname{Im}}

\newcommand{\ord}{\mathrm{O}}
\DeclareMathOperator{\Res}{Res}
\newcommand{\red}{\textcolor{red}}

\newif\ifshort
\shorttrue
\title{The Short Pulse equation on the half lines: a Riemann--Hilbert approach}
\author[Iryna Karpenko]{Iryna Karpenko}
\address{IK: Faculty of Mathematics\\ University of Vienna\\
Oskar-Morgenstern-Platz 1\\ 1090 Wien\\ Austria\\ and B. Verkin Institute for Low Temperature Physics and Engineering\\ 47, Nauky ave\\ 61103 Kharkiv\\ Ukraine}
\email{\href{mailto:iryna.karpenko@univie.ac.at}{iryna.karpenko@univie.ac.at}}

\author[Dmitry Shepelsky]{Dmitry Shepelsky}
\address{DS: B.Verkin Institute for Low Temperature Physics and Engineering, Kharkiv, Ukraine\\ and V.N.Karazin Kharkiv National University, Ukraine}
\email{\href{mailto:shepelsky@yahoo.com }{shepelsky@yahoo.com}}

\begin{document}

\begin{abstract} 
In this work, we analyze the initial boundary value (IBV) problems 
for the Short Pulse (SP) equation  $u_{xt}=u+\frac{1}{6}(u^3)_{xx}$
(classical solutions)
on the half-lines $x\ge 0$ and $x\le 0$ under the condition of 
sufficient decay as $x\to +\infty$ or as $x\to -\infty$.
Using the Riemann--Hilbert (RH) problem approach, we show
that these two IBV problems differ substantially 
regarding the boundary values required for problems to be well-posed. Namely, the solution of the ``right problem'' ($x\ge 0$)
is uniquely determined by the initial data $u(x,0)$, $x\ge 0$ alone, whereas for the ``left problem'' ($x\le 0$), two boundary
functions $u(0,t)$ and $u_x(0,t)$ have to be prescribed, together 
with the initial data, for a well-posed problem.

\end{abstract}

\maketitle

\section{Introduction}

Introduced approximately fifty years ago, the Inverse Scattering Transform (IST) provides a rigorous analytic framework for solving a distinguished class of one-dimensional nonlinear evolution equations, often referred to as \textit{integrable equations}. The IST method is fundamentally built on the analysis of an associated Lax pair, a pair of linear operators that typically depend on a spectral parameter and whose compatibility condition yields the nonlinear equation \cites{IST01,IST02,ablowitz1974inverse, ablowitz1991solitons}. Originally developed for initial-value problems on the real line with rapidly decaying data, the method reduces the direct spectral analysis to a classical scattering problem, while the inverse step, addressing the reconstruction of the solution from scattering data, is traditionally formulated via the Gelfand–Levitan–Marchenko (GLM) integral equations \cite{novikov1984theory} or, alternatively, as a Riemann–Hilbert (RH) factorization problem in the complex spectral plane \cite{trogdon2015riemann}. 

Extending the IST to nonlinear evolution equations posed on domains with boundaries, such as the half-line or finite intervals, introduces substantial analytical challenges stemming from the incorporation of boundary conditions. These difficulties were effectively addressed by the development of the Unified Transform Method (UTM), also known as the Fokas Method, introduced by Fokas in \cite{fokas1997} and subsequently advanced by many authors \cites{fokas2008unified,BS08, FIS,BFS06,KS25}. The UTM generalizes the IST by simultaneously analyzing both equations of the Lax pair, and a central component of the method is the \textit{global relation}, which links the spectral functions associated with the initial and boundary data.

A natural motivation for studying initial boundary value problems for integrable equations stems from their interpretation as models of wave propagation in physical media. In such contexts, spatial boundaries represent physical interfaces, where wave profiles are generated or prescribed. The imposed boundary conditions therefore characterize the input signals entering the medium, while the integrable dynamics govern their subsequent propagation, interaction, and nonlinear evolution.

A particularly important integrable model in nonlinear optics is the Short Pulse (SP) equation, which becomes relevant when the classical nonlinear Schrödinger (NLS) approximation is no longer adequate. The NLS equation successfully describes the slow modulation of broad, weakly nonlinear wave packets whose spectra are narrowly concentrated around a carrier frequency. However, modern high-speed fiber-optic communication relies on ultra-short pulses, a regime in which the assumptions underlying the NLS approximation break down. To address this, Schäfer and Wayne \cites{SW04,CJWS05} proposed the Short Pulse equation
\begin{equation}\label{SP_1}
u_{xt} = u + \frac{1}{6}(u^{3})_{xx},
\end{equation}
as an alternative model capable of accurately capturing the evolution of ultra-short intense infrared pulses in silica optics. Numerical simulations in \cite{CJWS05} demonstrated that the SP equation provides a reliable description even for pulses with broad spectra. 

The SP equation has been studied from several perspectives. Solitary-wave solutions have been constructed in \cites{Gambino2015,Matsuno2007,Sakovich2006}, and periodic solutions have been presented in \cite{Matsuno2008}. The well-posedness of the Cauchy problem has been studied in \cites{Coclite2015,Pelinovsky2010,SW04}, and sufficient conditions for wave breaking were derived in \cite{Liu2009}. 

The Short Pulse equation is formally integrable: it is the compatibility condition for the Lax pair (\cite{SS05}):
\begin{subequations}\label{Lax}
\begin{align}\label{Lax-x}
\Phi_x&=U\Phi,\\
\label{Lax-t}
\Phi_t&=V\Phi,
\end{align}
\end{subequations}
where the coefficients $U\equiv U(x,t,\lambda)$ and $V\equiv V(x,t,\lambda)$ are defined by
\begin{subequations}\label{Lax-UV}
\begin{align}\label{Lax-U}
U&=\begin{pmatrix} \lambda & \lambda u_x \\
\lambda u_x & -\lambda \end{pmatrix},\\
\label{Lax-V}
V&=\begin{pmatrix}\frac{\lambda}{2}u^2+\frac{1}{4\lambda} &
\frac{\lambda}{2}u^2u_x-\frac{1}{2}u\\
\frac{\lambda}{2}u^2u_x+\frac{1}{2}u & -\frac{\lambda}{2}u^2-\frac{1}{4\lambda}\end{pmatrix}.
\end{align}
\end{subequations}
The Riemann--Hilbert formalism for the Cauchy problem for \eqref{SP_1} on the zero background was developed in \cite{BSL17}.

The SP equation admits the equivalent ``conservation-law'' formulation
\begin{subequations}\label{SP}
\begin{align}\label{SP-2}
&q_t = \frac{1}{2}(u^2q)_x,\\
\label{tq}
&q=\sqrt{1+w^2},\\
\label{tw}
&w= u_x.
\end{align}
\end{subequations}

In what follows, we will denote 
\[
H^{2,2}(\Omega)=\{f(x)\in L^2(\Omega):~ x^2f(x), f'(x), x^2f'(x),f''(x), x^2f''(x) \in L^2(\Omega)\},
\]
\[
\tilde H^{2,2}(\Omega):=\{f(x)\in L^1_{loc}(\Omega): f'(x), x^2f'(x),f''(x), x^2f''(x) \in L^2(\Omega)\}.
\]

In the present paper, we consider two initial-boundary value (IBV) problems for  \eqref{SP}, 
on the half-lines $x\ge 0$ (Problem I) and $x\le 0$ (Problem II), 
searching for solutions $u(x,t)$ such that $u(\cdot, t)\in \tilde H^{2,2}(0,\infty)$ for Problem I or $u(\cdot, t)\in \tilde H^{2,2}(-\infty, 0)$
for Problem II 
for all $t\in[0,T]$ with $ 0 <T<\infty$.
We assume that 
the initial conditions
\begin{equation}\label{ic}
     u(x,0) = u_0(x), \quad \pm x \geq 0,
\end{equation}
are given, and our primary question is:
What are the conditions at the boundary  $x=0$ that make the problems well-posed?

The outline of the paper is as follows.

In Section \ref{sec:2}, we present the basic tools
for our study that stem from the Lax pair representations of the SP
equation: 
 Jost solutions (``eigenfunctions'') and the corresponding spectral functions (scattering coefficients).
 We discuss their  analytic properties in view of using them
 to construct the Riemann--Hilbert factorization problems,
 whose solutions can be used to obtain 
 representations of solutions of Problem I and Problem II. 

 In Section 3, we construct master RH problems and discuss the uniqueness of solutions for Problem I and Problem II. In particular,  we find
that the initial data alone determine uniquely its solution for Problem I.

 In Section 4, we discuss the compatibility conditions for the initial and boundary values (``global relations'').

 In Section 5, we describe the inverse spectral problems for the initial and boundary values in terms of solutions of the associated RH problems.

In Section 6, we first discuss 
the relationship between solutions of the SP equation in the original variables $(x,t)$
and solutions of the SP equation in variables $(y,t)$ suitable for the RH formalism.
Then we discuss 
how a function solving, locally, the SP equation 
in $(y,t)$ variables appears from the solution of a Riemann-Hilbert
problem parametrized by $y$ and $t$.

In Section 7, discuss the existence of solutions for Problem I and Problem II. In particular, 
 we reveal that for Problem II the boundary functions 
$u(0,t)$ and $u_x(0,t)$ can be given independently.

%-------------------%
\begin{notations*}
In what follows, $\sigma_1\coloneqq\left(\begin{smallmatrix}0&1\\1&0\end{smallmatrix}\right)$, $\sigma_2\coloneqq\left(\begin{smallmatrix}0&-\ii\\\ii&0\end{smallmatrix}\right)$, and $\sigma_3\coloneqq\left(\begin{smallmatrix}1&0\\0&-1\end{smallmatrix}\right)$ denote the standard Pauli matrices, 
$e^{\hat\sigma_3}A:=e^{\sigma_3}A e^{-\sigma_3}$,
$\mathbb{C}^+:=\{\lambda\in\mathbb{C}|\Im(\lambda)> 0\}$, and 
$\mathbb{C}^-:=\{\lambda\in\mathbb{C}|\Im(\lambda)< 0\}$. We also let $f^*(k)\coloneqq\overline{f(\bar k)}$ denote the Schwarz 
reflection of a function $f(k)$, $k\in\D{C}$.
\end{notations*}
%-------------------%

\section{Eigenfunctions and Spectral Functions}\label{sec:2}

\subsection{Problem I}\label{subsec:2.1}

Given a solution $u(x,t)$ of the SP equation on $x\in(0,\infty)$, $t\in(0,T)$, with $u(x,t)\to 0$ as $x\to +\infty$ fast enough for each fixed $t$, the first step is to analyze the analytic properties of eigenfunctions associated with the Lax pair in the complex spectral plane. This analysis provides the relations among the eigenfunctions that can be interpreted as a matrix Riemann--Hilbert factorization problem.

Notice that the coefficient matrices $U$ and $V$ given in \eqref{Lax-UV} are traceless, which implies that the determinant of any matrix solution of \eqref{Lax} (constructed from two linearly independent vector solutions) is independent of both $x$ and $t$.

To align the presentation with that used for other CH-type equations (see \cites{BS08,BKS20}), it is convenient to introduce the spectral parameter $k := \ii\lambda$.

Since $U$ and $V$ have singularities in the extended complex $\lambda$--plane at $k = 0$ and $k = \infty$, particular care is required when describing the asymptotic behavior of eigenfunctions near these points. An important observation is that at $\lambda=0$ the matrix $U$ becomes independent of $u$, a property shared by many Camassa--Holm--type systems. This suggests that the behavior as $k \to 0$ can be used to recover the solution of the nonlinear equation via the associated Riemann--Hilbert problem (cf. \cite{BSL17}).

\subsubsection{Eigenfunctions near $k=\infty$}\label{sec:2.1}
In order to control the large $k$ behavior, we introduce 
\[
\tilde \Phi_\infty(x,t,k):= P(x,t)\Phi(x,t,k),
\]
where
\begin{equation}\label{P}
    P(x,t):=\sqrt{\frac{1+q(x,t)}{2q(x,t)}}\begin{pmatrix}
1 & \frac{w(x,t)}{1+q(x,t)}  \\
-\frac{w(x,t)}{1+q(x,t)} & 1 \\
\end{pmatrix},
\end{equation}
where $w(x,t)=u_x(x,t)$ and $q(x,t)=\sqrt{1+w^2(x,t)}$. This transforms \eqref{Lax} into

\begin{subequations}\label{Lax-Q-form}
\begin{align}\label{Lax-Q-form_1}
&\tilde \Phi_{\infty x}+\ii k q\sigma_3\tilde \Phi_\infty = U_\infty \tilde \Phi_\infty,\\\label{Lax-Q-form_2}
&\tilde \Phi_{\infty t} +\ii k\left(\frac{u^2q}{2}-\frac{1}{4k^2}\right)\sigma_3\tilde \Phi_\infty = V_\infty \tilde \Phi_\infty,
\end{align}

where $U_\infty\equiv U_\infty(x,t,k)$ and $V_\infty\equiv V_\infty(x,t,k)$ are given by 

\begin{equation}\label{U-hat}
U_\infty=\frac{u_{xx}}{2q^2}
\begin{pmatrix}
0 & 1 \\
-1 & 0 \\
\end{pmatrix}
\end{equation}
and 
\begin{equation}\label{hat-V}
\begin{aligned}
 V_\infty&=
-\frac{1}{4\ii k q}
\begin{pmatrix}
1-q & -w \\-w & -1+q
\end{pmatrix}
+\frac{ u^2 u_{xx}}{4 q^2}  \begin{pmatrix}
0 & 1 \\
-1 & 0
\end{pmatrix}.
\end{aligned}
\end{equation}
\end{subequations}

Define $Q(x,t,k)$ by 
\begin{subequations}\label{Qp}
\begin{equation}\label{Q}
Q(x,t,k):= \ii k p(x,t,k)\sigma_3, 
\end{equation}
with
\begin{equation}\label{p}
p(x,t,k):=\int_0^{x} q(\xi,t)\dd\xi+\int_0^{t} \left(\frac{u^2q}{2}\right)(0,\tau)\dd\tau-\frac{t}{4 k^2}.
\end{equation}
\end{subequations}
Then $p$ satisfies
\[
p_x=q,\qquad
p_t=\frac{u^2q}{2}-\frac{1}{4k^2},\qquad p(0,0,k)=0,
\]
and 
equations in \eqref{Lax-Q-form} take the form
\begin{subequations}\label{phi-tilde}
   \begin{align}
&\tilde \Phi_{\infty x}+Q_x\tilde \Phi_\infty = U_\infty \tilde \Phi_\infty,\\
&\tilde \Phi_{\infty t} +Q_t\tilde \Phi_\infty = V_\infty \tilde \Phi_\infty,
\end{align} 
\end{subequations}
which is appropriate to control the large $k$ behavior of 
its solutions.

\begin{remark}
    Notice that in the case of Problem I, \eqref{SP} implies the ``conservation law''
\begin{equation} \label{cons_law}
\nu(t)=\nu(0)-\eta(t),  
\end{equation}
where

\[
\nu(t):=\int_0^\infty (q -1)(\xi,t)\dd \xi,\qquad 
    \eta(t):=\frac{1}{2}\int_0^t (u^2q)(0,\tau)\dd \tau.
\]
    
\end{remark}

\begin{remark}
Under the condition that $q(\cdot,t)-1\in L^1(0,\infty)$ for all $t$
(which holds for $u(\cdot,t)\in H^{2,2}(0,\infty)$),
the function $p(x,t,k)$ can be written as 
    \begin{equation}\label{p-alt}
        p(x,t,k)=-\int_x^{\infty} (q-1)(\xi,t)\dd\xi+x+\int_0^{\infty} (q-1)(\xi,0)\dd\xi-\frac{t}{4 k^2}.
    \end{equation}
   \end{remark}

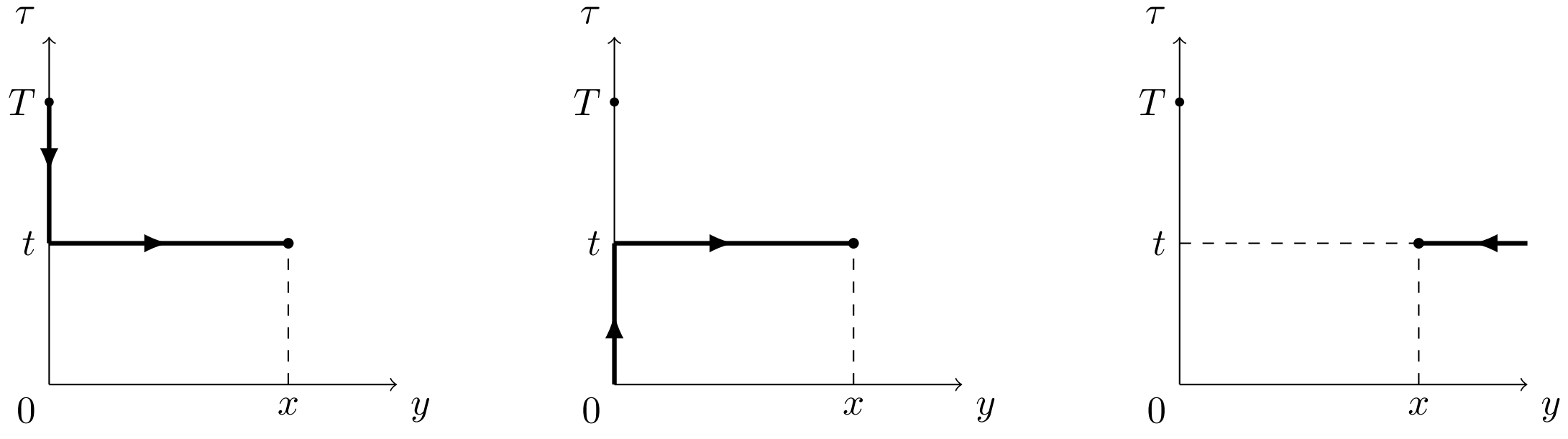
\begin{figure}[ht]
\centering
\begin{tikzpicture}[scale=1]
  % Tunables
  \def\x{2.2}  % x-location of the evaluation point
  \def\t{1.3}  % t level
  \def\T{2.6}  % mark on tau-axis
  \def\W{3.2}  % panel width

  % ---------- Panel 1 ----------
  \begin{scope}
    % axes
    \draw[->] (0,0) -- (\W,0) node[below right] {$y$};
    \draw[->] (0,0) -- (0,\W) node[above left] {$\tau$};

    % labels
    \node[below left] at (0,0) {$0$};
    \node[below]      at (\x,0) {$x$};
    \node[left]       at (0,\t) {$t$};

    % T mark
    \fill (0,\T) circle (1.2pt) node[left] {$T$};

    % dashed guide
    \draw[dashed] (\x,0) -- (\x,\t);

    % path: down then right (arrow in the middle)
    \draw[very thick,postaction={decorate},
          decoration={markings, mark=at position 0.5 with {\arrow{latex}}}] 
          (0,\T) -- (0,\t);
    \draw[very thick,postaction={decorate},
          decoration={markings, mark=at position 0.5 with {\arrow{latex}}}] 
          (0,\t) -- (\x,\t);

    % evaluation point
    \fill (\x,\t) circle (1.4pt);
  \end{scope}

  % ---------- Panel 2 ----------
  \begin{scope}[xshift=5.2cm]
    \draw[->] (0,0) -- (\W,0) node[below right] {$y$};
    \draw[->] (0,0) -- (0,\W) node[above left] {$\tau$};

    \node[below left] at (0,0) {$0$};
    \node[below]      at (\x,0) {$x$};
    \node[left]       at (0,\t) {$t$};

    \fill (0,\T) circle (1.2pt) node[left] {$T$};

    \draw[dashed] (\x,0) -- (\x,\t);

    % path: up then right (arrow in the middle)
    \draw[very thick,postaction={decorate},
          decoration={markings, mark=at position 0.5 with {\arrow{latex}}}] 
          (0,0) -- (0,\t);
    \draw[very thick,postaction={decorate},
          decoration={markings, mark=at position 0.5 with {\arrow{latex}}}] 
          (0,\t) -- (\x,\t);

    \fill (\x,\t) circle (1.4pt);
  \end{scope}

  % ---------- Panel 3 ----------
  \begin{scope}[xshift=10.4cm]
    \draw[->] (0,0) -- (\W,0) node[below right] {$y$};
    \draw[->] (0,0) -- (0,\W) node[above left] {$\tau$};

    \node[below left] at (0,0) {$0$};
    \node[below]      at (\x,0) {$x$};
    \node[left]       at (0,\t) {$t$};

    \fill (0,\T) circle (1.2pt) node[left] {$T$};

    \draw[dashed] (\x,0) -- (\x,\t);
    \draw[dashed] (0,\t) -- (\x,\t);

    % path: from right to (x,t), arrow in the middle pointing left
    \draw[very thick,postaction={decorate},
          decoration={markings, mark=at position 0.5 with {\arrow{latex}}}] 
          (\W,\t) -- (\x,\t);

    \fill (\x,\t) circle (1.4pt);
  \end{scope}
\end{tikzpicture}

\caption{Paths of integration for $\Phi_{\infty1}, \Phi_{\infty2},$ and $\Phi_{\infty3}$ 
($\Phi_{01}, \Phi_{02},$ and $\Phi_{03}$).}
\label{fig:integration-paths}
\end{figure}

We define the solutions $\tilde{\Phi}_{\infty j}$, $j=1,2,3$, of equation~\eqref{phi-tilde} in terms of the solutions $\Phi_{\infty j}$ of the associated integral equations. Each $\Phi_{\infty j}$ is determined by its corresponding initial integration point, and the two sets of functions are related by
\[
\tilde{\Phi}_{\infty j}(x,t,\lambda)
= \Phi_{\infty j}(x,t,\lambda)\, e^{-\ii kp(x,t,\lambda)\sigma_3}.
\]
The functions $\Phi_{\infty j}$ satisfy the integral equation
\begin{equation}\label{inteq_inf}
\Phi_{\infty}(x,t,k)=I+\int_{(x^*,t^*)}^{(x,t)}
	\eul^{\ii k(p(y,\tau,k)-p(x,t,k))\hat\sigma_3}( U_\infty\Phi_\infty \dd y+V_\infty\Phi_\infty \dd \tau)(y,\tau,k),
\end{equation}
where the initial point $(x^*,t^*)$ is chosen as $(0,T)$, $(0,0)$, or $(+\infty,t)$ for $j=1,2,3$, respectively. The integration paths are those shown in Figure~\ref{fig:integration-paths}; due to the compatibility of the Lax pair, the value of the integral is independent of the chosen path.

 More precisely,
\begin{equation}\label{inteq_inf1}
    \begin{aligned}
\Phi_{\infty 1}(x,t,k)=I & +
\int_{0}^{x}
	\eul^{-\ii k\int_y^x q(\xi,t)\dd \xi\hat\sigma_3}( U_\infty\Phi_{\infty1})(y,t,k) \dd y   \\
&-\eul^{-\ii k\int_0^x q(\xi,t)\dd \xi\hat\sigma_3}
\int_{t}^{T}
	\eul^{-\ii k\left(\int_\tau^t \left(\frac{u^2q}{2}\right )(0,s)\dd s-\frac{t-\tau}{4 k^2}\right)\hat\sigma_3}( V_\infty\Phi_{\infty1})(0,\tau,k) \dd \tau,
\end{aligned}
\end{equation}
    
\begin{equation}\label{inteq_inf2}
    \begin{aligned}
\Phi_{\infty 2}(x,t,k)=I & +
\int_{0}^{x}
	\eul^{-\ii k\int_y^x q(\xi,t)\dd \xi\hat\sigma_3}( U_\infty\Phi_{\infty2})(y,t,k) \dd y   \\
&+\eul^{-\ii k\int_0^x q(\xi,t)\dd \xi\hat\sigma_3}
\int_{0}^{t}
	\eul^{-\ii k\left(\int_\tau^t \left(\frac{u^2q}{2}\right )(0,s)\dd s-\frac{t-\tau}{4 k^2}\right)\hat\sigma_3}( V_\infty\Phi_{\infty2})(0,\tau,k) \dd \tau,
\end{aligned}
\end{equation}

    \begin{equation}\label{inteq_inf3}
\Phi_{\infty 3}(x,t,k)=I-
\int_{x}^{+\infty}
	\eul^{\ii k\int_x^y q(\xi,t)\dd \xi\hat\sigma_3}( U_\infty\Phi_{\infty3})(y,t,k) \dd y .
\end{equation}

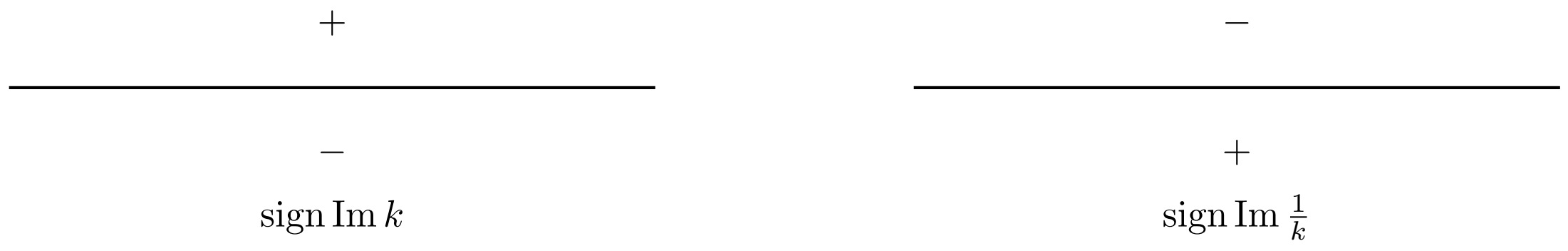
\begin{figure}[ht]
    \centering
\begin{tikzpicture}[scale=1.2]
    % Title

        \node at (0,-1) {$\mathrm{sign}\Im k$};
    
    % Horizontal axis
    \draw[thick] (-2.5,0) -- (2.5,0);
    
    % Signs
    \node at (0,0.5) {$+$};
    \node at (0,-0.5) {$-$};

    % Title
    \node at (7,-1) {$\mathrm{sign}\Im \frac{1}{k}$};
    
    % Horizontal axis
    \draw[thick] (4.5,0) -- (9.5,0);
    
    % Signs
    \node at (7,0.5) {$-$};
    \node at (7,-0.5) {$+$};
\end{tikzpicture}
    \caption{Signature tables for $p(x,t,k)$  }
    \label{fig:sign}
\end{figure}

The domains where the exponentials are bounded in the complex $k$-plane are separated by the real axis; see 
 Figure \ref{fig:sign}.

\begin{remark}
An important  difference from the classical Camassa--Holm equation and other CH-type models (see, e.g. \cites{BS08, KS25}) 
is that the quantity $u^{2}q$ is always nonnegative. Thus, in contrast to those 
equations, no branching based on the sign of $u^{2}q(0,t)$ is needed.
\end{remark}

Let $A^{(1)}$ and $A^{(2)}$ denote the columns of a $2\times 2$ matrix $A = \bigl( A^{(1)}\ \ A^{(2)} \bigr)$. With this notation, using Neumann series expansions, we obtain the following properties of $\Phi_{\infty i}^{(j)}(x,t,k)$:

\begin{enumerate}
    \item $\Phi_{\infty 1}^{(1)}(x,t,k)$ is analytic in $\mathbb{C}\setminus \{0\}$.
 Moreover, $\Phi_{\infty 1}^{(1)}(0,t,k)$=$\begin{pmatrix}
        1\\0
    \end{pmatrix}+O(\frac{1}{k})$ as $k\to\infty$ in $\mathbb{C}^-$; $\Phi_{\infty 1}^{(1)}(x,T,k)$=$\begin{pmatrix}
        1\\0
    \end{pmatrix}+O(\frac{1}{k})$ as $k\to\infty$ in $\mathbb{C}^+$.

    \item $\Phi_{\infty 1}^{(2)}(x,t,k)$ is analytic in $\mathbb{C}\setminus \{0\}$.
 Moreover, $\Phi_{\infty 1}^{(2)}(0,t,k)$=$\begin{pmatrix}
        0\\1
    \end{pmatrix}+O(\frac{1}{k})$ as $k\to\infty$ in $\mathbb{C}^+$; $\Phi_{\infty 1}^{(2)}(x,T,k)$=$\begin{pmatrix}
        0\\1
    \end{pmatrix}+O(\frac{1}{k})$ as $k\to\infty$ in $\mathbb{C}^-$.

     \item $\Phi_{\infty 2}^{(1)}(x,t,k)$ is analytic in $\mathbb{C}\setminus\{0\}$ and $\Phi_{\infty 2}^{(1)}(x,0,k)$ is analytic in $\mathbb{C}$.
 Moreover, $\Phi_{\infty 2}^{(1)}(x,t,k)$=$\begin{pmatrix}
        1\\0
    \end{pmatrix}+O(\frac{1}{k})$ as $k\to\infty$ in $\mathbb{C}^+$.

    \item $\Phi_{\infty 2}^{(2)}(x,t,k)$ is analytic in $\mathbb{C}\setminus\{0\}$ and $\Phi_{\infty 2}^{(2)}(x,0,k)$ is analytic in $\mathbb{C}$.
 Moreover, $\Phi_{\infty 2}^{(2)}(x,t,k)$=$\begin{pmatrix}
        0\\1
    \end{pmatrix}+O(\frac{1}{k})$ as $k\to\infty$ in $\mathbb{C}^-$.

      \item $\Phi_{\infty 3}^{(1)}(x,t,k)$ is analytic in $\mathbb{C}^-$ and continuous up to the boundary.

%\ma{Почему?}
      
 Moreover, $\Phi_{\infty 3}^{(1)}(x,t,k)$=$\begin{pmatrix}
        1\\0
    \end{pmatrix}+O(\frac{1}{k})$ as $k\to\infty$ in $\mathbb{C}^-$.

    \item $\Phi_{\infty 3}^{(2)}(x,t,k)$ is analytic in $\mathbb{C}^+$ and continuous up to the boundary.
 Moreover, $\Phi_{\infty 3}^{(2)}(x,t,k)$=$\begin{pmatrix}
        0\\1
    \end{pmatrix}+O(\frac{1}{k})$ as $k\to\infty$ in $\mathbb{C}^+$.       
\end{enumerate}

Since the coefficient matrices in the Lax pair \eqref{Lax-Q-form} are traceless, we first note that
\begin{equation}\label{jost_det}
  \det \Phi_{\infty j} \equiv 1, \qquad j=1,2,3.  
\end{equation}

Moreover, the matrices $U_\infty(k)\equiv U_\infty(x,t,k)$ and 
$V_\infty(k)\equiv V_\infty(x,t,k)$ satisfy the symmetries
\begin{subequations}\label{sym-UV}
\begin{alignat}{3}\label{sym-U}
U_\infty(- k)&=\overline{U_\infty(\bar k)},&\qquad&
U_\infty(-k)=\sigma_2U_\infty(k)\sigma_2,\\
V_\infty(- k)&=\overline{V_\infty(\bar k)},&&
V_\infty(-k)=\sigma_2V_\infty(k)\sigma_2
\end{alignat}
\end{subequations}
and the phase function $p(k)\equiv p(x,t,k)$ obeys
\begin{equation}\label{sym-p}
p^*(k)=p(k)=p(-k).
\end{equation}
These relations imply that the functions $\Phi_{\infty j}$ inherit the same symmetries:
\begin{equation}\label{sym-Phi}
\Phi_{\infty j}(-k)=\overline{\Phi_{\infty j}(\bar k)},\qquad
\Phi_{\infty j}(-k)=\sigma_2\Phi_{\infty j}(k)\sigma_2.
\end{equation}

The functions $\tilde\Phi_{\infty j}$,
being the solutions of both equations of \eqref{Lax-Q-form}, are related (where they are defined) by matrices independent of $x$ and $t$, and thus
\begin{subequations}\label{rel_inf}
    \begin{align}\label{rel_inf_1}
        &\Phi_{\infty 3}(x,t,k)=\Phi_{\infty 2}(x,t,k)\eul^{-\ii k p(x,t,k)\sigma_3}s(k)\eul^{\ii k p(x,t,k)\sigma_3},\\\label{rel_inf_2}
        &\Phi_{\infty 1}(x,t,k)=\Phi_{\infty 2}(x,t,k)\eul^{-\ii k p(x,t,k)\sigma_3}S(k)\eul^{\ii k p(x,t,k)\sigma_3}
    \end{align}
\end{subequations}
Due to the symmetries of $U_\infty$ and $V_\infty$,  $s(k)$ and $S(k)$
can be written as 
       \begin{equation}
            \label{s}
s(k)=\Phi_{\infty 3}(0,0,k)=\begin{pmatrix}
        a^*(k) & b(k)\\
        -b^*( k) & a(k)
    \end{pmatrix}
    \end{equation}
and
           \begin{equation}
        \label{S}
        S(k)=\Phi_{\infty 1}(0,0,k)=\begin{pmatrix}
        A^*(k) & B(k)\\
        -B^*(k) & A(k)
    \end{pmatrix}\end{equation}
    with some $a(k)$, $b(k)$, $A(k)$, and  $B(k)$.

Relation \eqref{rel_inf_1} can be seen as the scattering relation for the 
$x$-equation \eqref{Lax-Q-form_1}, \eqref{U-hat} (with $t$ being a parameter); similarly, relation \eqref{rel_inf_2} can be viewed 
as the scattering relation for the $t$-equation \eqref{Lax-Q-form_2},
\eqref{hat-V} (with $x$ being a parameter). In what follows,
we will analyze the associated direct and inverse scattering problems.

\subsubsection{Eigenfunctions near $k=0$}

In order to have good control as $k \to 0$, 
it is convenient to rewrite the original Lax pair \eqref{Lax} in the form

\begin{subequations}\label{Lax-2}
\begin{equation}\label{Lax-2-x}
\tilde\Phi_{0x}+\ii k \sigma_3\tilde\Phi_0 =  U_0 \tilde\Phi_0,
\end{equation}
\begin{equation}\label{Lax-2-t}
\tilde\Phi_{0t}+\frac{1}{4 \ii k}\sigma_3\tilde\Phi_0=V_0 \tilde\Phi_0,
\end{equation}
where
\begin{equation}\label{U0-hat}
 U_0(x,t,k):=-\ii k w
\begin{pmatrix}
0 & 1 \\
1 & 0 \\
\end{pmatrix},
\end{equation}
and
\begin{equation}\label{hat-V0}
 V_0(x,t,k):=-\frac{\ii k u^2}{2}\begin{pmatrix}
     1&w\\w&-1
 \end{pmatrix} +\frac{u}{2}\begin{pmatrix}
     0&-1\\1&0
 \end{pmatrix}.
\end{equation}
\end{subequations}
Notice that $U_0(x,t,0)\equiv 0$.

Further, introduce (compare with \eqref{p})
\begin{equation}\label{p_0mu}
p_0(x,t,k):=x-\frac{t}{4 k^2}
\end{equation}
and $Q_0(x,t,k):= \ii k p_0(x,t,k)\sigma_3$.

Next, we introduce the solutions $\tilde{\Phi}_{0j}(x,t,k)$, $j=1,2,3$, of 
\eqref{Lax-2} in the same way as $\tilde{\Phi}_{\infty j}(x,t,k)$ by defining
\[
\Phi_{0}(x,t,k):=\tilde \Phi_{0}(x,t,k)\eul^{\ii k p_0(x,t,k)\sigma_3},
\]
where the functions $\Phi_{0j}(x,t,k)$ are defined as the solutions of the 
integral equations
\begin{equation}\label{inteq_0}
\Phi_{0}(x,t,k)=I+\int_{(x^*,t^*)}^{(x,t)}
	\eul^{\ii k(p_0(y,\tau,k)-p_0(x,t,k))\hat\sigma_3}( U_0\Phi_0 \dd y+V_0\Phi_0 \dd \tau)(y,\tau,k).
\end{equation}
Namely,

\begin{equation}\label{inteq_01}
        \begin{aligned}
\Phi_{0 1}(x,t,k)=I & +
\int_{0}^{x}
	\eul^{-\ii k(x-y)\hat\sigma_3}( U_0\Phi_{01})(y,t,k) \dd y   \\
& - \eul^{-\ii k x\hat\sigma_3}
\int_{t}^{T}
	\eul^{\frac{i}{4 k}(t-\tau)\hat\sigma_3}( V_0\Phi_{01})(0,\tau,k) \dd \tau,
\end{aligned}
\end{equation}

  \begin{equation} \label{inteq_02}
    \begin{aligned}
\Phi_{02}(x,t,k)= I & +
\int_{0}^{x}
	\eul^{-\ii k(x-y)\hat\sigma_3}( U_0\Phi_{02})(y,t,k) \dd y   \\
&+\eul^{-\ii k x\hat\sigma_3}
\int_{0}^{t}
	\eul^{\frac{i}{4 k}(t-\tau)\hat\sigma_3}( V_0\Phi_{02})(0,\tau,k) \dd \tau,
\end{aligned}
\end{equation}

    \begin{equation}\label{inteq_03}
\Phi_{03}(x,t,k)=I-
\int_{x}^{+\infty}
	\eul^{\ii k (y-x)\hat\sigma_3}( U_0\Phi_{03})(y,t,k) \dd y .
\end{equation}

Using Neumann series expansions, we obtain the following properties of $\Phi_{0 i}^{(j)}(x,t,k)$:

\begin{enumerate}
    \item $\Phi_{0 1}^{(1)}(x,t,k)$ is analytic in $\mathbb{C}\setminus \{0\}$.
 Moreover, $\Phi_{01}^{(1)}(x,t,k)=\begin{pmatrix}
        1\\0
    \end{pmatrix}+O(k)$ as $k\to 0$ in $\mathbb{C}^-$.

    \item $\Phi_{0 1}^{(2)}(x,t,k)$ is analytic in $\mathbb{C}\setminus \{0\}$.
 Moreover, $\Phi_{0 1}^{(2)}(x,t,k)=\begin{pmatrix}
        0\\1
    \end{pmatrix}+O(k)$ as $k\to 0$ in $\mathbb{C}^+$.

     \item $\Phi_{0 2}^{(1)}(x,t,k)$ is analytic in $\mathbb{C}\setminus\{0\}$.
 Moreover, $\Phi_{0 2}^{(1)}(x,t,k)=\begin{pmatrix}
        1\\0
    \end{pmatrix}+O(k)$ as $k\to0$ in $\mathbb{C}^+$.

    \item $\Phi_{0 2}^{(2)}(x,t,k)$ is analytic in $\mathbb{C}\setminus\{0\}$.
 Moreover, $\Phi_{0 2}^{(2)}(x,t,k)=\begin{pmatrix}
        0\\1
    \end{pmatrix}+O(k)$ as  $k\to0$ in $\mathbb{C}^-$.

      \item $\Phi_{0 3}^{(1)}(x,t,k)$ is analytic in $\mathbb{C}^-$  and continuous up to the boundary.
 Moreover, $\Phi_{0 3}^{(1)}(x,t,k)$=$\begin{pmatrix}
        1\\0
    \end{pmatrix}+O(k)$ as $k\to 0$ in $\mathbb{C}^-$.

    \item $\Phi_{0 3}^{(2)}(x,t,k)$ is analytic in $\mathbb{C}^+$  and continuous up to the boundary.
 Moreover, $\Phi_{03}^{(2)}(x,t,k)=\begin{pmatrix}
        0\\1
    \end{pmatrix}+O(k)$ as $k\to0$ in $\mathbb{C}^+$.
\end{enumerate}

Since the coefficient matrices in the Lax pair \eqref{Lax-2} are traceless, we note that
\[
\det \Phi_{0 j} \equiv 1, \qquad j=1,2,3.
\]
Furthermore, since $U_0(x,t,k)$, $ V_0(x,t,k)$ satisfy the symmetries \eqref{sym-UV} and $p_0(x,t,k)$ satisfy the symmetries \eqref{sym-p}, it follows that $\Phi_{0 j}(x,t,k)$ satisfy the symmetries \eqref{sym-Phi}. 

The functions $\tilde\Phi_{0 j}$,
being the solutions of both equations of \eqref{Lax-2}, are related (where they are defined) by matrices independent of $x$ and $t$, and thus

\begin{subequations}\label{rel_0}
    \begin{align}\label{rel_0_1}
        &\Phi_{0 3}(x,t,k)=\Phi_{0 2}(x,t,k)\eul^{-\ii kp_0(x,t,k)\sigma_3}\tilde s(k)\eul^{\ii k p_0(x,t,k)\sigma_3},\\\label{rel_0_2}
        &\Phi_{0 1}(x,t,k)=\Phi_{0 2}(x,t,k)\eul^{-\ii k p_0(x,t,k)\sigma_3}\tilde S(k)\eul^{\ii k p_0(x,t,k)\sigma_3}.
    \end{align}
\end{subequations}
Due to the symmetries of $U_0$ and $V_0$,  $\tilde s(k)$ and $\tilde S(k)$
can be written as 
 \begin{equation}
        \label{tils}
    \tilde s(k)=\Phi_{0 3}(0,0,k)=\begin{pmatrix}
        \tilde  a^*(k) & \tilde b(k)\\
       -\tilde  b^*(k) & \tilde a(k)
    \end{pmatrix}
    \end{equation}
and
 \begin{equation}
        \label{tilS}\tilde S(k)=\Phi_{0 1}(0,0,k)=\begin{pmatrix}
        \tilde A^*(k)  &\tilde  B(k)\\
        -\tilde B^*(k) &\tilde  A(k)
    \end{pmatrix}\end{equation}
    with some $\tilde a(k)$, $\tilde b(k)$, $\tilde A(k)$, and  $\tilde B(k)$.

    Similarly to Section \ref{sec:2.1}, relations \eqref{rel_0}
    can be viewed as the scattering relations for 
    respectively equations \eqref{Lax-2-x}, \eqref{U0-hat} and 
    \eqref{Lax-2-t}, \eqref{hat-V0}.

\subsubsection{Relations between eigenfunctions of different Lax pairs}
Since $\Phi_0$ and $\Phi_\infty$ 
come from the 
 same system of ODEs \ref{Lax}, they are related as follows:
\begin{subequations}\label{Phi_0_inf_rel}
\begin{equation}\label{Phi_0_inf_rel_}
 \Phi_{\infty j}(x,t,k)=P(x,t)\Phi_{0 j}(x,t,k)\eul^{-\ii kp_0(x,t,k)\sigma_3}C_j(k)\eul^{\ii k p(x,t,k)\sigma_3}   
\end{equation}
with $C_j(k)$ independent of $x$ and $t$ given by:
\begin{align}\label{Phi_0_inf_coeff}
    C_1&=\eul^{-\frac{\ii T}{4 k}\sigma_3}P^{-1}(0,T)\eul^{-\ii k\left(\int_0^T(\frac{u^2q}{2})(0,\tau)\dd\tau-\frac{T}{4 k^2}\right)\sigma_3},\\\label{Phi_0_inf_coeff2}
    C_2&=P^{-1}(0,0),\\\label{Phi_0_inf_coeff3}
    C_3&=\eul^{-\ii k \nu(0)\sigma_3}.
\end{align}
\end{subequations}
In particular, substituting $(x,t)=(0,0)$ in \ref{Phi_0_inf_rel} yields 
        \begin{align}\label{s_via_til_s}
          &s(k)=P(0,0)\tilde s(k)\eul^{-\ii k \nu(0)\sigma_3},\\\label{S_via_til_S}
          &S(k)=P(0,0)\tilde S(k)C_1=P(0,0)\tilde S(k)\eul^{-\frac{\ii T}{4 k}\sigma_3}P^{-1}(0,T)\eul^{-\ii k\left(\eta(T)-\frac{T}{4 k^2}\right)\sigma_3}
        \end{align}

\begin{remark}\label{rem:simpl}
In the case $u_x(0,0)=0$, we have $P^{-1}(0,0)=I$ 
 and thus \eqref{Phi_0_inf_rel_} for $j=2,3$
reduces to 
\begin{subequations}\label{Phi_0_inf_reduced}
\begin{align}\label{Phi_0_inf_reduced_2}
\Phi_{\infty 2}(x,t,k)& =P(x,t)\Phi_{0 2}(x,t,k)
\eul^{\ii k h_2(x,t))\sigma_3},\\
\Phi_{\infty 3}(x,t,k)& =P(x,t)\Phi_{0 3}(x,t,k)
\eul^{\ii k h_3(x,t))\sigma_3}
\end{align}
\end{subequations}
with 
\[
h_2(x,t) = -\int_x^\infty (q(\xi,t)-1)d\xi + \int_0^\infty
(q(\xi,0)-1)d\xi,
\quad 
h_3(x,t) = -\int_x^\infty (q(\xi,t)-1)d\xi.
\]
Consequently, 
\[
s(k)=\tilde s(k)\eul^{-\ii k \nu(0)\sigma_3},
\]
or 
\[
a(k) = \tilde{a}(k)\eul^{\ii k \nu(0)},\quad
b(k) = \tilde{b}(k)\eul^{\ii k \nu(0)},
\]
which 
results in a significant simplification of the subsequent analysis
and constructions.
\end{remark}

In the general case, in what follows we use the notations
\[
P(0,0)=
\begin{pmatrix}
    \kappa_1^0 & \kappa_2^0 \\
    -\kappa_2^0 & \kappa_1^0
\end{pmatrix},
\]
where $\kappa_1^0$ and $\kappa_2^0$ are expressed in terms of a single quantity 
$u_x(0,0)$. Namely,
\begin{equation}\label{kappas-x}
  \kappa_1^0:=\sqrt{\frac{1+q(0,0)}{2q(0,0)}},
\qquad
\kappa_2^0:=\kappa_1^0\,\frac{u_x(0,0)}{1+q(0,0)},  
\end{equation}
where
\[
q(0,0)=\sqrt{\,1+u_x^2(0,0)\,}.
\]
Using this notation, \eqref{s_via_til_s} reads as
\begin{equation}\label{tilde-a--a}
\begin{aligned}
&a(k)=(-\kappa_2^0\tilde b(k)+\kappa_1^0\tilde a(k))\eul^{\ii k\nu(0)},\\
&b(k)=(\kappa_1^0\tilde b(k)+\kappa_2^0\tilde a(k))\eul^{\ii k\nu(0)}.
\end{aligned}
\end{equation}

Similarly, introducing \[
\kappa_1^T:=\sqrt{\frac{1+q(0,T)}{2q(0,T)}},
\qquad
\kappa_2^T:=\kappa_1^T\,\frac{u_x(0,T)}{1+q(0,T)},
\]
where
\[
q(0,T)=\sqrt{\,1+u_x^2(0,T)\,},
\]
\eqref{S_via_til_S} reads as
\begin{equation}\label{tilde-A--A}
\begin{aligned}
&\tilde A(k)=e^{\,i\left(k\eta(T)-\frac{T}{2k}\right)}
\,\kappa_2^T\big(\kappa_2^0 A^*(k)-\kappa_1^0 B^*(k)\big)
+
e^{-i k\eta(T)}
\,\kappa_1^T\big(\kappa_2^0 B(k)+\kappa_1^0 A(k)\big),\\
&\tilde B(k)=e^{\,i\left(k\eta(T)-\frac{T}{2k}\right)}
\,\kappa_2^T\big(\kappa_1^0 A^*(k)+\kappa_2^0 B^*(k)\big)
+
e^{-i k\eta(T)}
\,\kappa_1^T\big(\kappa_1^0 B(k)-\kappa_2^0 A(k)\big).
\end{aligned}
\end{equation}

\subsubsection{The direct $x$-spectral problem
$\{ w_0(x) \} \mapsto \{ \tilde{a}(k), \tilde{b}(k) \}$}

Starting with some $w_0(x)$ such that $w_0(x)\in L^1(0,\infty)$ (here we don't require  $w_0(x)$ to be the initial value of a solution of the SP equation), we define the direct $x$-spectral mapping
\[
\{ w_0(x) \} \longrightarrow \{ \tilde{a}(k), \tilde{b}(k) \}
\]
by considering the equation \eqref{inteq_03} for $t=0$:
\begin{equation}\label{inteq_03_t_0} 
\Phi_{0 3}(x,0,k)=I-
\int_{x}^{+\infty}
	\eul^{\ii k(y-x)\hat\sigma_3}( U_0\Phi_{03})(y,0,k) \dd y,\qquad x\ge 0
\end{equation}
where $U_0$ is given by \eqref{U0-hat} with $w(x,0)$ replaced by $w_0(x)$.
Then $\tilde{a}(k)$ and $\tilde{b}(k)$ are defined by (cf. \eqref{tils})
\[
\begin{pmatrix}
    \tilde{b}(k) \\ \tilde{a}(k)
\end{pmatrix} = \Phi_{03}^{(2)}(0,0,k).
\]

The analysis of the Volterra integral equation \eqref{inteq_03_t_0} yields  the following properties of $\tilde a(k)$ and $\tilde b(k)$:

\begin{enumerate}
    \item $\tilde a(k)$ and $\tilde b(k)$ are analytic in $\mathbb{C}^+$. Moreover, $\tilde a(k)=1+O(k^2)$ and $\tilde b(k)=\ii k \int_0^\infty w_0(y)\dd y-2k^2\int_0^\infty yw_0(y)\dd y+O(k^3)$ as $k\to 0$  in $\mathbb{C}^+$.

    \item Determinant relation (since the matrices $U_0$ and $V_0$ are traceless):
        \begin{equation}\label{detrel_tilab}
          \tilde a(k)\tilde a^*( k)+ \tilde b(k)\tilde b^*( k)=1, \qquad k\in \mathbb R;
        \end{equation}

    \item Symmetries 
    \begin{equation}\label{sym_tila_}
    \overline{\tilde a(\bar k)}=\tilde a(-k),\qquad
    \overline{\tilde b(\bar k)}=\tilde b(-k).
\end{equation}

\end{enumerate}

\subsubsection{The direct $x$-spectral problem
$\{ w_0(x) \} \mapsto \{ {a}(k), {b}(k) \}$}\label{sub:prop_a}

Alternatively, we can define another direct $x$-spectral 
mapping 
\[
\{ w_0(x) \} \longrightarrow \{ a(k),b(k) \}
\]
considering  equation \eqref{inteq_inf3} for $t=0$:
\begin{equation}\label{inteq_inf3_t_0}    
\Phi_{\infty 3}(x,0,k)=I-
\int_{x}^{+\infty}
	\eul^{\ii k \int_x^yq(\xi,0)\dd \xi\hat\sigma_3}( U_\infty\Phi_{\infty3})(y,0,k) \dd y.
\end{equation}
where $U_\infty$ is given by \eqref{U-hat} with $u_{xx}(x,0)$ replaced by $w_{0x}(x)$ and $q(x,0)$ replaced by $\sqrt{1+w_{0}^2(x)}$.
Then ${a}(k)$ and ${b}(k)$ are defined by (cf. \eqref{s})
\[
\begin{pmatrix}
    {b}(k) \\ {a}(k)
\end{pmatrix} = \Phi_{\infty 3}^{(2)}(0,0,k)
\]
and  have the following properties:

\begin{enumerate}
    \item $a(k)$ and $ b(k)$ are analytic in $\mathbb{C}^+$. Moreover, $ a(k)=1+O(\frac{1}{k})$ and $ b(k)=O(\frac{1}{k})$ as $k\to\infty$ in $\mathbb{C}^+$.

    \item Determinant relation:
        \begin{equation}\label{detrel_ab}
           a(k) a^*( k)+  b(k) b^*( k)=1,  \qquad k\in \mathbb R;
        \end{equation}

    \item Symmetries:

    \begin{equation}
        \label{sym_a_}
    \overline{a(\bar k)}=a(-k),\qquad
    \overline{b(\bar k)}=b(-k).
    \end{equation}

    \item Additionally, it follows from 
   \eqref{s_via_til_s} that as $k\to 0$, $k\in\mathbb{C}^+$:
   \begin{equation}\label{a_at_i}
       a(k)=\kappa_1^0+O(k), \qquad 
    b(k)=\kappa_2^0+O(k),
   \end{equation}
where (cf. \eqref{kappas-x}) 
\begin{equation*}
\kappa_1^0=\sqrt{\frac{1+q_0(0)}{2q_0(0)}},
\qquad
\kappa_2^0=\kappa_1^0\,\frac{w_0(0)}{1+q_0(0)}  
\end{equation*}
with $q_0(0)= \sqrt{\,1+w^2(0)\,}$.

\begin{remark}
    Pushing the expansion of $\eul^{\ii k\nu(0)}$ in \eqref{tilde-a--a} a step further, we obtain
\begin{align*}
&a(k)=\kappa_1^0+\ii k\left(\kappa_2^0 \int_0^\infty w_0(y)\dd y+\kappa_1^0\nu(0)\right)+O(k^2),\quad k\to 0,\quad k\in\mathbb{C}^+\\
&b(k)=\kappa_2^0+\ii k\left(-\kappa_1^0 \int_0^\infty w_0(y)\dd y+\kappa_2^0\nu(0)\right)+O(k^2),\quad k\to 0,\quad k\in\mathbb{C}^+.
\end{align*}
\end{remark}

\end{enumerate}

\subsubsection{The direct $t$-spectral problem $\{v_j(t) \}_0^1 \longrightarrow \{ A(k), B(k) \}$}

%\paragraph{\texorpdfstring{$\{v_j(t) \}_0^1 \longrightarrow \{ A(k), B(k) \}$}{vj(x) -> (A(mu), B(mu))}}

First, consider \eqref{Lax-Q-form} for $x=0$.
Notice that if $u(0,t)\neq0$ then \eqref{SP_1} implies that
\[
u^2(0,t)u_{xx}(0,t)=2\left(u_{xt}(0,t)-u(0,t)-u(0,t)u_x^2(0,t)\right)
\]
and thus $u_{xx}(0,t)$ is determined by $u(0,t)$ and $u_x(0,t)$.

Now, starting with some $v_j(t)$, $j=0,1$, we define the direct $t$-spectral mapping
\[
\{v_j(t) \}_{0}^1 \longrightarrow \{ A(k), B(k) \}
\]
by considering the equation \eqref{inteq_inf1} for $x=0$
\begin{equation}\label{inteq_inf1_x_0}    
\Phi_{\infty 1}(0,t,k)=I-\int_{t}^{T}
	\eul^{-\ii k\left(\int_\tau^t (\frac{u^2q}{2})(0,s)\dd s-\frac{1}{4k^2}(t-\tau)\right)\hat\sigma_3}( V_\infty\Phi_{\infty1})(0,\tau,k) \dd \tau,
\end{equation}
where $V_\infty$ is given by \eqref{hat-V} with $q$ replaced by $\sqrt{1+v_1^2(t)}$ and  $u^2 u_{xx}$ replaced by $2(v_{1t}(t)-v_0(t)-v_0(t)v_1^2(t))$.

Here we do not require $v_j(t)$ to be the boundary values of a solution of the SP equation. However,
if $v_0(t)\equiv 0$ on an interval $[T_1,T_2]$, then we require that $v_{1t}(t)\equiv 0$ on this interval,
 which is consistent with the SP equation with $u\equiv 0$.

 The analysis of the Volterra integral equation \eqref{inteq_inf1_x_0} allows us to define the eigenfunctions and the corresponding spectral functions $A(k)$ and $B(k)$. In particular, $A(k)$ and $B(k)$ satisfy the following properties:

\begin{enumerate}
    \item $ A(k)$ and $ B(k)$ are analytic in $\mathbb{C}\setminus\{0\}$. Moreover, $ A(k)=1+O(\frac{1}{k})+O(\frac{\eul^{2\ii k \eta(T)}}{k})$ and $ B(k)=O(\frac{1}{k})+O(\frac{\eul^{2\ii k \eta(T)}}{k})$ as $k\to\infty$.

    \item Determinant relation:
        \begin{equation}\label{detrel_AB}
           A(k) A^*( k)+  B(k) B^*(k)=1  
        \end{equation}

    \item Symmetries:
    \begin{equation}\label{sym_A_}
    \overline{A(\bar k)}=A(-k),\qquad
    \overline{B(\bar k)}=B(-k).
\end{equation}

\end{enumerate}

\subsubsection{The direct $t$-spectral problem $\{v_j(t) \}_0^1 \longrightarrow \{ \tilde{A}(k), \tilde{B}(k) \}$}\label{sub:prop_A}

Alternatively, we can define another direct $t$-spectral 
mapping 
\[
\{v_j(t) \}_0^1 \longrightarrow \{ \tilde{A}(k), \tilde{B}(k) \}
\]
considering  equation \eqref{inteq_01} for $x=0$:
\begin{equation}\label{inteq_01_x_0}    
\Phi_{0 1}(0,t,k)=I-
\int_{t}^{T}
	\eul^{-\frac{\ii}{4k}(\tau-t)\hat\sigma_3}( V_0\Phi_{01})(0,\tau,k ) \dd \tau,
\end{equation}
where $V_0$ is given by \eqref{hat-V0} with $u$ replaced by $v_0$ and $w$ replaced by $v_1$. Again, if $v_0(t)\equiv 0$ on an interval $[T_1,T_2]$, then we require that $v_{1t}(t)\equiv 0$ on this interval.

Using \eqref{inteq_01_x_0}, one can determine the eigenfunctions and the corresponding spectral functions associated with $\{v_j(t) \}_0^1$.
 In particular, $\tilde A(k)$ and $\tilde B(k)$ have the following properties:

\begin{enumerate}
    \item $\tilde A(k)$ and $\tilde B(k)$ are analytic in $\mathbb{C}\setminus\{0\}$. Moreover, $\tilde A(k)=1+O(k)+O\left(k\eul^{-\frac{\ii T}{2k}}
\right)$ and $\tilde B(k)=O(k)+O\left(k\eul^{-\frac{\ii T}{2k}}\right)$ as $k\to 0$.

    \item Determinant relation:
        \begin{equation}\label{detrel_tilAB}
          \tilde A(k)\tilde A^*(k)+ \tilde B(k)\tilde B^*(k)=1  
        \end{equation}

    \item Symmetries:
    \begin{equation}\label{sym_tilA_}
    \overline{\tilde A(\bar k)}=\tilde A(-k),\qquad 
    \overline{\tilde B(\bar k)}=\tilde B(-k).
\end{equation}

    \item Behavior at $\infty$ (follows from \eqref{tilde-A--A}):

\begin{equation}   \label{tilS_at_infty}
   \tilde A(k)=\kappa_1^T\kappa_1^0\eul^{-\ii k\eta(T)}+O(\frac{\eul^{-\ii k\eta(T)}}{k}),\quad \tilde B(k)=-\kappa_1^T\kappa_2^0\eul^{-\ii k\eta(T)}+O(\frac{\eul^{-\ii k\eta(T)}}{k}),\quad k\to\infty,\quad k\in\mathbb{C}^+,
\end{equation}

\end{enumerate}

\begin{remark}\label{rem_eta}
    Notice that $\eta(T)$, $\kappa_i^0$, and $\kappa_i^T$, $i=1,2$, can be obtained by evaluating $\tilde S(k)$ as $k\to \infty$. Therefore, having $\tilde A(k)$ and $\tilde B(k)$ we can obtain $A(k)$ and $B(k)$ using 
\eqref{S_via_til_S}. 
\end{remark}

\subsection{Problem II} 

As in case of Problem I, assuming that we are given a solution $u(x,t)$ of the SP equation on $x\in(-\infty,0)$, $t\in(0,T)$, with $u(x,t)\to 0$ as $x\to-\infty$ for each fixed $t$, we analyze the analytic properties of the eigenfunctions associated with the Lax pair.  

\begin{remark}
    Notice that in the case of Problem II, \eqref{SP} implies the ``conservation law''
\begin{equation} \label{cons_law_prII}
\nu(t)=\nu(0)-\eta(t),  
\end{equation}
where

\[
\nu(t):=-\int_{-\infty}^0 (q -1)(\xi,t)\dd \xi,\qquad 
    \eta(t):=\frac{1}{2}\int_0^t (u^2q)(0,\tau)\dd \tau.
\]
    
\end{remark}

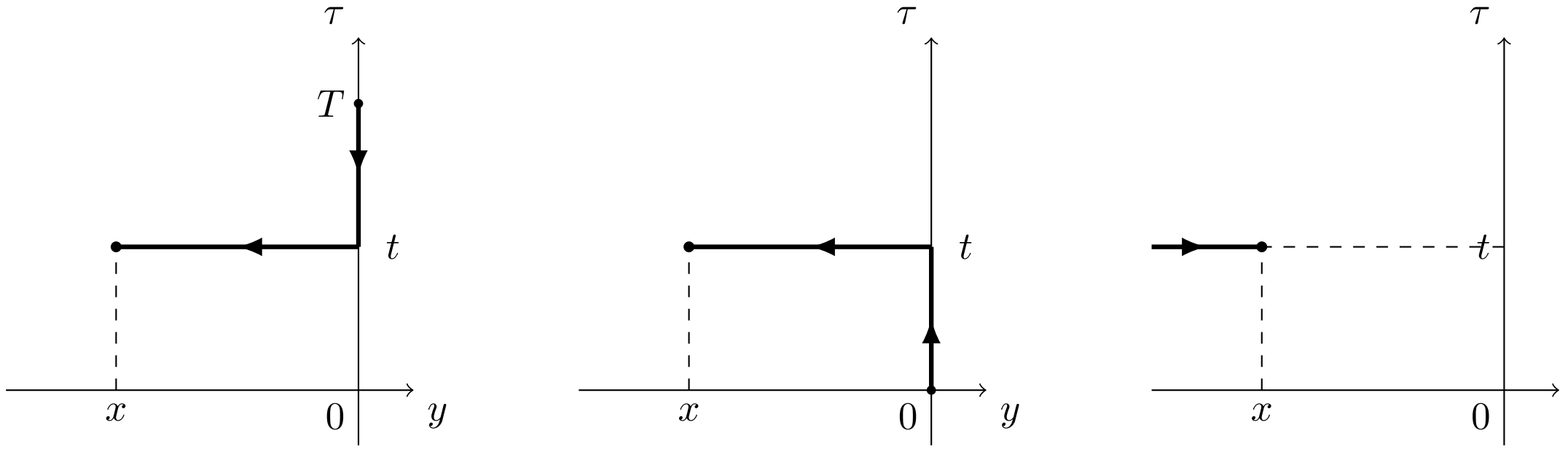
\begin{figure}[h]
\centering
\begin{tikzpicture}[scale=1]
  % Tunables
  \def\x{2.2}  % x-location of the evaluation point
  \def\t{1.3}  % t level
  \def\T{2.6}  % mark on tau-axis
  \def\W{3.2}  % panel width

  % ---------- Panel 1 ----------
  \begin{scope}
    % axes
    \draw[->] (-\W,0) -- (0.5,0) node[below right] {$y$};
    \draw[->] (0,-0.5) -- (0,\W) node[above left] {$\tau$};

    % labels
    \node[below left] at (0,0) {$0$};
    \node[below]      at (-\x,0) {$x$};
    \node[left]       at (0.5,\t) {$t$};

    % T mark
    \fill (0,\T) circle (1.2pt) node[left] {$T$};

    % dashed guide
    \draw[dashed] (-\x,0) -- (-\x,\t);

    % path: down then right (arrow in the middle)
    \draw[very thick,postaction={decorate},
          decoration={markings, mark=at position 0.5 with {\arrow{latex}}}] 
          (0,\T) -- (0,\t);
    \draw[very thick,postaction={decorate},
          decoration={markings, mark=at position 0.5 with {\arrow{latex}}}] 
          (0,\t) -- (-\x,\t);

    % evaluation point
    \fill (-\x,\t) circle (1.4pt);
  \end{scope}

  % ---------- Panel 2 ----------
  \begin{scope}[xshift=5.2cm]
    \draw[->] (-\W,0) -- (0.5,0) node[below right] {$y$};
    \draw[->] (0,-0.5) -- (0,\W) node[above left] {$\tau$};

    \node[below left] at (0,0) {$0$};
    \node[below]      at (-\x,0) {$x$};
    \node[left]       at (0.5,\t) {$t$};

    \fill (0,0) circle (1.2pt);

    \draw[dashed] (-\x,0) -- (-\x,\t);

    % path: up then right (arrow in the middle)
    \draw[very thick,postaction={decorate},
          decoration={markings, mark=at position 0.5 with {\arrow{latex}}}] 
          (0,0) -- (0,\t);
    \draw[very thick,postaction={decorate},
          decoration={markings, mark=at position 0.5 with {\arrow{latex}}}] 
          (0,\t) -- (-\x,\t);

    \fill (-\x,\t) circle (1.4pt);
  \end{scope}

  % ---------- Panel 3 ----------
  \begin{scope}[xshift=10.4cm]
    \draw[->] (-\W,0) -- (0.5,0) node[below right] {$y$};
    \draw[->] (0,-0.5) -- (0,\W) node[above left] {$\tau$};

    \node[below left] at (0,0) {$0$};
    \node[below]      at (-\x,0) {$x$};
    \node[left]       at (0,\t) {$t$};

    \draw[dashed] (-\x,0) -- (-\x,\t);
    \draw[dashed] (0,\t) -- (-\x,\t);

    % path: from right to (x,t), arrow in the middle pointing left
    \draw[very thick,postaction={decorate},
          decoration={markings, mark=at position 0.5 with {\arrow{latex}}}] 
          (-\W,\t) -- (-\x,\t);

    \fill (-\x,\t) circle (1.4pt);
  \end{scope}
\end{tikzpicture}
\caption{Paths of integration for $\Phi_{01}, \Phi_{02},$ and $\Phi_{04}$ 
($\Phi_{\infty1}, \Phi_{\infty2},$ and $\Phi_{\infty 4}$).}
\label{fig:integration-paths_}
\end{figure}

We define the eigenfunctions $\Phi_\infty$ as in \eqref{inteq_inf}. However, instead of $\Phi_{\infty 3}$ we introduce $\Phi_{\infty 4}$, for which the initial point
  $(x^*,t^*)$ is chosen to be $(-\infty,t)$. The integration paths are those shown in Figure \ref{fig:integration-paths_}.

More precisely,
\begin{equation}\label{inteq_inf1_}
    \begin{aligned}
\Phi_{\infty 1}(x,t,\lambda)=I&-
\int_{x}^{0}
	\eul^{-\ii k\int_y^x q(\xi,t)\dd \xi\hat\sigma_3}( U_\infty\Phi_{\infty1})(y,t,k) \dd y   \\
&-\eul^{\ii k\int_x^0 q(\xi,t)\dd \xi\hat\sigma_3}
\int_{t}^{T}
	\eul^{-\ii k\left(\int_\tau^t \left(\frac{u^2q}{2}\right )(0,s)\dd s-\frac{t-\tau}{4 k^2}\right)\hat\sigma_3}( V_\infty\Phi_{\infty1})(0,\tau,k) \dd \tau,
\end{aligned}
\end{equation}
    
\begin{equation}\label{inteq_inf2_}
    \begin{aligned}
\Phi_{\infty 2}(x,t,k)=I&-
\int_{x}^{0}
	\eul^{-\ii k\int_y^x q(\xi,t)\dd \xi\hat\sigma_3}( U_\infty\Phi_{\infty2})(y,t,k) \dd y  \\
&+\eul^{\ii k\int_x^0 q(\xi,t)\dd \xi\hat\sigma_3}
\int_{0}^{t}
	\eul^{-\ii k\left(\int_\tau^t \left(\frac{u^2q}{2}\right )(0,s)\dd s-\frac{t-\tau}{4 k^2}\right)\hat\sigma_3}( V_\infty\Phi_{\infty2})(0,\tau,k) \dd \tau,
\end{aligned}
\end{equation}

    \begin{equation}\label{inteq_inf3_}
\Phi_{\infty 4}(x,t,k)=I+
\int_{-\infty}^{x}
	\eul^{-\ii k\int_y^x q(\xi,t)\dd \xi\hat\sigma_3}( U_\infty\Phi_{\infty4})(y,t,k) \dd y,
\end{equation}
and, therefore, $\Phi^{(j)}_{\infty i}(x,t,k)$ have the following properties:

\begin{enumerate}
    \item $\Phi_{\infty 1}^{(1)}(x,t,k)$ is analytic in $\mathbb{C}\setminus \{0\}$.
 Moreover, $\Phi_{\infty 1}^{(1)}(x,t,k)$=$\begin{pmatrix}
        1\\0
    \end{pmatrix}+O(\frac{1}{k})$ as $k\to\infty$ in $\mathbb{C}^-$.

    \item $\Phi_{\infty 1}^{(2)}(x,t,k)$ is analytic in $\mathbb{C}\setminus \{0\}$.
 Moreover, $\Phi_{\infty 1}^{(2)}(x,t,k)$=$\begin{pmatrix}
        0\\1
    \end{pmatrix}+O(\frac{1}{k})$ as $k\to\infty$ in $\mathbb{C}^+$.

     \item $\Phi_{\infty 2}^{(1)}(x,t,k)$ is analytic in $\mathbb{C}\setminus\{0\}$ and $\Phi_{\infty 2}^{(1)}(x,0,k)$ is analytic in $\mathbb{C}$.
 Moreover, $\Phi_{\infty 2}^{(1)}(0,t,k)$=$\begin{pmatrix}
        1\\0
    \end{pmatrix}+O(\frac{1}{k})$ as $k\to\infty$ in $\mathbb{C}^+$;  $\Phi_{\infty 2}^{(1)}(x,0,k)$=$\begin{pmatrix}
        1\\0
    \end{pmatrix}+O(\frac{1}{k})$ as $k\to\infty$ in $\mathbb{C}^-$.

    \item $\Phi_{\infty 2}^{(2)}(x,t,k)$ is analytic in $\mathbb{C}\setminus\{0\}$ and $\Phi_{\infty 2}^{(2)}(x,0,k)$ is analytic in $\mathbb{C}$.
 Moreover, $\Phi_{\infty 2}^{(2)}(x,0,k)$=$\begin{pmatrix}
        0\\1
    \end{pmatrix}+O(\frac{1}{k})$ as $k\to\infty$ in $\mathbb{C}^-$;  $\Phi_{\infty 2}^{(2)}(0,t,k)$=$\begin{pmatrix}
        0\\1
    \end{pmatrix}+O(\frac{1}{k})$ as $k\to\infty$ in $\mathbb{C}^+$.

      \item $\Phi_{\infty 4}^{(1)}(x,t,k)$ is analytic in $\mathbb{C}^+$ and continuous up to the boundary.
 Moreover, $\Phi_{\infty 4}^{(1)}(x,t,k)$=$\begin{pmatrix}
        1\\0
    \end{pmatrix}+O(\frac{1}{k})$ as $k\to\infty$ in $\mathbb{C}^+$.

    \item $\Phi_{\infty 4}^{(2)}(x,t,k)$ is analytic in $\mathbb{C}^-$ and continuous up to the boundary.
 Moreover, $\Phi_{\infty 4}^{(2)}(x,t,k)$=$\begin{pmatrix}
        0\\1
    \end{pmatrix}+O(\frac{1}{k})$ as $k\to\infty$ in $\mathbb{C}^-$.       
\end{enumerate}

Similarly, we define the eigenfunctions $\Phi_0$ as the solutions of the following integral equations:

\begin{equation}\label{inteq_01_}
        \begin{aligned}
\Phi_{0 1}(x,t,k)=I&-
\int_{x}^{0}
	\eul^{-\ii k(x-y)\hat\sigma_3}( U_0\Phi_{01})(y,t,k) \dd y   \\
&-\eul^{-\ii k x\hat\sigma_3}
\int_{t}^{T}
	\eul^{\frac{i}{4 k}(t-\tau)\hat\sigma_3}( V_0\Phi_{01})(0,\tau,k) \dd \tau,
\end{aligned}
\end{equation}

  \begin{equation} \label{inteq_02_}
    \begin{aligned}
\Phi_{02}(x,t,k)=I&-
\int_{x}^{0}
	\eul^{-\ii k(x-y)\hat\sigma_3}( U_0\Phi_{02})(y,t,k) \dd y   \\
&+\eul^{-\ii k x\hat\sigma_3}
\int_{0}^{t}
	\eul^{\frac{i}{4 k}(t-\tau)\hat\sigma_3}( V_0\Phi_{02})(0,\tau,k) \dd \tau,
\end{aligned}
\end{equation}

    \begin{equation}\label{inteq_03_}
\Phi_{04}(x,t,k)=I+
\int_{-\infty}^{x}
	\eul^{-\ii k (x-y)\hat\sigma_3}( U_0\Phi_{04})(y,t,k) \dd y,
    \end{equation}
and using Neumann series expansions, we obtain the following properties of $\Phi_{0 i}^{(j)}(x,t,k)$:

\begin{enumerate}
    \item $\Phi_{0 1}^{(1)}(x,t,k)$ is analytic in $\mathbb{C}\setminus \{0\}$.
 Moreover, $\Phi_{01}^{(1)}(x,t,k)=\begin{pmatrix}
        1\\0
    \end{pmatrix}+O(k)$ as $k\to 0$ in $\mathbb{C}^-$.

    \item $\Phi_{0 1}^{(2)}(x,t,k)$ is analytic in $\mathbb{C}\setminus \{0\}$.
 Moreover, $\Phi_{0 1}^{(2)}(x,t,k)=\begin{pmatrix}
        0\\1
    \end{pmatrix}+O(k)$ as $k\to 0$ in $\mathbb{C}^+$.

     \item $\Phi_{0 2}^{(1)}(x,t,k)$ is analytic in $\mathbb{C}\setminus\{0\}$.
 Moreover, $\Phi_{0 2}^{(1)}(x,t,k)=\begin{pmatrix}
        1\\0
    \end{pmatrix}+O(k)$ as $k\to0$ in $\mathbb{C}^+$.

    \item $\Phi_{0 2}^{(2)}(x,t,k)$ is analytic in $\mathbb{C}\setminus\{0\}$.
 Moreover, $\Phi_{0 2}^{(2)}(x,t,k)=\begin{pmatrix}
        0\\1
    \end{pmatrix}+O(k)$ as  $k\to0$ in $\mathbb{C}^-$.

      \item $\Phi_{0 4}^{(1)}(x,t,k)$ is analytic in $\mathbb{C}^+$  and continuous up to the boundary.
 Moreover, $\Phi_{0 4}^{(1)}(x,t,k)$=$\begin{pmatrix}
        1\\0
    \end{pmatrix}+O(k)$ as $k\to 0$ in $\mathbb{C}^+$.

    \item $\Phi_{0 4}^{(2)}(x,t,k)$ is analytic in $\mathbb{C}^-$  and continuous up to the boundary.
 Moreover, $\Phi_{04}^{(2)}(x,t,k)=\begin{pmatrix}
        0\\1
    \end{pmatrix}+O(k)$ as $k\to0$ in $\mathbb{C}^-$.
\end{enumerate}
The  properties \eqref{jost_det} and \eqref{sym-Phi} continue to hold for $\Phi_{\infty j}(x,t,k)$ and $\Phi_{0 j}(x,t,k)$ for $j=1,2,4$.

We introduce the matrices $s(k)$, $S(k)$, $\tilde s(k)$, and $\tilde S(k)$ as in \eqref{s}, \eqref{S}, \eqref{tils} and \eqref{tilS}, respectively, except that $\Phi_{\infty 3}$ ($\Phi_{0 3}$) is replaced by  $\Phi_{\infty 4}$ ($\Phi_{0 4}$).

In contrast to the Problem I, the spectral functions $a(k)$,  $b(k)$, $\tilde a(k)$, and $\tilde b(k)$ satisfy the following properties:
\begin{enumerate}[(i)]
    \item $\tilde a(k)$ and $\tilde b(k)$ are analytic in $\mathbb{C}^-$. Moreover, $\tilde a(k)=1+O(k^2)$ and $\tilde b(k)=\ii k \int_0^\infty w_0(y)\dd y+O(k^2)$ as $k\to 0$  in $\mathbb{C}^-$.

        \item $a(k)$ and $ b(k)$ are analytic in $\mathbb{C}^-$. Moreover, $ a(k)=1+O(\frac{1}{k})$ and $ b(k)=O(\frac{1}{k})$ as $k\to\infty$ in $\mathbb{C}^-$. Moreover,
    \begin{align}
    \label{a_at_i_}
    &a(k)=\kappa_1^0+O(k), \quad k\to0,\quad k\in\mathbb{C}^-,\\ \label{b_at_i_}
    &b(k)=\kappa_2^0+O(k), \quad k\to 0,\quad k\in\mathbb{C}^-.
\end{align}
\end{enumerate}
The identities \eqref{detrel_ab}, \eqref{detrel_tilab}, \eqref{sym_a_} and \eqref{sym_tila_} continue to hold.

All properties of $A(k)$, $B(k)$, $\tilde A(k)$, $\tilde B(k)$  remain the same.

\section{The uniqueness result for Problem I}
\label{sec:3}

For later use, it is convenient to introduce the functions
\begin{equation}\label{psi_0}
    \Psi_{0j}(x,t,k)\coloneqq P(x,t)\Phi_{0j}(x,t,k)\eul^{\ii k(p-p_0)(x,t,k)\sigma_3}, 
\end{equation}
which will simplify the form of jump conditions.
Notice that in the case $u_x(0,0)=0$, \eqref{psi_0} simplifies  (see \eqref{Phi_0_inf_reduced}) to
\begin{equation}\label{psi_0_reduced}
\Psi_{02}(x,t,k)=\Phi_{\infty 2}(x,t,k), \qquad 
\Psi_{03}(x,t,k)=\Phi_{\infty 3}(x,t,k)\eul^{\ii k\nu(0)\sigma_3}.
\end{equation}

Furthermore, we will assume that $\epsilon>0$ is sufficiently small so that all 
zeros of $a(k)$ and $a^*(k)$ appearing in the denominators below lie outside the disks $\{|k|\leq\epsilon\}$.

\subsection{The pre-Riemann--Hilbert problem}

Assuming the IBV problem \eqref{SP-2}--\eqref{ic} possesses a solution $ u(x,t)$, using appropriate solutions of the Lax pair equations,
introduce the piecewise meromorphic matrix-valued function $M(x,t,k)$ (depending on $(x,t)$ as parameters)  defined in the domains separated by the contour depicted in Figure \ref{fig:contour_RH_x}:

\begin{equation}\label{M_(xt)}
M^{(xt)}(x,t,k)=\begin{cases}

\left( \frac{\Phi_{\infty 2}^{(1)}(x,t,k)}{a(k)},\Phi_{\infty 3}^{(2)}(x,t,k)\right),\quad k\in\mathbb{C}^+\cap\{|k| >\epsilon\},\\

\left( \frac{\Psi_{0 2}^{(1)}(x,t,k)}{\tilde a(k)},\Psi_{0 3 }^{(2)}(x,t,k)\right)\eul^{-\ii k\nu(0)\sigma_3},\quad k\in\mathbb{C}^+\cap\{|k| <\epsilon\},\\

\left( \Phi_{\infty 3}^{(1)}(x,t,k),\frac{\Phi_{\infty 2}^{(2)}(x,t,k)}{a^*( k)}\right),\quad k\in\mathbb{C}^-\cap\{|k| >\epsilon\},\\

\left( \Psi_{0 3}^{(1)}(x,t,k),\frac{\Psi_{0 2}^{(2)}(x,t,k)}{\tilde a^*( k)}\right)\eul^{-\ii k\nu(0)\sigma_3},\quad k\in\mathbb{C}^-\cap\{|k| <\epsilon\},

\end{cases}
\end{equation}
where the functions $\Psi_{0j}$
are defined by \eqref{psi_0}.

Then the function $M^{(xt)}(x,t,k)$ has the following properties:

\begin{enumerate}
    \item Jump relation across $\mathbb{R}\cup\{|k|=\epsilon\}$
    \begin{subequations}
        \label{jump_M_(xt)}
        \begin{equation}
           M_-^{(xt)}(x,t,k)=M_+^{(xt)}(x,t,k)J^{(xt)}(x,t,k),\quad k\in \mathbb{R}\cup\{|k|=\epsilon\} 
        \end{equation}
        where
         \begin{equation} \label{jump_M_(xt)_matr}
          J^{(xt)}(x,t,k)=\eul^{-\ii kp(x,t,k)\sigma_3} J^{(xt)}_0(k)\eul^{\ii kp(x,t,k)\sigma_3}
        \end{equation} 
        with
\begin{equation}\label{jump_M_(xt)_matr_0}
   J^{(xt)}_0(k)=\begin{cases}
       \begin{pmatrix}
       1&-r^*(k)\\-r(k)&1+r(k)r^*(k)
   \end{pmatrix},\quad k\in\mathbb{R}\cap\{|k| >\epsilon\},\\
   \begin{pmatrix}
       1+\tilde r(k)\tilde r^*(k)&\tilde r^*(k)\eul^{{2}\ii k\nu(0)}\\\tilde r(k)\eul^{-{2}\ii k\nu(0)}&1
   \end{pmatrix},\quad k\in\mathbb{R}\cap\{|k| <\epsilon\},\\
   \begin{pmatrix}
       1&0\\
       -\frac{\kappa_2^0\eul^{-\ii k\nu(0)}}{a(k)\tilde a(k)}&1
   \end{pmatrix},\quad k\in\mathbb{C}^+\cap\{|k| =\epsilon\},\\
   \begin{pmatrix}
       1&-\frac{\kappa_2^0\eul^{\ii k\nu(0)}}{a^*(k)\tilde a^*(k)}\\
       0&1
   \end{pmatrix},\quad k\in\mathbb{C}^-\cap\{|k| =\epsilon\}
   \end{cases} 
\end{equation}
\end{subequations}
and $r(k)=\frac{b^*(k)}{a(k)}$, $\tilde r(k)=\frac{\tilde b^*(k)}{\tilde a(k)}$.  
Moreover,  \eqref{s_via_til_s} implies
\begin{equation}
    \label{r_tilr_con}
    \tilde r(k)=\eul^{2\ii k\nu(0)}
\left(r(k)-\frac{\kappa_2 \eul^{-\ii k\nu(0)}}{a(k)\tilde a(k)}\right), \quad k\in\mathbb{R}.
\end{equation}
\begin{remark}\label{rem:x-RH_reduced}
    In the case 
    $w_0(0)=u_{0x}(0)=0$, the 
    definition of $M^{(xt)}(x,t,k)$ and the corresponding 
    jump conditions simplify,
    in view of \eqref{Phi_0_inf_reduced}, to the following:
     \begin{equation*}
           M_-^{(xt)}(x,t,k)=M_+^{(xt)}(x,t,k)J(x,t,k),\quad k\in \mathbb{R} 
        \end{equation*}
        where
\[
M^{(xt)}(x,t,k)=\begin{cases}
\left( \frac{\Phi_{\infty 2}^{(1)}(x,t,k)}{a(k)},\Phi_{\infty 3}^{(2)}(x,t,k)\right),\quad k\in\mathbb{C}^+,\\
\left( \Phi_{\infty 3}^{(1)}(x,t,k),\frac{\Phi_{\infty 2}^{(2)}(x,t,k)}{a^*( k)}\right),\quad k\in\mathbb{C}^-,
\end{cases}
\]
\begin{equation*}
          J(x,t,k)=\eul^{-\ii kp(x,t,k)\sigma_3} J_0(k)\eul^{\ii kp(x,t,k)\sigma_3}
        \end{equation*} 
        with  
\begin{equation*}
   J_0( k)=
      \begin{pmatrix}
       1&-r^*(k)\\-r(k)&1+r(k)r^*(k)
   \end{pmatrix},\quad k\in\mathbb{R}.
\end{equation*}
Here $\mathbb{R}$ is oriented from left to right.
    
\end{remark}

\item Behavior at $\infty$:
\begin{equation}\label{inf_M_(xt)}
     M^{(xt)}(x,t,k)=I+O(\frac{1}{k}),\quad k\to\infty.
\end{equation}

\item $\det M^{(xt)}(x,t,k)\equiv1$

\item Symmetry properties:
\begin{equation}\label{sym-M_(xt)}
M^{(xt)}( k)=\overline{M^{(xt)}(-\bar k)},\qquad M^{(xt)}(k)=\sigma_2 \overline{M^{(xt)}(\bar k)}\sigma_2
\end{equation}

\item Residue properties. Let $\{k_j\}$ be zeros of $a(k)$ in $\mathbb{C}^+$. Assume that they are all simple. Then

\begin{subequations}\label{res-M_(xt)}
\begin{align}\label{res-M+_(xt)}
\Res_{k_j}M^{(xt)(1)}(x,t,k)&=\frac{e^{2\ii k_jp(x,t,k_j)}}{\dot a(k_j)b(k_j)}M^{(xt)(2)}(x,t,k_j),\\
\label{res-M-_(xt)}
\Res_{\bar k_j}M^{(xt)(2)}(x,t,k)&=-\frac{e^{-2\ii\bar k_jp(x,t,\bar k_j)}}{\dot a^*(\bar k_j)b^*(\bar k_j)}M^{(xt)(1)}(x,t,\bar k_j).
\end{align}
\end{subequations}

Notice that $\dot a^*(\bar k_j)b^*(\bar k_j)=\overline{\dot a( k_j)b( k_j)}$

\item Behavior at $0$:
\begin{equation}\label{i_beh-M_(xt)}
M^{(xt)}(x,t,k)=P(x,t)\left(I-\ii k\begin{pmatrix}
     \tilde \gamma_1(x,t)&u(x,t)\\
     u(x,t)& \tilde \gamma_2(x,t)
\end{pmatrix}+O(k^2)\right),
\end{equation}
where $P(x,t)$ is given by \eqref{P} and $\tilde \gamma_j(x,t)$, $j=1,2$ are some functions, which we do not specify.    
\end{enumerate}

\subsection{The master Riemann--Hilbert problem parametrized by $y$ and $t$}

The construction of jump matrix involves $p(x,t,k)$ (see \eqref{p}), which in turn involves $u(x,t)$, suggesting the introduction of a new variable
\begin{equation}\label{y_}
\begin{aligned}
y(x,t)& =\int_0^xq (\xi,t)\dd\xi+\int_0^t\frac{u^2q}{2}(0,\tau)\dd \tau \\
& = x-\int_x^\infty(q-1)(\xi,t)\dd\xi+\int_0^\infty(q-1)(\xi,0)\dd\xi,
\end{aligned}
\end{equation}
in terms of which both the jump matrix and the residue conditions can be written
explicitly.
Notice that $y(t,x)$ satisfies  the differential equations 
\[
\frac{\partial y}{\partial x}=q(x,t), \qquad \frac{\partial y}  {\partial t}=\frac{u^2q}{2}(x,t)
\]
 (which are compatible due to \eqref{SP}), as well as normalization $y(0,0)=0$. Moreover, $y(0,t)=\eta(t)$.

Then we can formulate a RH problem for $\hat M^{(xt)}(y,t,k):=M^{(xt)}(x(y,t),t,k)$ by replacing $p(x,t,k)$ in \eqref{jump_M_(xt)_matr} and \eqref{res-M_(xt)} by 
\begin{equation}\label{hat_p}
 \hat p(y,t,k):=y-\frac{t}{4k^2}.   
\end{equation}

\textbf{The Riemann--Hilbert problem RH$^{(xt)}$:} Given $a(k)$, $b(k)$ for $k\in\mathbb{C}^+$ (satisfying the properties from Section \ref{sub:prop_a}), and the set $\{k_j\}_1^N\subset \mathbb{C}^+$, find a piece-wise meromorphic $2\times 2$  matrix valued function $\hat M^{(xt)}(y,t,k)$ that satisfies the following conditions:

\begin{enumerate}
    \item Jump condition across $\mathbb{R}\cup|k|=\epsilon$
    \begin{subequations}
        \label{jump_hatM_(xt)}
        \begin{equation}
           \hat M_-^{(xt)}(y, t,k)=\hat M_+^{(xt)}(y,t, k)\hat J^{(xt)}(y,t, k),\quad k\in\mathbb{R} \cup|k|=\epsilon
        \end{equation}
        where
         \begin{equation}\label{jump_hatM_(xt)_matr}
          \hat J^{(xt)}(y,t, k)=\eul^{-\ii k\left(y-\frac{t}{4k^2}\right)\sigma_3} J^{(xt)}_0( k)\eul^{\ii k\left(y-\frac{t}{4k^2}\right)\sigma_3}.
        \end{equation} 
        with $J^{(xt)}_0( k)$ given by \eqref{jump_M_(xt)_matr_0}.
    \end{subequations}

\item Behavior at $\infty$:
\begin{equation}\label{inf_hatM_(xt)}
     \hat M^{(xt)}(y, t,k)=I+O(\frac{1}{ k}),\quad k\to\infty.
\end{equation}

\item Residue conditions:

\begin{align}\label{res-hatM+_(xt)}
\Res_{ k_j}\hat M^{(xt)(1)}(y,t, k)&=\frac{e^{2\ii k_j\left(y-\frac{t}{4k_j^2}\right)}}{\dot a( k_j)b( k_j)}\hat M^{(xt)(2)}(y, t,k_j),\\
\label{res-hatM-_(xt)}
\Res_{\bar k_j}\hat M^{(xt)(2)}(y,t, k)&=-\frac{e^{-2\ii \bar k_j\left(y-\frac{t}{4\bar k_j^2}\right)}}{\dot a^*(\bar k_j)b^*(\bar k_j)}\hat M^{(xt)(1)}(y,t,\bar k_j).
\end{align}   
\end{enumerate}

Importantly, we have the following result concerning the solvability of this problem.

\begin{proposition}\label{prop:exist}
For all $y\ge 0$ and $t\ge 0$, there exists a unique solution $\hat M^{(xt)}(y,t,k)$ of the RH problem \eqref{jump_hatM_(xt)}--\eqref{res-hatM-_(xt)}. Moreover, 
$\det \hat M^{(xt)}\equiv1$ and it satisfies the symmetries
        \begin{equation}\label{sym-hatM_(xt)}
\hat M^{(xt)}( k)=\overline{\hat M^{(xt)}(-\bar k)},\qquad \hat M^{(xt)}(k)=\sigma_2\overline{\hat M^{(xt)}(\bar k)}\sigma_2.
\end{equation}

\end{proposition}

Particularly, the solvability of the Riemann–Hilbert problem \textbf{RH$^{(xt)}$}
is ensured by Zhou’s vanishing lemma (see \cite{Z89}), which is applicable since 
the jump matrix $\hat J^{(xt)}(y,t, k)$ defined in  
\eqref{jump_hatM_(xt)_matr} is positively defined on $\mathbb R$ and satisfies the symmetry \[\hat J^{(xt)}(y, t,k)=\overline{\hat J^{(xt)T}(y,t,\bar  k)}, \qquad k\in\{|k|=\epsilon\}.\]
For details, see Section \ref{sec:7}.

Using the direct spectral analysis above and combining it with Proposition \ref{prop:exist}, we arrive at the following result:

\begin{proposition}\label{uniq_I}
Let $ u(x,t)$ be a solution of the SP equation in the domain $x>0$, $0<t<T$
such that $u(\cdot,t)\in \tilde W^{1,2}((0,\infty))=\{f(x)\in L^1_{loc}((0,\infty)): f'(x),~f''(x) \in L^1((0,\infty))\}$ for all 
$t\in [0,T]$. 
Then $u(x,t)$ can be represented in terms of the unique solution of a Riemann--Hilbert problem \eqref{jump_hatM_(xt)}--\eqref{res-hatM-_(xt)}, for which the data (jump matrix and residue condition) are given in terms of the
initial value $ u(x,0)$, $x\in[0,\infty)$ alone via the associated spectral functions.
In particular, the data $u(x,0)$, $x\in[0,\infty)$
uniquely determine $u(x,t)$, $x\in[0,\infty)$, $t\in[0,T]$
as the solution of the SP equation belonging to 
$\tilde W^{1,2}((0,\infty))$ for all 
$t\in [0,T]$.

\end{proposition}

\begin{proof} The uniqueness of the solution to the RH problem implies, using \eqref{i_beh-M_(xt)}, the following procedure for representing
$ u(x,t)$ in terms of $\hat M^{(xt)}(y,t,k)$:

  \begin{enumerate}[Step 1.]
    \item Having $a(k)$, $b(k)$ compute $\kappa_1^0$, $\kappa_2^0$ and $\nu(0)$, and determine $\tilde a(k)$ and $\tilde b(k)$ via  \eqref{s_via_til_s} (see Remark \ref{rem:data_ab}). With these data, formulate the Riemann–Hilbert problem
 \eqref{jump_hatM_(xt)}--\eqref{res-hatM-_(xt)};
    \item 
    Solve this RH problem for all $y\ge 0$ and $t\in[0,T]$
    and 
evaluate its  solution $\hat M^{(xt)}(y,t,k)$  at $k=0$:
    \begin{equation}\label{hatM_xt__i}
       \hat M^{(xt)}(y,t,k)=\begin{pmatrix}
           \hat \alpha(y,t)&\hat \beta(y,t)\\
           -\hat \beta(y,t)&\hat \alpha(y,t)
       \end{pmatrix}\left( I+\ii k\begin{pmatrix}
           \hat f_1(y,t)&\hat f_2(y,t)\\\hat f_2(y,t)&-\hat f_1(y,t)
       \end{pmatrix}+O(k^2)\right). 
    \end{equation}
    \item 
    Define $\hat u(y,t)$ and $\hat q(y,t)$, $y\ge 0$ and $t\in[0,T]$, from this expansion as follows (cf. \eqref{i_beh-M_(xt)}):
    \begin{subequations}\label{sol_via_RH_xt}
    \begin{align}\label{hat_u_via_RH_xt}
     \hat{ u}(y,t)&=-\hat f_2(y,t), \\\label{hat_q_via_RH_xt}
     \hat q(y,t)&=  \frac{1}{\hat \alpha^2(y,t)-\hat \beta^2(y,t)},\\\label{hat_w_via_RH_xt}
     \hat w(y,t)&=\frac{2\hat \alpha\hat \beta(y,t)}{\hat \alpha^2(y,t)-\hat \beta^2(y,t)}.
    \end{align}
     \end{subequations}
    \item 
    Define $\eta(t)$ as the solution of the differential 
    equation 
    \begin{equation}
        \label{eta_diff}
        \eta'(t)=\frac{1}{2}\hat{ u}^2(\eta,t)\hat{q}(\eta,t), \quad \eta(0)=0 
    \end{equation}
    and 
    define $x(y,t)$ as the solution of the problem (cf. \eqref{y_})
\begin{subequations}\label{y_via_RH_xt}
      \begin{align}
     &\frac{\partial x}{\partial y}=\hat \alpha^2(y,t)-\hat \beta^2(y,t),\\
     &x(\eta(t),t)=0.
    \end{align}
\end{subequations}
\item
Define $u(x,t)$ and $q(x,t)$, $x\ge 0$, $t\in [0,T]$
 by
 \begin{subequations}\label{v_j_via_RH_xt}
\begin{align}
    u(x,t)&=
     \hat{ u}(y(x,t),t),\\
      q(x,t)&=
     \hat{ q}(y(x,t),t).
     \end{align}
    \end{subequations}

\end{enumerate}

\end{proof}

%\textcolor{red}{\begin{remark}
%    The jump matrix and the residue conditions of the master RH problem \textbf{RH$^{(xt)}$} are written in terms of the spectral data ($a(k)$, $b(k)$, $\tilde a(k)$, $\tilde b(k)$ and $\nu(0)$) associated with the initial condition, and do not involve the scattering functions associated with the boundary condition. However, the boundary values  $u(0,t)$ and $u_x(0,t)$ are involved in the change of variable \eqref{y_} through $\eta(t)$.
 %   \end{remark}}

\begin{remark}
     Notice that the Cauchy–Schwarz inequality implies
     \[
     \tilde H^{2,2}((0,\infty))\subset\tilde W^{1,2}((0,\infty)).
     \]
 \end{remark}

\begin{remark} In the case $u_x(0,0)=0$, the master RH problem has the same form as the RH problem associated with the SP equation on the line (cf. \cite{BSL17}).
\end{remark}

\subsection{Well-posedness of Problem I: mapping of RH problems}

The fact that the initial data on a half-line $x\in [l,+\infty)$ determine uniquely the solution of the SP equation in the corresponding quarter-plane $x\ge l$, $t\ge 0$ suggests that the RH problems associated with two initial data, 
$u_1(x,0)$ for $x\in [l_1,+\infty)$ and $u_2(x,0)$ for 
$x\in [l_2,+\infty)$ with $l_2<l_1$ such that 
$u_2(x,0)=u_1(x,0)$ for all $x\in [l_1,+\infty)$,
must be related for $x\in [l_1,+\infty)$, $t\ge 0$.

First, we establish a useful relation amongst the spectral functions of the problems, in the framework of the direct problems. In analogy with \eqref{p} (see also \eqref{p-alt}), 
introduce 
\begin{equation}\label{p-j}
p_j(x,t,k):=\int_{l_j}^{x} q_j(\xi,t)\dd\xi+\frac{1}{2}\int_0^{t} \left(u_j^2q_j\right)(l_j,\tau)\dd\tau-\frac{t}{4 k^2},\qquad j=1,2
\end{equation}
and notice that 
if $u_1(x,t)=u_2(x,t)$ for $x\ge l_1$, then 
\begin{equation}\label{p2-p1}
p_2(x,t,k)-p_1(x,t,k) = p_2(l_1,0,k)-p_1(l_1,0,k) = 
\int_{l_2}^{l_1} q_2(\xi,0)=:\mu>0
\end{equation}
for all $x\ge l_1$, $t\ge 0$.
Introduce the Jost solutions $\Phi_{\infty 2}$, $\Phi_{\infty 3}$ and $\check \Phi_{\infty 2}$
as follows: $\Phi_{\infty 3}$ is introduced by \eqref{inteq_inf3} whereas 
let $\Phi_{\infty 2}$ and $\check \Phi_{\infty 2}$ be the analogues of 
$\Phi_{\infty 2}$ introduced by \eqref{inteq_inf2} taking into account the left ends of the spatial intervals:
\begin{subequations}\label{inteq_inf2-4}
    \begin{align}
\Phi_{\infty 2}(x,t,k)=I & +
\int_{l_1}^{x}
	\eul^{-\ii k\int_y^x q(\xi,t)\dd \xi\hat\sigma_3}( U_\infty\Phi_{\infty2})(y,t,k) \dd y  \nonumber \\
&+\eul^{-\ii k\int_{l_1}^x q(\xi,t)\dd \xi\hat\sigma_3}
\int_{0}^{t}
	\eul^{-\ii k\left(\int_\tau^t \left(\frac{u^2q}{2}\right )(l_1,s)\dd s-\frac{t-\tau}{4 k^2}\right)\hat\sigma_3}( V_\infty\Phi_{\infty2})(l_1,\tau,k) \dd \tau,\\
\check\Phi_{\infty 2}(x,t,k)=I & +
\int_{l_2}^{x}
	\eul^{-\ii k\int_y^x q(\xi,t)\dd \xi\hat\sigma_3}( U_\infty\check\Phi_{\infty2})(y,t,k) \dd y  \nonumber \\
&+\eul^{-\ii k\int_{l_2}^x q(\xi,t)\dd \xi\hat\sigma_3}
\int_{0}^{t}
	\eul^{-\ii k\left(\int_\tau^t \left(\frac{u^2q}{2}\right )(l_2,s)\dd s-\frac{t-\tau}{4 k^2}\right)\hat\sigma_3}( V_\infty\check\Phi_{\infty2})(l_2,\tau,k) \dd \tau.    
\end{align}
\end{subequations}
Here $u(x,t)$ is a solution of SP equation on $x\in[l_2,+\infty)$.

Introduce the scattering matrices $s^{(j)}(k)$, $j=1,2$
(with the corresponding entries $a^{(j)}(k)$ and $b^{(j)}(k)$)
associated with the initial data 
$\{u_j(x,0), x\in [l_j,+\infty)\}$, which are determined by the scattering 
relations similar to \eqref{rel_inf_1}:
\begin{subequations}\label{rel_inf-234}
    \begin{align}\label{rel_inf_32}
        &\Phi_{\infty 3}(x,t,k)=\Phi_{\infty 2}(x,t,k)
        \eul^{-\ii k p_1(x,t,k)\sigma_3}s^{(1)}(k)\eul^{\ii k p_1(x,t,k)\sigma_3},\\\label{rel_inf_34}
        &\Phi_{\infty 3}(x,t,k)=\check\Phi_{\infty 2}(x,t,k)\eul^{-\ii k p_2(x,t,k)\sigma_3}s^{(2)}(k)\eul^{\ii k p_2(x,t,k)\sigma_3},
    \end{align}
\end{subequations}
where 
\begin{subequations}
            \label{s1}
            \begin{equation}
s^{(1)}(k)=\Phi_{\infty 3}(l_1,0,k)=:\begin{pmatrix}
        a^{(1)*}(k) & b^{(1)}(k)\\
        -b^{(1)*}( k) & a^{(1)}(k)
    \end{pmatrix}
    \end{equation}
        and 
    \begin{equation}
            \label{s2}
s^{(2)}(k)=\Phi_{\infty 3}(l_2,0,k)=:\begin{pmatrix}
        a^{(2)*}(k) & b^{(2)}(k)\\
        -b^{(2)*}( k) & a^{(2)}(k)
    \end{pmatrix}.
    \end{equation}
    \end{subequations}
It follows that  $\Phi_{\infty 2}$ and $\check\Phi_{\infty 2}$ can be related 
(taking into account \eqref{p2-p1}) by
\begin{equation}\label{rel-42}
\check\Phi_{\infty 2}(x,t,k) = \Phi_{\infty 2}(x,t,k)
\eul^{-\ii k p_1(x,t,k)\sigma_3}s^{(1)}(k)
\eul^{-\ii k \mu \sigma_3}[s^{(2)}]^{-1}(k)
\eul^{\ii k p_2(x,t,k)\sigma_3}.
\end{equation}

Now we notice that by construction, the first column of 
$\check\Phi_{\infty 2}(x,t,k)$ is analytic in $\mathbb {C}^+$ and, moreover, its (21) entry
is $O(1/k)$ as $k\to \infty$. On the other hand, the (21) entry of the
r.h.s. of \eqref{rel-42} evaluated at $x=l_1$, $t=0$ reads in terms of 
$a^{(j)}$ and $b^{(j)}$ as 
$-b^{(1)*}(k)a^{(2)}(k) + b^{(2)*}(k)a^{(1)}(k)\eul^{2\ii k \mu }$. Consequently, we have 
\begin{equation}\label{gr-ini}
    -b^{(1)*}(k)a^{(2)}(k) + b^{(2)*}(k)a^{(1)}(k)\eul^{2\ii k\mu } = 
    O\left(\frac{1}{k}\right), \qquad k\to \infty, \quad 
    k\in {\mathbb C}^+.
\end{equation}
\begin{remark}\label{rem:gr-ini}
    In general, the spectral functions $b^{(j)*}(k)$ are determined  
    for $k\in {\mathbb C}^-$, but their combination as in the l.h.s.
    of \eqref{gr-ini} turns out to be analytically continued into 
    ${\mathbb C}^+$.
\end{remark}

\begin{proposition}
    Let $M^{(2)}(x,t,k)$, $x\in [l_2, +\infty)$, $t\ge 0$ satisfy the pre-RH problem (see Section 3.1),
    with $p(x,t,k)$ replaced by (cf. \eqref{p-alt}) 
 \begin{equation}\label{p2-alt}
        p_2(x,t,k)=-\int_x^{\infty} (q-1)(\xi,t)\dd\xi+x-l_2 +\int_{l_2}^{\infty} (q-1)(\xi,0)\dd\xi-\frac{t}{4 k^2}
    \end{equation}
    and the spectral functions $a^{(2)}(k)$, $b^{(2)}(k)$
associated with the initial data $u(x,0)$, $x\in [l_2, +\infty)$.
    Define $\tilde p_1(x,t,k):=p_2(x,t,k)-\int_{l_2}^{l_1}q(\xi,0)d\xi 
    \equiv p_2(x,t,k)-\mu$
    and define $\tilde M^{(1)}(x,t,k)$ for $x\in [l_1, +\infty)$, $t\ge 0$ as 
    follows:
    \begin{equation}
     \tilde M^{(1)}(x,t,k) := M^{(2)}(x,t,k) \begin{cases}
         \begin{pmatrix}
             1 &   0 \\
             X^{(1)}(k)\eul^{2\ii k \tilde p_1(x,t,k)} - 
             X^{(2)}(k)\eul^{2\ii k p_2(x,t,k)} & 1
         \end{pmatrix}, & k\in {\mathbb C}^+, \\
         \begin{pmatrix}
             1 &  X^{(2)*}(k)\eul^{-2\ii k p_2(x,t,k)} - 
             X^{(1)*}(k)\eul^{-2\ii k \tilde p_1(x,t,k)} \\
             0 & 1
         \end{pmatrix}, & k\in {\mathbb C}^-, \\
              \end{cases}
    \end{equation}
    where 
    \begin{equation}
       X^{(j)}(k) = \begin{cases}
        r^{(j)}(k), &   k\in {\mathbb C}^+\cap \{|k|>\epsilon\} \\
        r^{(j)}(k) - \frac{\kappa_2^{(0j)}\eul^{-\ii k \nu_j(0)}}{a^{(j)}(k) \tilde a^{(j)}(k)}, & k\in {\mathbb C}^+\cap \{|k|<\epsilon\}
       \end{cases}
    \end{equation}
    with $r^{j}(k)=\frac{b^{(j)*}(k)}{a^{(j)}(k)}$, 
    $\nu_j(0) = \int_{l_j}^\infty(q(\xi,0)-1)d\xi$, and 
    $\kappa_2^{(0j)}=b^{(j)}(0)$.
    
    Then for $x\in [l_1, +\infty)$ and $t\ge 0$,
    $\tilde M^{(1)}(x,t,k) = M^{(1)}(x,t,k)$, where $M^{(1)}(x,t,k)$  satisfies the pre-RH problem associated with the initial data $u(x,0)$ restricted to $x\in [l_1, +\infty)$. 
   
\end{proposition}

\begin{proof}
    In view of the construction of pre-RH problems and the 
    way of extracting the solution of the SP equation from the 
    solution of the associated pre-RH problem, it is sufficient to show that for $k\in {\mathbb C}^+$ and for all $x\ge l_1$, $t\ge 0$,
    \begin{enumerate}[(i)]
        \item 
    $X^{(1)}(k)\eul^{2\ii k \tilde p_1(x,t,k)} - 
             X^{(2)}(k)\eul^{2\ii k p_2(x,t,k)} = O(1/k) $ as $k\to \infty$;
\item      $X^{(1)}(k)\eul^{2\ii k\tilde p_1(x,t,k)} - 
             X^{(2)}(k)\eul^{2\ii k p_2(x,t,k)} = O(k) $ as $k\to 0$.
\end{enumerate}

Regarding item (i), notice that for $k\in {\mathbb C}^+\cap \{|k|>\epsilon\}$,
\[
X^{(1)}(k)\eul^{2\ii k \tilde p_1(x,t,k)} - 
             X^{(2)}(k)\eul^{2\ii k p_2(x,t,k)} = -
             \frac{e^{2\ii k \tilde p_1(x,t,k)}}{a^{(1)}(k)a^{(2)}(k)}
             \left(-b^{(1)*}(k)a^{(2)}(k) + b^{(2)*}(k)a^{(1)}(k)\eul^{2\ii k\mu }\right)
\]
and thus, taking into account  \eqref{gr-ini},
the analyticity of 
$a^{(j)}(k)$ in ${\mathbb C}^+$ with $a^{(j)}(k)=1+O\left(\frac{1}{k}\right)$, and that $\Im\tilde p_1(x,t,k)\ge 0$ for all $x\ge l_1$, $t\ge 0$,
the statement of (i) follows.

Regarding item (ii), we notice that the estimate $X^{(1)}(k)\eul^{2\ii k \tilde p_1(x,t,k)} - 
             X^{(2)}(k)\eul^{2\ii k p_2(x,t,k)}=O(k)$
follow from \eqref{a_at_i}, the estimates $\tilde a^{(j)}(k)
= 1 + O(k^2)$
as $k\to 0$, and the fact that $X^{(1)}(k)\eul^{2\ii k \tilde p_1(x,t,k)} - 
             X^{(2)}(k)\eul^{2\ii k p_2(x,t,k)}$ is analytic in ${\mathbb C}^+$ including near $k=0$ 
(see Remark \ref{rem:gr-ini}).

\end{proof}

\section{The uniqueness result for Problem II}
\label{sec:4}

We will assume that $\epsilon>0$ is sufficiently small so that all 
zeros of $\tilde d(k)$ and $\tilde d^*(k)$ appearing in the denominators below lie outside the disks $\{|k|\leq\epsilon\}$.

\subsection{The pre-Riemann--Hilbert problem}
Introduce 
\begin{align*}
    d(k)=a(k)A^*(k)+b(k)B^*(k), \quad \tilde d(k)=\tilde a(k)\tilde A^*(k)+\tilde b(k)\tilde B^*(k),
\end{align*}
and the piecewise meromorphic matrix-valued function $M(x,t,k)$ (depending on $(x,t)$ as parameters)  defined in the domains separated by the contour depicted in Figure \ref{fig:contour_RH_x}:

\begin{equation}\label{M_(xt)__}
M^{(xt)}(x,t,k)=\begin{cases}

\left(\Phi_{\infty 4}^{(1)}(x,t,k),\frac{\Phi_{\infty 1}^{(2)}(x,t,k)}{d^*(k)}\right),\quad k\in\mathbb{C}^+\cap\{|k| >\epsilon\},\\

\left(\Psi_{0 4}^{(1)}(x,t,k),\frac{\Psi_{0 1}^{(2)}(x,t,k)}{\tilde d^*(k)}\right)\eul^{-\ii k\nu(0)\sigma_3},\quad k\in\mathbb{C}^+\cap\{|k| <\epsilon\},\\

\left( \frac{\Phi_{\infty 1}^{(1)}(x,t,k)}{d(k)},\Phi_{\infty 4}^{(2)}(x,t,k)\right)\eul^{-\ii k\nu(0)\sigma_3},\quad k\in\mathbb{C}^-\cap\{|k| >\epsilon\},\\

\left( \frac{\Psi_{0 1}^{(1)}(x,t,k)}{\tilde d(k)},\Psi_{0 4}^{(2)}(x,t,k)\right),\quad k\in\mathbb{C}^-\cap\{|k| <\epsilon\},

\end{cases}
\end{equation}
where the functions $\Psi_{0j}$
are defined by \eqref{psi_0}.

Then the function $M^{(xt)}(x,t,k)$ satisfy the following properties:
\begin{enumerate}
    \item Jump relation across $\mathbb{R}\cup\{|k|=\epsilon\}$:
    \begin{subequations}
        \label{jump_M_(xt)__}
        \begin{equation}
           M_-^{(xt)}(x,t,k)=M_+^{(xt)}(x,t,k)J^{(xt)}(x,t,k),\quad k\in\Sigma,
        \end{equation}
        where
         \begin{equation}
          J^{(xt)}(x,t,k)=\eul^{-p(x,t,k)\sigma_3} J^{(xt)}_0(k)\eul^{p(x,t,k)\sigma_3}
        \end{equation} 
        with       
\begin{equation}\label{jump_M_(xt)_0__}
   J^{(xt)}_0(k)=\begin{cases}      \begin{pmatrix}
          \frac{1}{dd^*}&-\frac{aB-bA}{d^*}\\
          -\frac{a^*B^*-b^*A^*}{d}&1
       \end{pmatrix},\quad \mu\in\mathbb{R}\setminus\{|k|<\epsilon\},\\

\begin{pmatrix}          1&\frac{\tilde a\tilde B-\tilde b\tilde A}{\tilde d^*}\eul^{2\ii k\nu(0)}\\
          \frac{\tilde a^*\tilde B^*-\tilde b^*\tilde A^*}{\tilde d}\\eul^{-2\ii k\nu(0)}&\frac{1}{\tilde d\tilde d^*}
       \end{pmatrix},\quad \mu\in\mathbb{R}\cap\{|k|<\epsilon\},\\

      \begin{pmatrix}
          1&\frac{\kappa_1(A\tilde B-B\tilde A)+\kappa_2(B \tilde B +A \tilde A)}{\tilde d^* d^*}\eul^{\ii k\nu(0)}\\
          0&1
       \end{pmatrix},\quad \mu\in\mathbb{C}^+\cap\{|k|=\epsilon\},\\

            \begin{pmatrix}
          1&0\\
         \frac{\kappa_1(A^*\tilde B^*-B^*\tilde A^*)+\kappa_2(B^*\tilde B^*+A^* \tilde A^*)}{\tilde d d}\eul^{-\ii k\nu(0)}&1
       \end{pmatrix},\quad \mu\in\mathbb{C}^-\cap\{|k|=\epsilon\}.  
   \end{cases}
\end{equation}
    \end{subequations}

\begin{remark}
    In case $u_{x}(0,0)=0$ the matrix $J_0^{(xt)}(k)$ simplifies to the following:

\begin{equation}\label{jump_M_(xt)_0__simpl}
J^{(xt)}_0(k)=\begin{cases}      \begin{pmatrix}
          \frac{1}{dd^*}&-\frac{aB-bA}{d^*}\\
          -\frac{a^*B^*-b^*A^*}{d}&1
       \end{pmatrix},\quad \mu\in\mathbb{R}\setminus\{|k|<\epsilon\},\\
\begin{pmatrix}
          1&\frac{\tilde a\tilde B-\tilde b\tilde A}{\tilde d^*}\eul^{2\ii k\nu(0)}\\
          \frac{\tilde a^*\tilde B^*-\tilde b^*\tilde A^*}{\tilde d}\eul^{-2\ii k\nu(0)}&\frac{1}{\tilde d\tilde d^*}
       \end{pmatrix},\quad \mu\in\mathbb{R}\cap\{|k|<\epsilon\},\\
      \begin{pmatrix}
         1&\frac{A\tilde B-B\tilde A}{\tilde d^* d^*}\eul^{\ii k\nu(0)}\\
          0&1
       \end{pmatrix},\quad \mu\in\mathbb{C}^+\cap\{|k|=\epsilon\},\\
            \begin{pmatrix}
          1&0\\
         \frac{A^*\tilde B^*-B^*\tilde A^*}{\tilde d d}\eul^{-\ii k\nu(0)}&1
       \end{pmatrix},\quad \mu\in\mathbb{C}^-\cap\{|k|=\epsilon\}. 
   \end{cases}
\end{equation}
    
\end{remark}

\item Behavior at $\infty$:
\begin{equation}\label{inf_M_(xt)__}
     M^{(xt)}(x,t,k)=I+O(\frac{1}{k}),\quad k\to\infty.
\end{equation}

\item $\det M^{(xt)}(x,t,k)\equiv1$

\item Symmetry properties:
\begin{equation}\label{sym-M_(xt)__}
M^{(xt)}(- k)=\overline{M^{(xt)}(\bar k)},\qquad M^{(xt)}(k)=\sigma_2\overline{M^{(xt)}(\bar k)}\sigma_2
\end{equation}

\item Residue properties. Let $\{\eta_j\}$ be zeros of $d(k)$ in $\mathbb{C}^-$. Assume that they are all simple. Then

\begin{subequations}\label{res-M_(xt)__}
\begin{align}\label{res-M+_(xt)__}
\Res_{\eta_j}M^{(xt)(1)}(x,t,k)&=-\frac{e^{2\ii \eta_jp(x,t,\eta_j)}B^*(\eta_j)}{\dot d(\eta_j)a(\eta_j)}M^{(xt)(2)}(x,t,\eta_j),\\
\label{res-M-_(xt)__}
\Res_{\bar \eta_j}M^{(xt)(2)}(x,t,k)&=\frac{e^{-2\ii\bar \eta_jp(x,t,\bar \eta_j)}B(\eta_j)}{\dot d^*(\bar \eta_j)a^*(\bar \eta_j)}M^{(xt)(1)}(x,t,\bar \eta_j).
\end{align}
\end{subequations}

\item Behavior at $0$:
\begin{equation}\label{i_beh-M_(xt)__}
M^{(xt)}(x,t,k)=P(x,t)\left(I-\ii k\begin{pmatrix}
     \tilde \gamma_1(x,t)&u(x,t)\\
     u(x,t)& \tilde \gamma_2(x,t)
\end{pmatrix}+O(k^2)\right) .\end{equation}

\end{enumerate}

\subsection{The master Riemann--Hilbert problem parametrized by $y$ and $t$}

Introducing a new variable
\begin{equation}\label{y___}
y(x,t)=\int_0^x q (\xi,t)\dd\xi+\int_0^t\frac{u^2q}{2}(0,\tau)\dd \tau, 
\end{equation}
in terms of which both the jump matrix and the residue conditions can be written
explicitly, we obtain the master RH problem:

\textbf{The Riemann--Hilbert problem RH$^{(xt)}$:} Given $a(k)$, $b(k)$ for $k\in\mathbb{C}^-$, and $A(k)$, $B(k)$ for $k\in\mathbb{C}$, and the set $\{\eta_j\}_1^N\subset \mathbb{C}^-$, find a piece-wise meromorphic $2\times 2$  matrix valued function $\hat M^{(xt)}(y,t,\mu)$ that satisfies the following conditions:

\begin{enumerate}
    \item Jump condition across $\mathbb{R}\cup|k|=\epsilon$
    \begin{subequations}
        \label{jump_hatM_(xt)__}
        \begin{equation}
           \hat M_-^{(xt)}(y, t,k)=\hat M_+^{(xt)}(y,t, k)\hat J^{(xt)}(y,t, k),\quad k\in\mathbb{R} \cup|k|=\epsilon
        \end{equation}
        where
\begin{equation}\label{jump_hatM_(xt)_matr__}
          \hat J^{(xt)}(y,t, k)=\eul^{-\ii k\left(y-\frac{t}{4k^2}\right)\sigma_3} J^{(xt)}_0( k)\eul^{\ii k\left(y-\frac{t}{4k^2}\right)\sigma_3}.
        \end{equation} 
        with $J^{(xt)}_0( k)$ given by \eqref{jump_M_(xt)_0__}.
    \end{subequations}

\item Behavior at $\infty$:
\begin{equation}\label{inf_hatM_(xt)__}
     \hat M^{(xt)}(y, t,k)=I+O(\frac{1}{ k}),\quad k\to\infty.
\end{equation}

\item Residue conditions:

\begin{align}\label{res-hatM+_(xt)__}
\Res_{ \eta_j}\hat M^{(xt)(1)}(y,t, k)&=-\frac{e^{2\ii \eta_j\left(y-\frac{t}{4\eta_j^2}\right)}B^*(\eta_j)}{\dot d( \eta_j)a( \eta_j)}\hat M^{(xt)(2)}(y, t,\eta_j),\\
\label{res-hatM-_(xt)__}
\Res_{\bar \eta_j}\hat M^{(xt)(2)}(y,t, k)&=\frac{e^{-2\ii \bar \eta_j\left(y-\frac{t}{4\bar \eta_j^2}\right)}B(\bar\eta_j)}{\dot d^*(\bar \eta_j)a^*(\bar \eta_j)}\hat M^{(xt)(1)}(y,t,\bar \eta_j).
\end{align}   
\end{enumerate}

Evaluating the  solution $\hat M^{(xt)}(y,t, k)$ of this RH problem at $k=0$:

    \begin{equation}\label{hatM_xt__i__}
       \hat M^{(xt)}(y,t,k)=\begin{pmatrix}
           \hat \alpha(y,t)&\hat \beta(y,t)\\
           -\hat \beta(y,t)&\hat \alpha(y,t)
       \end{pmatrix}\left( I+\ii k\begin{pmatrix}
           \hat f_1(y,t)&\hat f_2(y,t)\\\hat f_2(y,t)&-\hat f_1(y,t)
       \end{pmatrix}+O(k^2)\right), 
    \end{equation}
define $\hat u(y,t)$ and $\hat q(y,t)$ from this expansion follows (cf. \eqref{i_beh-M_(xt)}):
    \begin{subequations}\label{sol_via_RH_xt__}
    \begin{align}\label{hat_u_via_RH_xt__}
     \hat{ u}(y,t)&=-\hat f_2(y,t), \\\label{hat_q_via_RH_xt__}
     \hat q(y,t)&=  \frac{1}{\hat \alpha^2(y,t)-\hat \beta^2(y,t)},\\\label{hat_w_via_RH_xt__}
     \hat w(y,t)&=\frac{2\hat \alpha\hat \beta(y,t)}{\hat \alpha^2(y,t)-\hat \beta^2(y,t)}.
    \end{align}
    \end{subequations}
    Now, define $x(y,t)$ as the solution of the problem (cf. \eqref{y___})
\begin{subequations}\label{y_via_RH_xt__}
      \begin{align}
     &\frac{\partial x}{\partial y}=\hat \alpha^2(y,t)-\hat \beta^2(y,t),\\
     &x(\eta(t),t)=0.
    \end{align}
\end{subequations}

Then,

\begin{subequations}\label{v_j_via_RH_xt__}
\begin{align}
    u(x,t)&=
     \hat{ u}(y(x,t),t),\\
      q(x,t)&=
     \hat{ q}(y(x,t),t).
     \end{align}
    \end{subequations}

Similarly to Problem I, we get the following result:

%\begin{proposition} Let $u(x,t)$ be a solution of the SP equation in the domain $x<0$, $0<t<T$ \ma{such that $u(\cdot,t)\in \tilde H^{2,2}((-\infty,0))$ for all $t\in [0,T]$}. Then $u(x,t)$ can be represented in terms of the unique solution of a Riemann--Hilbert problem \eqref{jump_hatM_(xt)__}--\eqref{res-hatM-_(xt)__}, for which the data (jump matrix and residue condition) are given \textcolor{red}{in terms of the initial value $ u(x,0)$ and boundary values $u(t,0)$, $u_x(t,0)$} via the associated spectral functions. \ma{In particular, the data $u(x,0)$, $u(0,t)$, $u_x(0,t)$ $x\in(-\infty,0]$ uniquely determine $u(x,t)$, $x\in(-\infty,0]$, $t\in[0,T]$ as the solution of the SP equation belonging to $\tilde H^{2,2}((-\infty,0))$ for all $t\in [0,T]$.} \end{proposition}

\begin{proposition}
Let $u(x,t)$ be a solution of the SP equation in the domain $x<0$, $0<t<T$ such that $u(\cdot,t)\in \tilde W^{1,2}((-\infty,0))=\{f(x)\in L^1_{loc}((-\infty,0)): f'(x),~f''(x) \in L^1((-\infty,0))\}$ for all 
$t\in [0,T]$. Then $u(x,t)$ can be represented in terms of the unique solution of a Riemann--Hilbert problem \eqref{jump_hatM_(xt)__}--\eqref{res-hatM-_(xt)__}, for which the data (jump matrix and residue condition) are given in terms of the
initial value $ u(x,0)$ and boundary values $u(t,0)$, $u_x(t,0)$ via the associated spectral functions. In particular, the data $u(x,0)$, $u(0,t)$, $u_x(0,t)$ $x\in(-\infty,0]$
uniquely determine $u(x,t)$, $x\in(-\infty,0]$, $t\in[0,T]$
as the solution of the SP equation belonging to 
$\tilde W^{1,2}((-\infty,0))$ for all 
$t\in [0,T]$.
 
\end{proposition}

\begin{remark}
     Notice that the Cauchy–Schwarz inequality implies
     \[
     \tilde H^{2,2}((-\infty,0))\subset\tilde W^{1,2}((-\infty,0)).
     \]
 \end{remark}

\section{Global relations}\label{sec:5}

First, we notice that the initial and boundary traces $u(x,0)$, $u(0,t)$,
$u_x(0,t)$, $u_{xx}(0,t), \dots$, involved in the determination of the spectral functions needed in the Riemann-Hilbert formalism for an initial boundary value problem for an integrable PDE,  cannot be prescribed independently for a well-posed IBV problem (see \cites{fokas2008unified,BS08, FIS,BFS06,KS25}). 
However, their compatibility can be characterized in the spectral terms, in the form of the so-called
\emph{Global Relations}.

\subsection{Problem I}

Specifically, in the case of Problem I for the SP equation,  using \eqref{rel_inf} to express $\Phi_{\infty3}(0,T,k)$ in terms of $\Phi_{\infty1}(0,T,k)$ and recalling that $\Phi_{\infty1}(0,T,k)=I$, we obtain
    \begin{equation}\label{for_GR}
        \Phi_{\infty3}(0,T,k)=\eul^{-\ii kp(0,T,k)\sigma_3}S^{-1}(k)s(k)\eul^{\ii kp(0,T,k)\sigma_3}.
    \end{equation}
Now, noticing that $\Phi_{\infty3}(0,T,k)$ can be estimated from
\begin{equation*}
  \Phi_{\infty 3}(0,T,k)=I-
\int_{0}^{+\infty}
	\eul^{\ii k\int_0^y q(\xi,T)\dd \xi\hat\sigma_3}( U_\infty\Phi_{\infty3})(y,T,k) \dd y   
\end{equation*}
we conclude, in particular, that
$ \Phi^{(12)}_{\infty3}(0,T,k)=(\frac{1}{k})$ as $k\to\infty$ in $\mathbb{C}^+$.
Then, considering the element $(12)$ in \eqref{for_GR}, we deduce the global relation 
involving the spectral functions
$a(k)$, $b(k)$, $A(k)$, and $B(k)$:
\begin{equation}\label{relations_Phi_inf3}
(A(k)b(k)-B(k)a(k))\eul^{-2\ii k\int_0^T\frac{u^2q}{2}(0,\tau)\dd\tau}=O(\frac{1}{k}),\quad k\to\infty, \quad k\in\mathbb{C}^+. 
   \end{equation}

Similarly, from the relations  \eqref{rel_0} between $\Phi_{0 j}$, we conclude that
    \begin{equation}\label{for_GR_2}
        \Phi_{03}(0,T,k)=\eul^{-\ii k p_0(0,T,k)\sigma_3}\tilde S^{-1}(k)\tilde s(k)\eul^{\ii kp_0(0,T,k)\sigma_3}.
    \end{equation}
In particular, the $(12)$ entry takes the form
\begin{equation*}
    \Phi^{(12)}_{03}(0,T,k)=(\tilde A(k)\tilde b(k)-\tilde B(k)\tilde a(k))\eul^{\frac{\ii}{2k}T},
\end{equation*}
and, recalling that
 $ \Phi^{(12)}_{03}(0,T,k)=O(k)$ as $k\to 0$, we arrive at the global relation
 involving the spectral functions
$\tilde a(k)$, $\tilde b(k)$, $\tilde A(k)$, and $\tilde B(k)$:
\begin{equation}\label{relations_Phi_03_i}
      (\tilde A(k)\tilde b(k)-\tilde B(k)\tilde a(k))\eul^{\frac{\ii}{2k}T}=O(k),\quad k\to0 ,\quad k\in\mathbb{C}^+.
    \end{equation} 

  We observe that \eqref{relations_Phi_inf3} and 
  \eqref{relations_Phi_03_i} indeed provide restrictions on the behavior of the involved spectral functions, because the exponentials in the left hand sides of \eqref{relations_Phi_inf3} and 
  \eqref{relations_Phi_03_i} are growing at the respective points.
  
\subsection{Problem II} 
Proceeding as in Problem I, we obtain
\begin{equation}\label{relations_Phi_inf3__}
(A(k)b(k)-B(k)a(k))\eul^{-2\ii k\int_0^T\frac{u^2q}{2}(0,\tau)\dd\tau}=O(\frac{1}{k}),\quad k\to\infty, \quad k\in\mathbb{C}^- 
   \end{equation}
   and
\begin{equation}\label{relations_Phi_03_i__}
      (\tilde A(k)\tilde b(k)-\tilde B(k)\tilde a(k))\eul^{\frac{\ii}{2k}T}=O(k),\quad k\to0 ,\quad k\in\mathbb{C}^-.
    \end{equation}

However, in contrast with the case of Problem I, these relations follow from the general properties of $a(k)$, $b(k)$, $\tilde a(k)$, $\tilde b(k)$,  $A(k)$, $B(k)$, $\tilde A(k)$, and $\tilde B(k)$ listed in Section \ref{sec:2}, since the respective exponentials are decaying.

\section{Inverse Spectral Mappings}\label{sec:6}

In this section, we present the inverse spectral mappings, formulated through the solutions of the associated Riemann–Hilbert problems. The jump matrices for these problems are constructed from the corresponding spectral functions.

Furthermore, we will assume that $\epsilon>0$ is sufficiently small so that all 
zeros of $\tilde a(k)$, $\tilde a^*(k)$, $\tilde A(k)$, and $\tilde A^*(k)$ appearing in the denominators below lie outside the disks $\{|\mu\pm\ii|\leq\epsilon\}$.

\subsection{Problem I}\label{subsec:6.1}

\subsubsection{The inverse $x$-spectral mapping}

The inverse $x$-spectral mapping
 \[ \{ a(k),b(k) \}\longrightarrow \{ w_0(x) \}\]
can be formulated in terms of the solution of a Riemann–Hilbert problem whose jump matrix is constructed from the spectral functions $a$ and $b$. The construction of this RH problem follows from the relations
between the functions $\Phi_{\infty2}(x,0,k)$, $\Phi_{\infty3}(x,0,k)$, $\Psi_{02}(x,0,k)$ and $\Psi_{03}(x,0,k)$ obtained in the direct spectral analysis.

\begin{remark}\label{rem:data_ab}
 Notice that $\kappa_1^0$, $\kappa_2^0$, and $\nu(0)$ can be obtained by expanding $a(k)$ and $b(k)$ as $k\to 0$. Therefore, having $a(k)$ and $b(k)$, we can determine $\tilde a(k)$ and $\tilde b(k)$ via \eqref{s_via_til_s}.
\end{remark}

\begin{figure}[h]
    \centering

\begin{tikzpicture}[scale=2]

% --- outlines ---
\draw[red] (0,0) circle (1);
\draw[thick] (-2,0) -- (2,0);

\draw[thick, postaction={decorate},
      decoration={markings, mark=at position 0.5 with {\arrow[scale=2]{>}}}]
      (-2,0) -- (-1,0);

\draw[thick, postaction={decorate},
      decoration={markings, mark=at position 0.7 with {\arrow[scale=2]{<}}}]
      (-1,0) -- (0,0);  

\draw[thick, postaction={decorate},
      decoration={markings, mark=at position 0.5 with {\arrow[scale=2]{>}}}]
      (1,0) -- (2,0);

\draw[thick, red, postaction={decorate},
      decoration={markings, mark=at position 0.25 with {\arrow[scale=2]{<}}}] 
      (0,0) circle (1);

\draw[thick, red, postaction={decorate},
      decoration={markings, mark=at position 0.75 with {\arrow[scale=2]{>}}}] 
      (0,0) circle (1);

% marks
\node at (0,0.1) {\textcolor{red}{$0$}};
\fill[red] (0,0) circle(0.02);

\node at (1.1,0.1) {\textcolor{blue}{$\epsilon$}};
\fill[blue] (1,0) circle(0.02);

\end{tikzpicture}

    \caption{The oriented contour for the Riemann--Hilbert problems \textbf{RH$^{(x)}$}, \textbf{RH$^{(t)}$}, and \textbf{RH$^{(xt)}$}}
    \label{fig:contour_RH_x}
\end{figure}
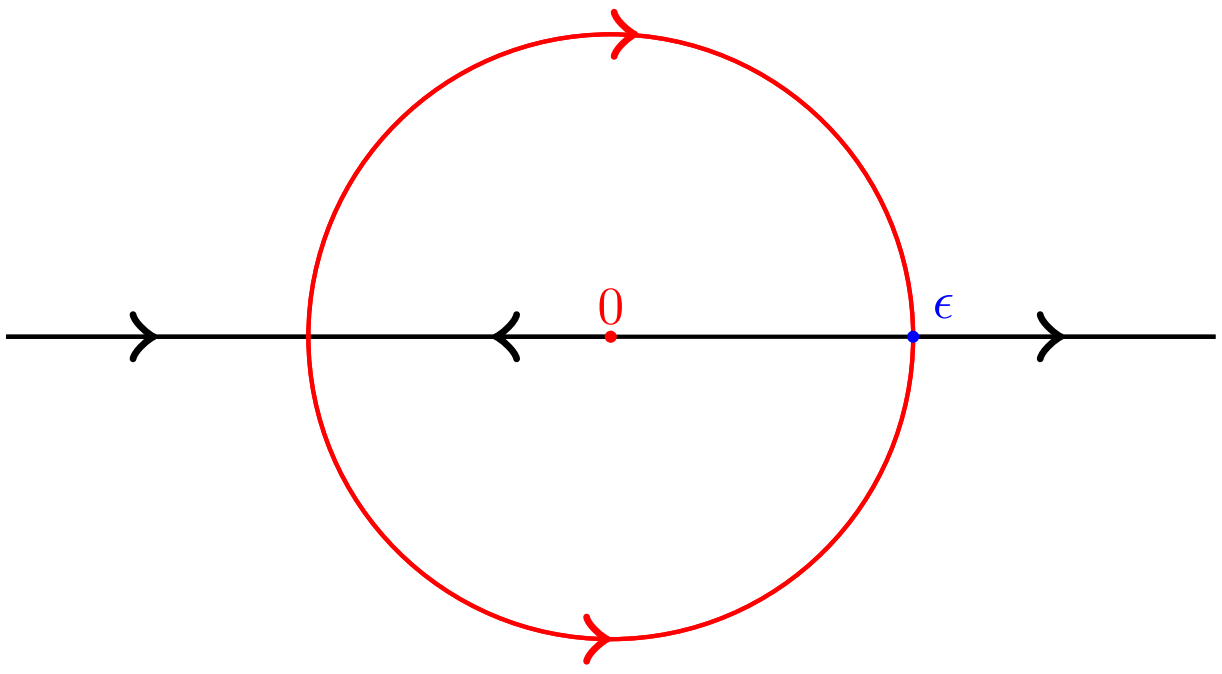

Now, consider the following matrix valued function $M^{(x)}(x,k)$ with $\det M^{(x)}\equiv1$ defined in the domains separated by the contour in Figure \ref{fig:contour_RH_x}:

\begin{equation}\label{M_(x)}
M^{(x)}(x,k)=\begin{cases}

\left( \frac{\Phi_{\infty 2}^{(1)}(x,0,k)}{a(k)},\Phi_{\infty 3}^{(2)}(x,0,k)\right),\quad k\in\mathbb{C}^+\cap\{|k| >\epsilon\},\\

\left( \frac{\Psi_{0 2}^{(1)}(x,0,k)}{\tilde a(k)},\Psi_{0 3}^{(2)}(x,0,k)\right)\eul^{-\ii k\nu(0)\sigma_3},\quad k\in\mathbb{C}^+\cap\{|k| <\epsilon\},\\

\left( \Phi_{\infty 3}^{(1)}(x,0,k),\frac{\Phi_{\infty 2}^{(2)}(x,0,k)}{a^*( k)}\right),\quad k\in\mathbb{C}^-\cap\{|k| >\epsilon\},\\

\left( \Psi_{0 3}^{(1)}(x,0,k),\frac{\Psi_{0 2}^{(2)}(x,0,k)}{\tilde a^*( k)}\right)\eul^{-\ii k\nu(0)\sigma_3},\quad k\in\mathbb{C}^-\cap\{|k| <\epsilon\},
\end{cases}
\end{equation}
where the functions $\Psi_{0j}$
are defined by \eqref{psi_0}.

Here, the functions $\Phi_{\infty2}(x,0,k)$, $\Phi_{\infty3}(x,0,k)$, $\Phi_{02}(x,0,k)$, $\Phi_{03}(x,0,k)$ are to be understood as defined by \eqref{inteq_inf2}, \eqref{inteq_inf3}, \eqref{inteq_02}, and \eqref{inteq_03},respectively,  for $t=0$, where $u(x,0)$ in $U_\infty(x,0)$ and $U_0(x,0)$ is replaced by $u_0(x)$.

\begin{remark}
    Notice that the contour depicted in Figure \ref{fig:contour_RH_x} is symmetric w.r.t.  $ k\mapsto\bar k$ and $ k\mapsto - k$.
\end{remark}

Then the function $M^{(x)}(x,k)$ has the following properties (c.f. \cite{BSL17}):
\begin{enumerate}
    \item Jump relation across $\mathbb{R}\cup\{|k|=\epsilon\}$:
    \begin{subequations}
        \label{jump_M_(x)}
        \begin{equation}
           M_-^{(x)}(x,k)=M_+^{(x)}(x,k)J(x,k),\quad k\in \mathbb{R}\cup\{|k|=\epsilon\} 
        \end{equation}
        where
         \begin{equation}
          J(x,k)=\eul^{-\ii kp(x,0,k)\sigma_3} J_0(k)\eul^{\ii kp(x,0,k)\sigma_3}
        \end{equation} 
        with $p(x,0,k)= \int_0^x q_0(\xi)\dd\xi 
        = \int_0^x \sqrt{1+w_0^2(\xi)}d\xi $ and
\begin{equation}\label{J0}
   J_0(k)=\begin{cases}
       \begin{pmatrix}
       1&-r^*(k)\\-r(k)&1+r(k)r^*(k)
   \end{pmatrix},\quad k\in\mathbb{R}\cap\{|k| >\epsilon\},\\
   \begin{pmatrix}
       1+\tilde r(k)\tilde r^*(k)&\tilde r^*(k)\eul^{2\ii k\nu(0)}\\\tilde r(k)\eul^{-2\ii k\nu(0)}&1
   \end{pmatrix},\quad k\in\mathbb{R}\cap\{|k| <\epsilon\},\\
   \begin{pmatrix}
       1&0\\
       -\frac{\kappa_2^0\eul^{-\ii k\nu(0)}}{a(k)\tilde a(k)}&1
   \end{pmatrix},\quad k\in\mathbb{C}^+\cap\{|k| =\epsilon\},\\
   \begin{pmatrix}
       1&-\frac{\kappa_2^0\eul^{\ii k\nu(0)}}{a^*(k)\tilde a^*(k)}\\
       0&1
   \end{pmatrix},\quad k\in\mathbb{C}^-\cap\{|k| =\epsilon\},
   \end{cases} 
\end{equation}
where $r(k)=\frac{b^*(k)}{a(k)}$ and $\tilde r(k)=\frac{\tilde b^*(k)}{\tilde a(k)}$. 
    \end{subequations}

\item Behavior at $\infty$:
\begin{equation}\label{inf_M_(x)}
     M^{(x)}(x,k)=I+O(\frac{1}{k}),\quad k\to\infty.
\end{equation}

\item $\det M^{(x)}(x,k)\equiv1$

\item Symmetry properties:
\begin{equation}\label{sym-M_(x)}
M^{(x)}(- k)=\overline{M^{(x)}(\bar k)},\qquad M^{(x)}(-k)=\sigma_2M^{(x)}(k)\sigma_2
\end{equation}

\item Residue properties. Let $\{k_j\}$ be zeros of $a(k)$ in $\mathbb{C}^+$. Assume that they are all simple. Then

\begin{align}\label{res-M+_(x)}
\Res_{k_j}M^{(x)(1)}(x,k)&=\frac{e^{2\ii k_jp(x,0,k_j)}}{\dot a(k_j)b(k_j)}M^{(x)(2)}(x,k_j),\\
\label{res-M-_(x)}
\Res_{\bar k_j}M^{(x)(2)}(x,k)&=-\frac{e^{-2\ii\bar k_jp(x,0,\bar k_j)}}{\dot a^*(\bar k_j)b^*(\bar k_j)}M^{(x)(1)}(x,\bar k_j).
\end{align}
\begin{remark}
 Notice that $\dot a^*(\bar k_j)b^*(\bar k_j)=\overline{\dot a( k_j)b( k_j)}$.   
\end{remark}

\item Behavior at $0$:
\begin{equation}\label{i_beh-M_(x)}
M^{(x)}(x,k)=P_0(x)+O(k), \quad k\to 0, 
\end{equation}
where 
\[
P_0(x)=\sqrt{\frac{1+q_0(x)}{2q_0(x)}}\begin{pmatrix}
1 & \frac{w_0(x)}{1+q_0(x)}  \\
-\frac{w_0(x)}{1+q_0(x)} & 1 \\
\end{pmatrix}
\]
with $q_0(x)=\sqrt{1+w_0^2(x)}$.
    
\end{enumerate}

The properties of $M^{(x)}(x,k)$ described above may be viewed as a family of Riemann–Hilbert factorization problems parametrized by $(x,t)$. In the construction of the corresponding jump matrix and residue conditions, however, 
\[
p(x,0, k)= \int_0^x q_0(\xi)\dd\xi
\] 
is involved, which, in turn, involves $q_0(x)$. This motivates the introduction of a new variable
\begin{equation}\label{y_(x)}
y(x)=\int_0^x q(\xi,0)\dd\xi, 
\end{equation}
 which makes the RH problem  explicitly dependent
on $y$ as a parameter:

\textbf{The Riemann--Hilbert problem for $\hat M^{(x)}(y, k)$ (RH$^{(x)}$):} Given $a( k)$, $b( k)$ for $ k\in\mathbb{C}^+$ (satisfying the properties from Section \ref{sub:prop_a}) and the set $\{k_j\}$ in $\mathbb{C}^+$, find a  piece-wise meromorphic $2\times 2$ matrix valued function $\hat M^{(x)}(y, k)$ that satisfies the following conditions:

\begin{enumerate}
    \item Jump condition across $\mathbb{R}\cup\{|k|=\epsilon\}$
    \begin{subequations}
        \label{jump_hatM_(x)}
        \begin{equation}
           \hat M_-^{(x)}(y, k)=\hat M_+^{(x)}(y, k)\hat J(y, k),\quad k\in \mathbb{R}\cup\{|k|=\epsilon\} 
        \end{equation}
        where
         \begin{equation}\label{jump_hatM_(x)_matrix}
          \hat J(y, k)=\eul^{-\ii ky\sigma_3} J_0( k)\eul^{\ii ky\sigma_3}
        \end{equation} 
    \end{subequations}
with $J_0( k)$ given by \eqref{J0}.

\item Behavior at $\infty$:
\begin{equation}\label{inf_hatM_(x)}
     \hat M^{(x)}(y, k)=I+O(\frac{1}{ k}),\quad k\to\infty.
\end{equation}

\item Residue conditions:

\begin{align}\label{res-hatM+_(x)}
\Res_{ k_j}\hat M^{(x)(1)}(y, k)&=\frac{e^{2\ii k_jy}}{\dot a( k_j)b( k_j)}\hat M^{(x)(2)}(y, k_j),\\
\label{res-hatM-_(x)}
\Res_{\bar k_j}\hat M^{(x)(2)}(y, k)&=-\frac{e^{-2\ii\bar  k_jy}}{\dot a^*(\bar k_j)b^*(\bar k_j)}\hat M^{(x)(1)}(y,\bar k_j).
\end{align}   
\end{enumerate}

\begin{proposition}
    \begin{enumerate}
        \item $\det \hat M^{(x)}\equiv1$
    
        \item If the solution of the RH problem \eqref{jump_hatM_(x)}--\eqref{res-hatM-_(x)} exists, it is unique.

        \item If the solution of the RH problem \eqref{jump_hatM_(x)}--\eqref{res-hatM-_(x)} exists, it satisfies the symmetries
        \begin{equation}\label{sym-hatM_(x)}
\hat M^{(x)}( k)=\overline{\hat M^{(x)}(-\bar k)},\qquad \hat M^{(x)}(k)=\sigma_2 \overline{\hat M^{(x)}(\bar k)}\sigma_2
\end{equation}

    \end{enumerate}
\end{proposition}

\begin{proof} \noindent (1) The conditions for $\hat M^{(x)}$ imply that $\det\hat M^{(x)}$ has neither a jump across $\mathbb{R}\cup\{|k|=\epsilon\}$ no singularities at $k_j$. Since $\det\hat M^{(x)}$ tends to $1$ as $k \to\infty$, Liouville’s theorem yields $\det \hat M^{(x)}\equiv1$.

\noindent (2) The uniqueness of the solution follows
using standard arguments based on the application of Liouville's theorem
to the ratio $\hat M_1^{(x)}(\hat M_2^{(x)})^{-1}$ of two 
potential solutions, $\hat M_1^{(x)}$ and $\hat M_2^{(x)}$.

\noindent (3) The symmetries follow from the corresponding symmetries of the jump matrix together with the uniqueness of the solution.

\end{proof}

The uniqueness of the solution to the Riemann–Hilbert problem \textbf{RH$^{(x)}$}, together with the properties of $ M^{(x)}(x,k)$, leads to the following procedure for the inverse mapping
\[
    \{a(k),\, b(k)\} \longrightarrow \{w_0(x)\}
\]
for the $x$-problem:

\begin{enumerate}[Step 1.]
    \item Given $a(k)$ and $b(k)$, compute $\kappa_1^0$, $\kappa_2^0$ and $\nu(0)$, and determine $\tilde a(k)$ and $\tilde b(k)$ via  \eqref{s_via_til_s} (see Remark \ref{rem:data_ab}). With these data, formulate the Riemann–Hilbert problem \textbf{RH$^{(x)}$};

    \item Solve the RH problem constructed in the previous step;

    \item Evaluate the solution $ \hat M^{(x)}(y,k)$ of this RH problem  at $ k=0$:

    \begin{equation}
       \hat M^{(x)}(y, k)=\begin{pmatrix}
           \hat \alpha(y)&\hat \beta(y)\\
           -\hat \beta(y)&\hat \alpha(y)
       \end{pmatrix}+O(k).
    \end{equation}

    \item Then $w_0(x)$ is determined implicitly as follows:
    \begin{subequations}\label{m_0_via_RH_x_par}
        
    \begin{align}\label{m_0_via_RH_x}
    &\hat{w}_0(y)=\frac{2\hat \alpha(y)\hat \beta(y)}{\hat \alpha^2(y)-\hat \beta^2(y)},\\\label{M_(x)_x_y}
        &x(y)=\int_0^y \left(\hat \alpha^2(y)-\hat \beta^2(y)\right) \dd y,\\
       &w_0(x)=\hat{w}_0(y(x)).
    \end{align}
    \end{subequations}
\end{enumerate}

\begin{remark}
    Having $u_x(x,0)$, one can reconstruct $u_0(x)$ by
\begin{equation}\label{Fund_u0}
   u_0(x)=u(0,0)+ \int_0^x u_x(\xi,0)\dd\xi
   \  = -\int_x^\infty u_x(\xi,0)\dd\xi .
\end{equation}
\end{remark}

\begin{remark}\label{rem:jump_coins}
    An important specific feature of the SP equation 
is that the construction of the pre-RH problem uses the same 
data as the RH problem for $M^{(x)}(x,t)$, i.e., the 
scattering functions associated with the initial condition alone.
Particularly, the construction simplifies in the case $u_x(0,0)=0$,
see Remark \ref{rem:x-RH_reduced}.
\end{remark}

\subsubsection{The Inverse $t$-Spectral Mapping}

The inverse $t$-spectral mapping
 \[ \{\tilde A(k),\tilde B(k) \}\longrightarrow \{ v_0(t), v_1(t)  \}\]
 is described in terms of the solution of the RH problem, whose jump matrix is determined by $\tilde A$ and $\tilde B$. The construction of the RH problem follows from the relation
between the eigenfunctions $\Phi_{\infty 1}(0,t,k)$, $\Phi_{\infty 2}(0,t, k)$, $\Psi_{01}(0,t, k)$, and $\Psi_{02}(0,t, k)$ in the direct spectral analysis.

Notice that having $\tilde A( k)$ and $\tilde B( k)$, we can obtain $\eta(T)$, $\kappa_i^0$, and $\kappa_i^T$, $i=1,2$, and determine $A( k)$ and $B(k)$ (see Remark \ref{rem_eta}).

We introduce the matrix valued function $M^{(t)}(t,k)$, with $\det M^{(t)}\equiv1$, defined in
the domains separated by the contour in Figure \eqref{fig:contour_RH_x}:

\begin{equation}\label{M_(t)}
M^{(t)}(t, k)=\begin{cases}

\left( \Phi_{\infty 2}^{(1)}(0,t, k),\frac{\Phi_{\infty 1}^{(2)}(0,t, k)}{A( k)}\right),\quad  k\in\mathbb{C}^+\cap\{| k| >\epsilon\},\\

\left(\Psi_{0 2}^{(1)}(0,t, k),\frac{\Psi_{0 1}^{(2)}(0,t, k)}{\tilde A( k)}\right),\quad  k\in\mathbb{C}^+\cap\{| k| <\epsilon\},\\

\left(\frac{\Phi_{\infty 1}^{(1)}(0,t, k)}{A^*( k)},\Phi_{\infty 2}^{(2)}(0,t, k)\right),\quad  k\in\mathbb{C}^-\cap\{| k| >\epsilon\}\\

\left(\frac{\Psi_{0 1}^{(1)}(0,t, k)}{\tilde A^*( k)},\Psi_{0 2}^{(2)}(0,t, k)\right),\quad  k\in\mathbb{C}^-\cap\{| k| <\epsilon\},
\end{cases}
\end{equation}
where the functions $\Psi_{0j}$
are defined by \eqref{psi_0}.

The function $M^{(t)}(t, k)$ has the following properties:

\begin{enumerate}
    \item Jump relation across $\mathbb{R}\cup\{|k|=\epsilon\}$
    \begin{subequations}
        \label{jump_M_(t)}
        \begin{equation}
           M_-^{(t)}(t, k)=M_+^{(t)}(t, k)J^{(t)}(t, k),\quad k\in\mathbb{R}\cup\{|k|=\epsilon\}
        \end{equation}
        where
         \begin{equation}
          J^{(t)}(t, k)=\eul^{-\ii kp(0,t, k)\sigma_3} J^{(t)}_0( k)\eul^{\ii kp(0,t, k)\sigma_3}
        \end{equation} 
        with $p(0,t, k)= \eta(t)-\frac{t}{4 k^2}$ and
\begin{equation}\label{J_0_jump_M_(t)}
   J^{(t)}_0( k)=\begin{cases}
       \begin{pmatrix}
          1+R( k)R^*( k)&-R( k)\\
          -R^*( k)&1
       \end{pmatrix},\quad  k\in\mathbb{R}\cap\{| k|>\epsilon\},\\
       \begin{pmatrix}
          1&\tilde R( k)\\
        \tilde R^*( k)&1+\tilde R( k)\tilde R^*( k)
       \end{pmatrix},\quad  k\in\mathbb{R}\cap\{| k|<\epsilon\},\\
       
       \begin{pmatrix}
        \kappa_1^0+R(k)\kappa_2^0&\tilde R(k)\kappa_1^0+\tilde R(k)R(k)\kappa_2^0+\kappa_2^0-R(k)\kappa_1^0\\
          -\kappa_2^0 & \kappa_1^0-\tilde R(k)\kappa_2^0
       \end{pmatrix},\quad  k\in\mathbb{C}^+\cap\{| k|=\epsilon\},\\

              \begin{pmatrix}
          \kappa_1^0+ R^*(k)\kappa_2^0&-\kappa_2^0\\
          \tilde R^*(k)\kappa_1^0+\tilde R^*(k)R^*(k)\kappa_2^0+\kappa_2^0-R^*(k)\kappa_1^0&\kappa_1^0-\tilde R^*(k)\kappa_2^0
       \end{pmatrix},\quad  k\in\mathbb{C}^-\cap\{| k|=\epsilon\}
       
   \end{cases}
\end{equation}
and $R( k):=\frac{B( k)}{A( k)}$ and $\tilde R( k):=\frac{\tilde B( k)}{\tilde A( k)}$.
    \end{subequations}

\begin{remark}
    In the  case $P(0,0)=I$ (with $\kappa_1^0=1$ and $\kappa_2^0=0$),
    the jump matrix simplifies to
\begin{equation}
   J^{(t)}_0( k)=\begin{cases}
       \begin{pmatrix}
          1+R( k)R^*( k)&-R( k)\\
          -R^*( k)&1
       \end{pmatrix},\quad  k\in\mathbb{R}\cap\{| k|>\epsilon\},\\
       \begin{pmatrix}
          1&\tilde R( k)\\
        \tilde R^*( k)&1+\tilde R( k)\tilde R^*( k)
       \end{pmatrix},\quad  k\in\mathbb{R}\cap\{| k|<\epsilon\},\\
       
       \begin{pmatrix}
        1&\tilde R(k)-R(k)\\
          0 & 1
       \end{pmatrix},\quad  k\in\mathbb{C}^+\cap\{| k|=\epsilon\},\\

              \begin{pmatrix}
          1&0\\
          \tilde R^*(k)-R^*(k)&1
       \end{pmatrix},\quad  k\in\mathbb{C}^-\cap\{| k|=\epsilon\}
       
   \end{cases}
\end{equation}

\end{remark}

\item Behavior at $\infty$:
\begin{equation}\label{inf_M_(t)}
     M^{(t)}(t, k)=I+O(\frac{1}{ k}),\quad k\to\infty.
\end{equation}

\item $\det M^{(t)}(t, k)\equiv1$

\item Symmetry properties:
\begin{equation}\label{sym-M_(t)}
M^{(t)}( k)=\overline{M^{(t)}(-\bar k)},\qquad M^{(t)}( k)=\sigma_2M^{(t)}( -k)\sigma_2
\end{equation}

\item Residue properties: Let $\{\mu_j\}$ be zeros of $A( k)$ in $\mathbb{C}^+$. Assume that they are all real. Then, using the symmetries \eqref{sym_A_}, it follows that 

\begin{align}\label{res-M+_(t)}
\Res_{\mu_j}M^{(t)(2)}(t, k)&=\frac{B(\mu_j)}{\dot A(\mu_j)}e^{-2\ii\mu_jp(0,t,\mu_j)}M^{(t)(1)}(t,\mu_j),\\
\label{res-M-_(t)}
\Res_{\bar\mu_j}M^{(t)(1)}(t, k)&=-\frac{B^*(\bar\mu_j)}{\dot A^*(\bar \mu_j)}e^{2\ii\bar\mu_jp(0,t,\bar\mu_j)}M^{(t)(2)}(t,-\mu_j).
\end{align}

\item Behavior at $0$. Expanding $\Phi_{02}^{(1)}(0,t, k)$ and $\Phi_{01}^{(2)}(0,t, k)$ at $ k=0$ via the Neumann series and using \eqref{psi_0} yields
\begin{equation}\label{i_beh-M_(t)}
M^{(t)}(t, k)=P^{x=0}(t)\left(I-\ii k\begin{pmatrix}
     \tilde \gamma_1(t)&v_0(t)\\
     v_0(t)& \tilde \gamma_2(t)
\end{pmatrix}+O(k^2)\right),
\end{equation}
where
\[
P^{x=0}(t)=\sqrt{\frac{1+q^0(t)}{2q^0(t)}}\begin{pmatrix}
1 & \frac{v_1(t)}{1+q^0(t)}  \\
-\frac{v_1(t)}{1+q^0(t)} & 1 \\
\end{pmatrix}
\]
with $q^0(t)=\sqrt{1+v_1^2(t)}$.

\end{enumerate}

The properties of $M^{(t)}(t, k)$ outlined above can be interpreted as a family of Riemann–Hilbert factorization problems parametrized by $(x,t)$, yet the construction of the jump matrix and residue conditions involves $p(0,t, k)$, which in turn involves 
\[
\eta(t)=\frac{1}{2}\int_0^t v_0^2(\tau)\sqrt{1+v_1^2(\tau)}\dd\tau.
\]
This suggests the introduction of a new parameter for the RH problem,
$z$, in terms of which $\hat J^{(t)}$ can be explicitly written: \[\hat p(z,t,k)=
z-\frac{t}{4k^2},\] so that,
substituting 
\begin{equation}\label{y_(t)}
z = z(t)=\eta(t),
\end{equation}
we have
$p(0,t,k)=\hat p(\eta(t),t,k)$.

Notice that $\frac{\dd \eta(t)}{\dd t}=\frac{1}{2}v_0^2(t)\sqrt{1+v_1^2(t)}$.

Consider the following family of RH problems parametrized by $z$ and $t$:

\textbf{The Riemann--Hilbert problem for $\hat M^{(t)}(z,t, k)$ (RH$^{(t)}$):} Given $\tilde A( k)$, $\tilde B( k)$ for $ k\in\mathbb{C}\setminus\{0\}$ (satisfying the properties from Section \ref{sub:prop_A}), and the set $\{\mu_j\}$ in $\mathbb{C}^+$,
find  a piece-wise meromorphic $2\times 2$ matrix valued function $\hat M^{(t)}(z,t, k)$ that satisfies the following conditions:

\begin{enumerate}
    \item Jump relation across $\mathbb{R}\cup\{|k|=\epsilon\}$
    \begin{subequations}
        \label{jump_hatM_(t)}
        \begin{equation}
           \hat M_-^{(t)}(z,t, k)=\hat M_+^{(t)}(z,t, k)\hat J^{(t)}(z,t, k),\quad k\in\mathbb{R}\cup\{|k|=\epsilon\}
        \end{equation}
        where
         \begin{equation}
          \hat J^{(t)}(z,t, k)=\eul^{-\ii k\left(z-\frac{t}{4 k^2}\right)\sigma_3}J_0^{(t)}( k)\eul^{\ii k\left(z-\frac{t}{4 k^2}\right)\sigma_3}
        \end{equation} 
    \end{subequations}
    and $J_0^{(t)}( k)$ is defined in\eqref{J_0_jump_M_(t)}.

\item Behavior at $\infty$:
\begin{equation}\label{inf_hatM_(t)}
     \hat M^{(t)}(z,t, k)=I+O(\frac{1}{ k}),\quad k\to\infty.
\end{equation}

\item Residue conditions: 

\begin{align}\label{res-hatM+_(t)}
\Res_{\mu_j} \hat M^{(t)(2)}(t, k)&=\frac{B(\mu_j)}{\dot A(\mu_j)}e^{-2\ii\mu_j\left(z-\frac{t}{4 \mu_j^2}\right)} \hat M^{(t)(1)}(t,\mu_j),\\
\label{res-hatM-_(t)}
\Res_{-\mu_j} \hat M^{(t)(1)}(t, k)&=-\frac{B(\mu_j)}{\dot A(\mu_j)}e^{-2\ii\mu_j\left(z-\frac{t}{4 \mu_j^2}\right)} \hat M^{(t)(2)}(t,-\mu_j).
\end{align}
\end{enumerate}

Similarly to the $x$-problem, the following proposition holds:

\begin{proposition}
    \begin{enumerate}
        \item $\det \hat M^{(t)}\equiv1$
    
        \item If the solution of the RH problem \eqref{jump_hatM_(t)}--\eqref{res-hatM-_(t)} exists, it is unique.

        \item If the solution of the RH problem \eqref{jump_hatM_(t)}--\eqref{res-hatM-_(t)} exists, it satisfies the symmetries
        \begin{equation}\label{sym-hatM_(t)}
\hat M^{(t)}( k)=\overline{\hat M^{(t)}(-\bar k)},\qquad \hat M^{(t)}(k)=\sigma_2\overline{\hat M^{(t)}(\bar k)}\sigma_2
\end{equation}
 
    \end{enumerate}
\end{proposition}

The uniqueness of the solution of the Riemann--Hilbert problem \textbf{RH$^{(t)}$} justifies the following procedure for the inverse mapping
\[
    \{\tilde A(k),\,\tilde B(k) \} \longrightarrow \{v_0(t), v_1(t)\}
\]
for the $t$-problem:

\begin{enumerate}[Step 1.]
    \item Given $\tilde A(k)$ and $\tilde B(k)$, compute $\kappa_i^0$, $\kappa_i^T$, $i=1,2$, and $\eta(T)$, and determine $A(k)$ and $B(k)$ via  \eqref{S_via_til_S} (see Remark \ref{rem_eta}). Having these data, formulate the Riemann–Hilbert problem \textbf{RH$^{(t)}$};

    \item Solve the RH problem constructed in the previous step;

    \item Evaluate the solution $\hat M^{(t)}(z,t, k)$ of this RH problem at $ k=0$:

    \begin{equation}
       \hat M^{(t)}(z,t, k)=\begin{pmatrix}
           \hat \alpha(z,t)&\hat \beta(z,t)\\
           -\hat \beta(z,t)&\hat \alpha(z,t)
       \end{pmatrix}\left( I+\ii k\begin{pmatrix}
           \hat f_1(z,t)&\hat f_2(z,t)\\\hat f_2(z,t)&-\hat f_1(z,t)
       \end{pmatrix}+O(k^2)\right). 
    \end{equation}

    \item Define $\hat v_j(z,t)$, $j=0,1$ from this expansions in the following way (cf. \eqref{i_beh-M_(t)}):
    \begin{align}\label{hat_v_j_via_RH_t}
     \hat v_0(z,t)&=-\hat f_2(z,t), \\
     \hat v_1(z,t)&=\frac{2\hat\beta(z,t)\hat\alpha(z,t)}{\hat \alpha^2(z,t)-\hat \beta^2(z,t)}.
    \end{align}
    Notice that $\hat q^0(z,t)= \frac{1}{\hat \alpha^2(z,t)-\hat \beta^2(z,t)}.$

    Now, define $z(t)$ as a solution of the following differential equation (cf. \eqref{y_(t)}):

    \begin{align}\label{y_(t)_via_RH_t}
     \frac{\dd z}{\dd t}&=\frac{1}{2}\hat f_2^2(z,t)\frac{1}{\hat \alpha^2(z,t)-\hat \beta^2(z,t)},\\
     z(0)&=0.
    \end{align}

Then,

\begin{equation}\label{v_j_via_RH_t}
     v_j(t)=
     \hat v_j(z(t),t), \quad j=0,1.
    \end{equation}
\end{enumerate}

\subsection{Problem II}\label{subsec:6.2}

\subsubsection{The Inverse $x$-Spectral Mapping}
Recalling the analytic properties of eigenfunctions and scattering coefficients, we introduce the matrix-valued function  $M^{(x)}(x,k)$ defined in the domains separated by the contour in Figure \ref{fig:contour_RH_x}:

\begin{equation}\label{M_(x)_}
M^{(x)}(x,k)=\begin{cases}

\left(\Phi_{\infty 4}^{(1)}(x,0,k), \frac{\Phi_{\infty 2}^{(2)}(x,0,k)}{a^*(k)}\right),\quad k\in\mathbb{C}^+\cap\{| k| >\epsilon\},\\

\left(\Psi_{0 4}^{(1)}(x,0,k), \frac{\Psi_{0 2}^{(2)}(x,0,k)}{\tilde a^*(k)}\right)\eul^{-\ii k\nu(0)\sigma_3},\quad k\in\mathbb{C}^+\cap\{| k| <\epsilon\},\\

\left( \frac{\Psi_{0 2}^{(1)}(x,0,k)}{\tilde a( k)},\Psi_{0 4}^{(2)}(x,0,k)\right)\eul^{-\ii k\nu(0)\sigma_3},\quad k\in\mathbb{C}^-\cap\{| k| <\epsilon\},\\

\left( \frac{\Phi_{\infty 2}^{(1)}(x,0,k)}{a( k)},\Phi_{\infty 4}^{(2)}(x,0,k)\right),\quad k\in\mathbb{C}^-\cap\{| k| >\epsilon\}.

\end{cases}
\end{equation}

Then the function $M^{(x)}(x,k)$ has the following properties:
\begin{enumerate}
    \item Jump relation across $\mathbb{R}\cup\{|k|=\epsilon\}$:
    \begin{subequations}
        \label{jump_M_(x)_left}
        \begin{equation}
           M_-^{(x)}(x,k)=M_+^{(x)}(x,k)J(x,k),\quad k\in \mathbb{R}\cup\{|k|=\epsilon\} 
        \end{equation}
        where
         \begin{equation}
          J(x,k)=\eul^{-\ii kp(x,0,k)\sigma_3} J_0(k)\eul^{\ii kp(x,0,k)\sigma_3}
        \end{equation} 
        with $p(x,0,k)= \int_0^x q_0(\xi)\dd\xi 
        = \int_0^x \sqrt{1+w_0^2(\xi)}d\xi $ and
\begin{equation}\label{J_0__}
   J_0(k)=\begin{cases}
       \begin{pmatrix}
       1+r(k)r^*(k)&r^*(k)\\r(k)&1
   \end{pmatrix},\quad k\in\mathbb{R}\cap\{|k| >\epsilon\},\\
   \begin{pmatrix}
     1 &-\tilde r^*(k)\eul^{2\ii k\nu(0)}\\-\tilde r(k)\eul^{-2\ii k\nu(0)}& 1+\tilde r(k)\tilde r^*(k)
   \end{pmatrix},\quad k\in\mathbb{R}\cap\{|k| <\epsilon\},\\
\begin{pmatrix}
       1&\frac{\kappa_2^0 \eul^{\ii k\nu(0)}}{a^*(k)\tilde a^*(k)}\\
       0&1
   \end{pmatrix}
,\quad k\in\mathbb{C}^+\cap\{|k| =\epsilon\},\\
    \begin{pmatrix}
       1&0\\
    \frac{\kappa_2^0 \eul^{-\ii k\nu(0)}}{a(k)\tilde a(k)}&1
   \end{pmatrix},\quad k\in\mathbb{C}^-\cap\{|k| =\epsilon\},
   \end{cases} 
\end{equation}
where $r(k)=\frac{b^*(k)}{a(k)}$ and $\tilde r(k)=\frac{\tilde b^*(k)}{\tilde a(k)}$. 
    \end{subequations}

\item Behavior at $\infty$:
\begin{equation}\label{inf_M_(x)__}
     M^{(x)}(x,k)=I+O(\frac{1}{k}),\quad k\to\infty.
\end{equation}

\item $\det M^{(x)}(x,k)\equiv1$

\item Symmetry properties:
\begin{equation}\label{sym-M_(x)__}
M^{(x)}(- k)=\overline{M^{(x)}(\bar k)},\qquad M^{(x)}(-k)=\sigma_2M^{(x)}(k)\sigma_2
\end{equation}

\item Residue properties. Let $\{k_j\}$ be zeros of $a(k)$ in $\mathbb{C}^-$. Assume that they are all simple. Then

\begin{align}\label{res-M+_(x)__}
\Res_{k_j}M^{(x)(1)}(x,k)&=\frac{e^{2\ii k_jp(x,0,k_j)}}{\dot a(k_j)b(k_j)}M^{(x)(2)}(x,k_j),\\
\label{res-M-_(x)_}
\Res_{\bar k_j}M^{(x)(2)}(x,k)&=-\frac{e^{-2\ii\bar k_jp(x,0,\bar k_j)}}{\dot a^*(\bar k_j)b^*(\bar k_j)}M^{(x)(1)}(x,\bar k_j).
\end{align}

\item Behavior at $0$:
\begin{equation}\label{i_beh-M_(x)__}
M^{(x)}(x,k)=P_0(x)+O(k) , \quad k\to 0,  
\end{equation}
where
\[
P_0(x)=\sqrt{\frac{1+q_0(x)}{2q_0(x)}}\begin{pmatrix}
1 & \frac{w_0(x)}{1+q_0(x)}  \\
-\frac{w_0(x)}{1+q_0(x)} & 1 \\
\end{pmatrix}
\]
with $q_0(x)=\sqrt{1+w_0^2(x)}$.

\end{enumerate}

\begin{remark}\label{rem:x-RH_reduced_left}
    In the case 
    $w_0(0)=u_{0x}(0)=0$, the 
    definition of $M^{(x)}(x,k)$ and the corresponding 
    jump conditions simplify,
    in view of \eqref{Phi_0_inf_reduced}, to the following:
     \begin{equation*}
           M_-^{(x)}(x,k)=M_+^{(x)}(x,k)J(x,k),\quad k\in \mathbb{R} 
        \end{equation*}
        where
\[
M^{(x)}(x,k)=\begin{cases}
\left(\Phi_{\infty 4}^{(1)}(x,0,k), \frac{\Phi_{\infty 2}^{(2)}(x,0,k)}{a^*(k)}\right),\quad k\in\mathbb{C}^+,\\
\left( \frac{\Phi_{\infty 2}^{(1)}(x,0,k)}{a( k)},\Phi_{\infty 4}^{(2)}(x,0,k),\right),\quad k\in\mathbb{C}^-.
\end{cases}
\]
\begin{equation*}
          J(x,k)=\eul^{-\ii kp(x,0,k)\sigma_3} J_0^{(x)}(k)\eul^{\ii kp(x,0,k)\sigma_3}
        \end{equation*} 
        with  
\begin{equation}\label{J_0___simpl}
   J^{(x)}_0( k)=
      \begin{pmatrix}
       1+r(k)r^*(k)&r^*(k)\\r(k)&1
   \end{pmatrix},\quad k\in\mathbb{R}.
\end{equation}
Here $\mathbb{R}$ is oriented from left to right.
\end{remark}

Introducing a new variable
\begin{equation}\label{y_(x)__}
y(x)=\int_0^x q(\xi,0)\dd\xi 
\end{equation}
we obtain the RH problem  explicitly dependent
on $y$ as a parameter:

\textbf{The Riemann--Hilbert problem for $\hat M^{(x)}(y, k)$ (RH$^{(x)}$):} Given $a( k)$, $b( k)$ for $ k\in\mathbb{C}^-$ (satisfying the properties from Section \ref{sub:prop_a}) and the set $\{k_j\}$ in $\mathbb{C}^-$, find a  piece-wise meromorphic $2\times 2$ matrix valued function $\hat M^{(x)}(y, k)$ that satisfies the following conditions:

\begin{enumerate}[(1)]
    \item Jump condition across $\mathbb{R}$
    \begin{subequations}
        \label{jump_hatM_(x)__}
        \begin{equation}
           \hat M_-^{(x)}(y, k)=\hat M_+^{(x)}(y, k)\hat J(y, k),\quad k\in \mathbb{R}\cup\{|k|=\epsilon\} 
        \end{equation}
        where
         \begin{equation}\label{jump_hatM_(x)_matrix__}
          \hat J(y, k)=\eul^{-\ii ky\sigma_3} J_0^{(x)}( k)\eul^{\ii ky\sigma_3}
        \end{equation} 
    \end{subequations}
with $J_0^{(x)}( k)$ given by \eqref{J_0__}.

\item Behavior at $\infty$:
\begin{equation}\label{inf_hatM_(x)__}
     \hat M^{(x)}(y, k)=I+O(\frac{1}{ k}),\quad k\to\infty.
\end{equation}

\item Residue conditions:

\begin{align}\label{res-hatM+_(x)__}
\Res_{ k_j}\hat M^{(x)(1)}(y, k)&=\frac{e^{2\ii k_jy}}{\dot a( k_j)b( k_j)}\hat M^{(x)(2)}(y, k_j),\\
\label{res-hatM-_(x)__}
\Res_{\bar k_j}\hat M^{(x)(2)}(y, k)&=-\frac{e^{-2\ii\bar  k_jy}}{\dot a^*(\bar k_j)b^*(\bar k_j)}\hat M^{(x)(1)}(y,\bar k_j).
\end{align}   
\end{enumerate}

Then evaluating the solution $ \hat M^{(x)}(y,k)$ of the Riemann–Hilbert problem \textbf{RH$^{(x)}$}  at $ k=0$:
    \begin{equation}
       \hat M^{(x)}(y, k)=\begin{pmatrix}
           \hat \alpha(y)&\hat \beta(y)\\
           -\hat \beta(y)&\hat \alpha(y)
       \end{pmatrix}+O(k), 
    \end{equation}
we obtain
\begin{subequations}\label{m_0_via_RH_x_par__}
    \begin{align}\label{m_0_via_RH_x__}
    &\hat{w}_0(y)=\frac{2\hat \alpha(y)\hat \beta(y)}{\hat \alpha^2(y)-\hat \beta^2(y)},\\\label{M_(x)_x_y__}
        &x(y)=\int_0^y \left(\hat \alpha^2(y)-\hat \beta^2(y)\right) \dd y,\\
       &w_0(x)=\hat{w}_0(y(x)).
    \end{align}
    \end{subequations}

\subsubsection{The Inverse $t$-Spectral Mapping}

The formulation of the Riemann–Hilbert problem \textbf{RH$^{(t)}$} is the same as in the case of Problem I (see \eqref{jump_hatM_(t)}--\eqref{res-hatM-_(t)}). In particular, the same procedure of the inverse mapping holds.

\section{The SP equation,  associated Lax pair, and  RH formalism in $y,t$ variables}\label{sec:7}

\subsection{The SP equation and the  associated Lax pair in $y,t$ variables}

Recall that to have the data for the master RH problem to 
be explicitly determined by the initial and boundary data only, 
we introduced the new space variable $y=y(x,t)$ by \eqref{y_}. It is a solution of the system of differential equations 
\begin{subequations}\label{y-xt}
\begin{align}
\label{y_x} y_x & = q,\\
   \label{y_t} y_t &=\frac{u^2q}{2}
\end{align}
    \end{subequations}
    which is compatible due to \eqref{SP}, fixed by the condition $y(0,0)=0$.
Then, by differentiating the identity $x(y(x,t),t)=x$, we derive the equations that describe the 
inverse mapping $x=x(y,t)$:
\begin{subequations}\label{x-yt}
\begin{align}
\label{x_y} x_y & = \frac{1}{\hat q},\\
   \label{x_t} x_t &=-\frac{\hat u^2}{2},
\end{align}
    \end{subequations}
where $\hat u (y,t) = u(x(y,t),t)$ and  $\hat q (y,t) = q(x(y,t),t)$.

\begin{proposition} (SP equation in the $(y,t)$ variables) Let $u(x,t)$, $w(x,t)$  and $q(x,t)$ satisfy \eqref{SP} and let $y(x,t)$ be 
such that \eqref{y-xt} hold. 
Then the SP equation \eqref{SP} in the $(y,t)$ variables reads as the following system of equations
for  $\hat q(y,t):=q(x(y,t),t)$, $\hat u(y,t):=u(x(y,t),t)$, and 
$\hat w(y,t):=w(x(y,t),t)$:

\begin{subequations}\label{SP_in_y}
\begin{align}\label{SP_in_y-1}
\hat q_t&=\hat u \hat q\hat w,\\\label{SP_in_y-2}
\hat q&=\sqrt{\hat w^2+1},\\\label{SP_in_y-3}
\hat u_y&=\frac{\hat w}{\hat q}.
\end{align}
\end{subequations}
    
\end{proposition}

\begin{proof}
As we discussed above, \eqref{y-xt} implies \eqref{x-yt}. 
 Substituting \eqref{SP-2} and \eqref{x_t} into the equality
\[
\hat q_t=\left(q_xx_t+q_t\right)\mid_{x=x(y,t)},
\]
we get \eqref{SP_in_y-1}.

Now, substituting \eqref{x_y} into the equality
$\hat u_y(y,t)=\left(u_x(x,t)x_y(y,t)\right)\mid_{x=x(y,t)}$,
we get \eqref{SP_in_y-3}.

Finally, $\hat q(y,t)=\sqrt{w^2(x(y,t),t)+1}=\sqrt{\hat w^2(y,t)+1}$ and thus \eqref{SP_in_y-2} follows.
\end{proof}

\begin{remark}\label{rem:cons_law}
Notice that we can rewrite \eqref{SP_in_y-1}--\eqref{SP_in_y-3} as
\begin{equation}\label{SP_in_y_alt}
    \left(\frac{1}{\hat q}\right)_t=-\frac{1}{2}\left(\hat u^2\right)_y, \quad \hat q=\sqrt{\hat w^2+1}, \quad \hat u_y=\frac{\hat w}{\hat q}.
\end{equation}
Particularly, assuming  that $\frac{1}{\hat q(\cdot,t)}-1\in L^1(0,\infty)$ for all $t$, the first equation in \eqref{SP_in_y_alt} implies the ``conservation law''
\begin{equation}\label{cons-y}
    \int_0^\infty \left(\frac{1}{\hat q}-1\right)(\xi,t)\dd \xi-\int_0^\infty \left(\frac{1}{\hat q}-1\right)(\xi,0)\dd \xi=\frac{1}{2}\int_0^t \hat u^2(0,s)ds
\end{equation}
regardless the sign of $\hat q$ (understanding \eqref{SP_in_y-2}
as $\hat q^2 = \hat w^2 +1 $).

\end{remark}
\begin{remark} 
   Equations \eqref{SP_in_y-2} and \eqref{SP_in_y-3} imply that $\hat q^2=\frac{1}{1-\hat u_y^2}$, and differentiating \eqref{SP_in_y-2} yields $\hat w_t=\hat u\hat q^2$ under condition $\hat w\not\equiv 0$. Then, differentiating \eqref{SP_in_y-3} w.r.t. $t$ we get
\begin{equation}\label{SP_yt_}
       \hat u_{yt}=\hat u\sqrt{1-\hat u_y^2}     
    \end{equation}
 provided $1-\hat u_y^2>0$. 
 
Vice versa, starting from \eqref{SP_yt_} and introducing $\hat q$ and $\hat w$ by
\[
\hat q=\frac{1}{\sqrt{1-\hat u_y^2}},\quad \hat w=\frac{\hat u_y}{\sqrt{1-\hat u_y^2}}
\]
one recovers \eqref{SP_in_y}.
\end{remark}

\begin{remark}
Introducing $Q$ such that $\sin Q:=\hat u_y$, \eqref{SP_yt_} reduces to $Q_{yt}=\sin Q$, which is the sine–Gordon equation in the light-cone coordinates.
The transformation $x\mapsto y$ exploited in the present Riemann-–Hilbert framework is essentially the same as the one introduced by A. Sakovich and S. Sakovich in \cite{SS05}. 
\end{remark}

Similarly, the reverse  change of variable $(y,t)\mapsto(x,t)$ 
reduces \eqref{SP_in_y} to \eqref{SP} with 
$u(x,t):=\hat u(y(x,t),t)$, $w(x,t):=\hat w(y(x,t),t)$ and $q(x,t):=\hat q(y(x,t),t)$. 
\begin{proposition}\label{prop:hatu_u}
    Let $\hat u(y,t)$, $\hat w(y,t)$,  and $\hat m(y,t)$  satisfy 
    \eqref{SP_in_y} and let $x(y,t)$ be 
such that \eqref{x-yt} holds. Define $u(x,t):=\hat u(y(x,t),t)$, $w(x,t):=\hat w(y(x,t),t)$ and $q(x,t):=\hat q(y(x,t),t)$.
Then \eqref{SP} holds for $u(x,t)$, $w(x,t)$ and $q(x,t)$.
\end{proposition}

Now, let us reformulate the original Lax pair equations in the $(y,t)$ variables.

\begin{proposition}
The Lax pair \eqref{Lax-Q-form} in the variables $(y,t)$ takes the form 
\begin{subequations}\label{Lax_y}
\begin{align}\label{Lax_y_y}
    &\hat \Psi_{ y}=\left(-\ii k \sigma_3+\frac{\hat w_y}{2\hat q^2}\begin{pmatrix}
        0&1\\-1&0
    \end{pmatrix}\right)\hat \Psi,\\\label{Lax_y_t}
   &\hat\Psi_{ t}=-\frac{1}{4\ii k}\begin{pmatrix}
        \frac{1}{\hat q}&-\frac{\hat w}{\hat q}\\-\frac{\hat w}{\hat q}&-\frac{1}{\hat q}
    \end{pmatrix} \hat\Psi.
\end{align}
\end{subequations}
\end{proposition}

\begin{proof}
  Introducing $\hat\Psi(y,t) = \tilde\Phi_\infty (x(y,t),t)$ and taking into account \eqref{x_y} and \eqref{x_t}, the Lax pair \eqref{Lax-Q-form} in the variables $(y,t)$ takes the form \eqref{Lax_y}.
\end{proof}

\subsection{The  RH formalism in $y,t$ variables}

In this subsection, we describe the derivation of a (local) solution to the SP equation in the $(y,t)$ variables, starting from an RH problem parametrized by $y$ and $t$, motivated by the considerations in Section \ref{sec:6}.

Let $\Gamma$ denote an oriented contour in the complex plane, invariant under the symmetries $k\mapsto-k$ and $k\mapsto\bar k$. Suppose that ${\mathbb C}\setminus\Gamma = D_1\cup D_2$ so that $\Gamma$ is the counterclockwise boundary of $D_1$ (the clockwise boundary of $D_2$).

Consider the following \textbf{RH problem parametrized by $y$ and $t$}: Given the $2\times 2$-matrix valued function $J_0(k)$ for $k\in\Gamma$ such that 
\begin{subequations}\label{jump_cond}
    \begin{align}\label{det_jump}
    &\det J_0(k)\equiv 1,\\
    \label{jump_at_inf}
    &J_0(k)=I+O(\frac{1}{k}),\quad k\to\infty,\\
    \label{jump_at_0}
    &J_0(k)=I+O(k),\quad k\to0,
    \end{align}
\end{subequations}
and the sets $\{k_j,c_j\}_1^N$ with $k_j\in D_1$, and $\{\eta_j,d_j\}_1^M$ with $\eta_j\in D_2$,
find a piece-wise meromorphic (w.r.t. $\Gamma$) $2\times 2$-matrix valued function $\hat M(y,t,k)$ satisfying:
\begin{enumerate}[\textbullet]
\item
\emph{Jump} condition
\begin{equation}\label{jump-y_loc}
\hat M_+(y,t,k)=\hat M_-(y,t,\lambda) \hat J(y,t,k),\qquad k\in\Gamma,
\end{equation}
where  
$\hat J(y,t,k)=\eul^{-ik\hat p(y,t,k)\sigma_3}J_0(k)\eul^{ik\hat p(y,t,k)\sigma_3}$ with  
$\hat p(y,t,k)$ defined in \eqref{hat_p}.

\item  \emph{Normalization} condition:
\begin{equation}\label{norm-m-hat_loc}
\hat M(y,t,k)=I+\ord(\frac{1}{k}), \quad k\to\infty.
\end{equation}

\item \emph{Residue} conditions:
\begin{align}\label{res_hatM_loc}
\Res_{ k_j}\hat M^{(1)}(y,t,k)&= c_j \hat M^{(2)}(y,t, k_j) \eul^{-2\ii  k_j \hat p(y,t, k_j)}, ~  k_j\in D_1, \\\label{res_hatM__loc}
\Res_{\eta_j}\hat M^{(2)}(y,t,k)&=d_j \hat M^{(1)}(y,t,\eta_j)\eul^{2\ii \eta_j \hat p(y,t,\eta_j)} , ~ \eta_j\in D_2.
\end{align}
\end{enumerate}

\begin{proposition}
 Assume that $\hat M $ is a solution of RH problem \eqref{jump-y_loc}--\eqref{res_hatM__loc}. Then $\det \hat M=1$.
\end{proposition}
\begin{proof}
The fact that $\det J_0(k)\equiv 1$ implies that $\det \hat M$ has neither a jump across $\Gamma$ nor singularities at the points $k_j$ and $\eta_j$, and thus is an entire function. Therefore, \eqref{jump_at_inf}
together with Liouville's theorem implies $\det \hat M=1$.
\end{proof}

\begin{proposition}
     If a solution of the RH problem \eqref{jump-y_loc}--\eqref{res_hatM__loc} exists, it is unique.
\end{proposition}

\begin{proof}
    The statement follows by considering two solutions $\hat M_1$ and $\hat M_2$ and applying Liouville's theorem to their ratio $\hat M_1 (\hat M_2)^{-1}$.
\end{proof}

Now, assume that the RH problem \eqref{jump-y_loc}--\eqref{res_hatM__loc} 
has a solution $\hat M(y,t,k)$ that satisfies the \emph{symmetries}
\begin{equation}
 \label{sym_1}
 \hat M(k)=\overline{\hat M(-\bar k)}
 \end{equation}
and 
 \begin{equation}
 \label{sym_2}
\hat M(k)=\begin{pmatrix}
    0&-\ii\\\ii&0
\end{pmatrix}\overline{\hat M(\bar k)}\begin{pmatrix}
    0&-\ii\\\ii&0
\end{pmatrix}
 \end{equation} 
and  is differentiable w.r.t. $y$ and $t$.

Now, our goal is to show that, by evaluating the solution $\hat M(y,t,k)$ of the RH problem \eqref{jump-y_loc}--\eqref{res_hatM__loc} at appropriate points of the extended complex plane $\overline{\mathbb C}$, one can 
    obtain a solution of the SP equation in the $y,t$ variables \eqref{SP_in_y}.
For this purpose,
starting from $\hat M(y,t,k)$, we define $2\times 2$-matrix valued functions 
\[ 
\hat \Psi (y,t,k):= \hat M(y,t,k)\eul^{-\ii k \hat p(y,t,k)\sigma_3}, \qquad \hat p(y,t,k)=y-\frac{t}{4k^2},  
\] 
and show that $\hat\Psi(y,t,k)$ satisfies the system of differential equations:
\begin{equation}\label{Lax-hat-hat}
\begin{split}
\hat\Psi_y&=\doublehat{U}\hat\Psi, \\
\hat \Psi_t&=\doublehat{V}\hat\Psi,
\end{split}
\end{equation}
where $\doublehat{U}$ and $\doublehat{V}$ have the same (rational) dependence on $k$ as in \eqref{Lax_y_y} and \eqref{Lax_y_t}, with coefficients given in terms of $\hat M(y,t,k)$ evaluated at appropriate values of $k$.

\begin{proposition}\label{Prop_Lax_y}
    Let $\hat M(y,t,k)$ be the solution of the RH problem \eqref{jump-y_loc}--\eqref{res_hatM__loc} that satisfies symmetries \eqref{sym_1} and \eqref{sym_2}. Define $\hat \Psi (y,t,k):= \hat M(y,t,k)\eul^{-\ii k \hat p(y,t,k)\sigma_3}$. 

    Then $\hat\Psi(y,t,k)$ satisfies the differential equation
\begin{equation}\label{Lax_y_Psi}
    \hat\Psi_y=\doublehat{U}\hat\Psi
\end{equation}
with 
\begin{equation}\label{hat_hat_U}
    \doublehat{U}(y,t,k)=-\ii k \sigma_3 + 2\hat\xi(y,t)\begin{pmatrix}
    0&-1\\1&0
\end{pmatrix},
\end{equation}
where 
$\hat \xi(y,t)\in\mathbb{R}$ can be obtained from the large $k$ expansion of $\hat M(y,t,k)$: 
\begin{equation}\label{xi}
\hat \xi(y,t)=\hat m_2^\infty(y,t),
\end{equation}
where
\begin{equation}\label{xi-M}
\hat M(y,t,k)=I+\frac{\ii}{k}\begin{pmatrix}
    \hat m_1^\infty & \hat m_2^\infty\\
    \hat m_2^\infty& -\hat m_1^\infty
\end{pmatrix}(y,t)+O(\frac{1}{k^2}),\qquad k\to\infty.
\end{equation}

\end{proposition}

\begin{proof}
     First, notice that $\hat\Psi(y,t,\lambda)$ satisfies the jump condition
\[
\hat\Psi^+(y,t,k)=\hat\Psi^-(y,t,k)J_0(k)
\]
with the jump matrix $J_0(k)$ independent of $y$. Hence, $\hat\Psi_y(y,t,k)$ satisfies the same jump condition.

As for the residue conditions, we notice that the symmetry assumptions  
\eqref{sym_1} and \eqref{sym_2} imply that if $ k_j$ is pole of $\hat M^{(1)}$ then $-\bar k_j$ is also pole of $\hat M^{(1)}$ with $c_{ k_j}=-\overline{c_{-\bar k_j}}$. Similarly, if $\eta_j$ is pole of $\hat M^{(2)}$ then $-\bar\eta_j$ is also pole of $\hat M^{(2)}$ with $d_{\eta_j}=-\overline{d_{-\bar\eta_j}}$. Moreover,
\begin{equation}\label{res-sym}
    \eta_j=- k_j, \quad d_j=c_j.
\end{equation}
This implies, $d_{\bar  k_j}=-\overline{c_{ k_j}}$.

Since $\{c_j\}$ are independent of $y$, 
$\hat\Psi_y(y,t,k)$ satisfies the same residue conditions as $\hat\Psi(y,t,k)$ does:

\begin{align*}
\Res_{ k_j}\hat \Psi^{(1)}&= c_j \hat \Psi^{(2)}( k_j) , ~  k_j\in D_1, \\
\Res_{\bar k_j}\hat \Psi^{(2)}&=-\bar{ c}_j \hat \Psi^{(1)}(\bar k_j) , ~ \bar k_j\in D_2.
\end{align*}
with constants $c_j$ independent of $y$.
Consequently,  $\hat \Psi_y \hat \Psi^{-1}=\hat M_y \hat M^{-1}-\ii k \hat p_y\hat M \sigma_3 \hat M ^{-1}$ (with $\hat p_y=1$)  has neither jump no singularities at $ k_j$ and thus it is a meromorphic function, with possible singularity at $k=\infty$.

Let us analyze the behavior of $\hat \Psi_y \hat \Psi^{-1}$ as $k\to\infty$.

\begin{enumerate}[(i)]
    \item As $k\to\infty$, we have due to the symmetry \eqref{sym_1}

    \begin{equation}
   \hat M = I+\frac{\ii}{k}\begin{pmatrix}
    \hat m_{11}^\infty & \hat m_{12}^\infty\\
    \hat m_{21}^\infty& \hat m_{22}^\infty
\end{pmatrix}+O(\frac{1}{k^2}), \quad k\in\mathbb{C}^+
   \end{equation}
with $\hat m_{+ij}^\infty \in\mathbb{R}$. Moreover, the symmetry \eqref{sym_2} implies
    \begin{equation}
   \hat M = I-\frac{\ii}{k}\begin{pmatrix}
    \hat m_{22}^\infty & -\hat m_{21}^\infty\\
    -\hat m_{12}^\infty& \hat m_{11}^\infty
\end{pmatrix}+O(\frac{1}{k^2}), \quad k\in\mathbb{C}^-.
   \end{equation}
Consequently, $
\hat M_y=O(\frac{1}{k})$  and thus 
 $\hat M_y \hat M^{-1}=O(\frac{1}{k})$, which leads to 
the following expansion for $\hat\Psi_y\hat\Psi^{-1}$:

\begin{equation}\label{Psi_y_Psi_inf_C+}
  \hat \Psi_y \hat \Psi^{-1}=-\ii k \sigma_3 + \begin{pmatrix}
    \hat m_{11}^\infty+\hat m_{22}^\infty&-2\hat m_{12}^\infty\\2\hat m_{21}^\infty&-\hat m_{11}^\infty-\hat m_{22}^\infty
\end{pmatrix}+O(\frac{1}{k}),  \quad k\to \infty, \quad k\in\mathbb{C}^+,
\end{equation}
\begin{equation}\label{Psi_y_Psi_inf_C-}
  \hat \Psi_y \hat \Psi^{-1}=-\ii k \sigma_3 - \begin{pmatrix}
    \hat m_{11}^\infty+\hat m_{22}^\infty&2\hat m_{21}^\infty\\-2\hat m_{12}^\infty&-\hat m_{11}^\infty-\hat m_{22}^\infty
\end{pmatrix}+O(\frac{1}{k}),  \quad k\to \infty, \quad k\in\mathbb{C}^-.
\end{equation}
Since $\hat \Psi_y \hat \Psi^{-1}$ is meromorphic in $\mathbb{C}$ this expansions should coincide and thus, 
$\hat m_{11}^\infty=-\hat m_{22}^\infty$ and $\hat m_{12}^\infty=\hat m_{21}^\infty$.

Denoting $\hat m_{1}^\infty=\hat m_{11}^\infty$ and $\hat m_{2}^\infty=\hat m_{12}^\infty$, we get
\begin{equation}\label{Psi_y_Psi_inf}
 \hat \Psi_y \hat \Psi^{-1}=-\ii k \sigma_3 + \begin{pmatrix}
    0&-2\hat m_{2}^\infty\\2\hat m_{2}^\infty&0
\end{pmatrix}+O(\frac{1}{k}),  \quad k\to \infty.
\end{equation}

By Liouville's theorem, \eqref{Psi_y_Psi_inf}  imply \eqref{Lax_y_Psi}--\eqref{xi-M}.
\end{enumerate}
\end{proof}

\begin{remark}\label{rem:expM_1}
    On the other hand, the asymptotic behavior $J_0(k)=I+O(k)$ as $k\to0$ together with the symmetries \eqref{sym_1} and  \eqref{sym_2},  leads to the following expansions for $\hat M$ near $k=0$ in $\mathbb{C}^+$ and $\mathbb{C}^-$:

\begin{align}\label{hatM_at_0_pl}
    &\hat M =B\left(I+\ii k\Gamma_++O(k^2)\right),\quad k\to 0,\quad k\in\mathbb{C}^+,\\\label{hatM_at_0_min}
    &\hat M =B\left(I-\ii k\Gamma_-+O(k^2)\right),\quad k\to 0,\quad k\in\mathbb{C}^-,
\end{align}
where
\begin{equation}
    B=\begin{pmatrix}
        \hat \alpha&\hat \beta\\
        -\hat \beta&\hat \alpha
    \end{pmatrix},\quad 
    \Gamma_+=\begin{pmatrix}
        \hat \gamma_{11}&\hat \gamma_{12}\\
        \hat \gamma_{21}&\hat \gamma_{22}
    \end{pmatrix},\quad 
    \Gamma_-=\begin{pmatrix}
    \hat \gamma_{22}&-\hat \gamma_{21}\\
    -\hat \gamma_{12}&\hat \gamma_{11}
    \end{pmatrix}.
\end{equation}
with $\hat \alpha\in\mathbb{R}$, $\hat \beta\in\mathbb{R}$ and $\hat \gamma_{ij}\in\mathbb{R}$

Substituting these expansions into $ \Psi_y\Psi^{-1}=\hat M_y \hat M^{-1}-\ii k \hat p_y\hat M \sigma_3 \hat M ^{-1}$, we arrive at
\begin{align}
    \Psi_y\Psi^{-1}=B_yB^{-1}+\ii k\left(\pm B_y \Gamma_{\pm}^\sharp B^{-1}\pm B_y\Gamma_{\pm}B^{-1}\pm B\Gamma_{\pm y}B^{-1}-B\sigma_3B^{-1}\right)+O(k^2),\quad k\to 0, \quad k\in\mathbb{C}^\pm,
\end{align}
where, for a $2\times 2$ matrix $A=\begin{pmatrix}
    a&b\\c&d
\end{pmatrix}$, we denote its adjugate by $A^\sharp=\begin{pmatrix}
    d&-b\\-c&a
\end{pmatrix}$.

Therefore,  comparing this with \eqref{hat_hat_U}, we obtain
\begin{align}
 &B_yB^{-1}=2\hat \xi\begin{pmatrix}
     0&-1\\
     1&0
 \end{pmatrix},\\\label{BGamma_rel}
 & \pm B_y\Gamma_{\pm}^\sharp B^{-1}\pm B_y\Gamma_{\pm}B^{-1}\pm B\Gamma_{\pm y}B^{-1}-B\sigma_3B^{-1}=-\sigma_3
\end{align}
      
\end{remark}

\begin{proposition}\label{Prop_Lax_t}
$\hat\Psi(y,t,k)$ defined as in Proposition \ref{Prop_Lax_y} satisfies the differential equation
\begin{equation}\label{Lax_t_Psi}
    \hat\Psi_t=\doublehat{V}\hat\Psi
\end{equation}
with 
\begin{equation}\label{hat_hat_V}
    \doublehat{V}=-\frac{1}{4\ii k} \begin{pmatrix}
        \hat \alpha^2-\hat \beta^2&-2\hat \alpha \hat \beta \\
        -2\hat \alpha \hat \beta&-\hat \alpha^2+\hat \beta^2
    \end{pmatrix},
\end{equation}
where 
\begin{enumerate}[(i)]
    \item 
$\hat \alpha(y,t)\in\mathbb{R}$ and $\hat \beta(y,t)\in\mathbb{R}$ can be obtained from the from the expansion of $\hat M(y,t,k)$ as $k\to 0$: 
\begin{equation}\label{alpha_beta-M}
    \hat M(y,t,k)=\begin{pmatrix}
    \hat \alpha &\hat \beta\\
   -\hat \beta&\hat \alpha
\end{pmatrix}+
O(k),\quad k\to 0.
\end{equation}
\end{enumerate}

\end{proposition}

\begin{proof}
    Similarly to Proposition \ref{Prop_Lax_y}, we notice that $\hat\Psi_t \hat\Psi^{-1}=\hat M_t\hat M^{-1}-\ii k\hat p_t\hat M\sigma_3\hat M^{-1}$ has neither jump nor poles at $ k_j$, and thus it is a meromorphic function, with possible singularities at $k=\infty$ and $k=0$, the latter being due to the singularity of $\hat p_t$ at $k=0$:

Evaluating $\hat\Psi_t\hat\Psi^{-1}$ near these points, we have the following.
\begin{enumerate}[(i)]

\item As $k\to\infty$, we have $\hat M_t\hat M^{-1}=O(\frac{1}{k})$, $p_t(k)=O(\frac{1}{k})$ and thus
\begin{equation}\label{psi-inf-t}
\hat\Psi_t\hat\Psi^{-1}(k)=\ord(k^{-1}),\qquad k\to\infty.
\end{equation}

\item As $k\to 0$, using \eqref{hatM_at_0_pl} and \eqref{hatM_at_0_min}, we have

\begin{align}\label{psi_t_psi_0}
    \Psi_t\Psi^{-1}=\frac{\ii}{4 k}B\sigma_3B^{-1}+B_tB^{-1}-\frac{1}{4}B\left(\pm \sigma_3\Gamma_{\pm}^{\sharp}\pm\Gamma_{\pm}\sigma_3\right)B^{-1}+O(k),\quad k\to 0, \quad k\in\mathbb{C}^\pm
\end{align}
In particular,
\begin{equation}\label{Psi_t_Psi_0}
  \hat \Psi_t \hat \Psi^{-1}=-\frac{1}{4\ii k} \begin{pmatrix}
\hat \alpha^2- \hat \beta^2&-2 \hat \alpha \hat \beta\\-2 \hat \alpha \hat \beta&-\hat \alpha^2+\hat \beta^2
\end{pmatrix}+O(1),  \quad k\to 0.
\end{equation}

Now, \eqref{psi-inf-t} and \eqref{Psi_t_Psi_0} together with the Liouville's theorem imply \eqref{Lax_t_Psi}--\eqref{alpha_beta-M}.

\end{enumerate}

\end{proof}

\begin{remark}\label{rem:expM_2} Since $\hat \Psi_t \hat \Psi^{-1}$ is meromorphic in $\mathbb{C}$, its expansions near $k=0$ in $\mathbb{C}^+$ and $\mathbb{C}^-$ should coincide. Consequently, \[ \sigma_3\Gamma_{+}^{\sharp}+\Gamma_{+}\sigma_3=- \sigma_3\Gamma_{-}^{\sharp}-\Gamma_{-}\sigma_3,\] which implies
   $ \hat \gamma_{11}=-\hat \gamma_{22}$, $
    \hat \gamma_{12}=\hat \gamma_{21}$, and thus
    \begin{equation}
        \label{Gamma_pm}
    \Gamma\coloneqq\Gamma_+=-\Gamma_-,\quad \Gamma^{\sharp}=-\Gamma.
    \end{equation}
In particular, denoting $\hat f_1:=\hat \gamma_{11}$ and $\hat f_2:=\hat \gamma_{12}$, we have 
\begin{equation}\label{M_at_0_exp}
    \hat M(y,t,k)=\begin{pmatrix}
    \hat \alpha &\hat \beta\\
   -\hat \beta&\hat \alpha
\end{pmatrix}\left(I+\ii k\begin{pmatrix}
    \hat f_1& \hat f_2\\
    \hat f_2& -\hat f_1
\end{pmatrix}+
O(k^2)\right),\quad k\to 0.
\end{equation}

Furthermore, comparing \eqref{psi_t_psi_0} with \eqref{psi-inf-t}, we conclude that
\begin{equation}\label{BGamma_rel_t}
  B_tB^{-1}-\frac{1}{4}B\left(\pm \sigma_3\Gamma_{\pm}^{\sharp}\pm\Gamma_{\pm}\sigma_3\right)B^{-1}=0.  
\end{equation}

\end{remark}

\begin{proposition}
    Let $\hat \alpha(y,t)$, $\hat \beta(y,t)$, $\hat f_1(y,t)$ and $\hat f_2(y,t)$ be the functions determined in terms of $\hat M(y,t,k)$ as in Propositions \ref{Prop_Lax_y} and \ref{Prop_Lax_t}. Then they satisfy the following equations:
    \begin{subequations}\label{exp_rel}
        \begin{align}
            \label{exp_rel_1}
&\hat f_{1y}=1-(\hat\alpha^2-\hat \beta^2),\\\label{exp_rel_2_}
&\hat f_{2y}=-2\hat\alpha\hat \beta,\\\label{exp_rel_4}
&\hat \alpha_{t}=\frac{\hat f_2}{2}\hat \beta,\\\label{exp_rel_5}
&\hat \beta_{t}=-\frac{\hat f_2}{2}\hat \alpha,\\
\label{exp_rel_3}
&\hat \alpha^2+\hat \beta^2=1.
        \end{align}
    \end{subequations}
\end{proposition}

\begin{proof}
First, notice that \eqref{exp_rel_3}  follows directly from $\det \hat M=1$.

Next, substituting \eqref{Gamma_pm} into \eqref{BGamma_rel}, yields \eqref{exp_rel_1} and \eqref{exp_rel_2_}, while substituting \eqref{Gamma_pm} into \eqref{BGamma_rel_t}, yields \eqref{exp_rel_4} and \eqref{exp_rel_5}
\end{proof}

\begin{proposition}
    \label{reduce} Let $\hat\alpha^2-\hat\beta^2>0$ (then  \eqref{exp_rel_3} implies that $0<\hat\alpha^2-\hat\beta^2\leq1$).
Introducing  $\hat u(y,t)$, $\hat w(y,t)$, and 
$\hat q(y,t)$
in terms of $\hat \alpha$, $\hat \beta$,  $\hat f_1$, and $\hat f_2$ by
\begin{equation}\label{hmbbgg}
\hat u=-\hat f_2,\quad \hat w = \frac{2\hat\alpha\hat\beta}{\hat\alpha^2-\hat\beta^2},\quad \hat q = \frac{1}{\hat\alpha^2-\hat\beta^2},
\end{equation}
 equations  \eqref{exp_rel} -- \eqref{exp_rel_3} reduce to \eqref{SP_in_y}.
\end{proposition}
    
\begin{proof}
    First, notice that \eqref{exp_rel_4} and   \eqref{exp_rel_5} imply that $\hat f_2= \frac{2\alpha_t}{\beta}=-\frac{2\beta_t}{\alpha}$.
    Now, differentiating  $\hat q$ w.r.t. $t$, we get
    \[\hat q_t=\frac{-2(\hat \alpha\hat \alpha_t-\hat \beta\hat \beta_t)}{(\hat \alpha^2-\hat \beta^2)^2}=\frac{-2\hat \alpha\hat \beta\hat f_2}{(\hat \alpha^2-\hat \beta^2)^2}=\hat u\hat w\hat q,\]
    which gives \eqref{SP_in_y-1}.

    Now, let us compute $\hat u_y$ using \eqref{exp_rel}:
    \[\hat u_y=-\hat f_{2y}=2\hat\alpha\hat\beta=\frac{\hat w}{\hat q}\]
    which gives \eqref{SP_in_y-3}.

    In order to get \eqref{SP_in_y-2} notice that due to \eqref{exp_rel_3} we have $1+\hat w^2=1+\frac{4\hat \alpha^2\hat \beta^2}{(\hat \alpha^2-\hat \beta^2)^2}=\hat q^2$.
\end{proof}

\begin{remark}
    Since $\hat u_y=2\alpha\beta$, 
    the inequality 
    $2|ab|\le a^2 + b^2$ together with 
    \eqref{exp_rel_3} implies that $-1\leq \hat u_y\leq1$.
\end{remark}

\begin{remark}
    If $\hat\alpha^2-\hat\beta^2<0$ (then  \eqref{exp_rel_3} implies that $0>\hat\alpha^2-\hat\beta^2\geq -1$), introducing  $\hat u(y,t)$, $\hat w(y,t)$, and 
$\hat q(y,t)$
in terms of $\hat \alpha$, $\hat \beta$,  $\hat f_1$, and $\hat f_2$ by \eqref{hmbbgg}
 equations  \eqref{exp_rel} -- \eqref{exp_rel_3} reduce to the following system
\begin{align*}
    \hat q_t&=\hat u \hat q\hat w,\\
\hat q&=-\sqrt{\hat w^2+1},\\
\hat u_y&=\frac{\hat w}{\hat q}.
\end{align*}

 Equations $\hat q=-\sqrt{\hat w^2+1}$ and $\hat u_y=\frac{\hat w}{\hat q}$ imply that $\hat q^2=\frac{1}{1-\hat u_y^2}$, and differentiating $\hat q=\red{-}\sqrt{\hat w^2+1}$ w.r.t. $t$ yields $\hat w_t=\frac{\hat q_t\hat q}{\hat w}=\hat u\hat q^2$ under the condition $\hat w\neq 0$. Then, differentiating $\hat u_y=\frac{\hat w}{\hat q}$ w.r.t. $t$ we get
\begin{equation*}
       \hat u_{yt}=\frac{\hat w_t\hat q-\hat w\hat q_t}{\hat q^2}=\hat u\left(\hat q-\frac{\hat w^2}{\hat q}\right)=\frac{\hat u}{\hat q}=\red{-}\hat u\sqrt{1-\hat u_y^2}     
    \end{equation*}
 provided $1-\hat u_y^2\neq0$. 

Therefore, we arrive at 
\begin{equation*}
    \hat u_{yt}=\begin{cases}
        \hat u\sqrt{1-\hat u_y^2}    , \quad \hat\alpha^2-\hat\beta^2>0\\
        -\hat u\sqrt{1-\hat u_y^2}   , \quad \hat\alpha^2-\hat\beta^2<0.  
    \end{cases}
\end{equation*}

 Introduce a new function $Q(y,t)$ by 
 \begin{equation*}
     Q(y,t): \begin{cases}
         \sin Q:=\hat u_y, \quad \cos Q=\sqrt{1-\hat u_y^2}, \quad\text{ if } \quad \hat\alpha^2-\hat\beta^2>0,\\
         \sin Q:=\hat u_y, \quad \cos Q=-\sqrt{1-\hat u_y^2}, \quad\text{ if } \quad \hat\alpha^2-\hat\beta^2<0,
     \end{cases}
 \end{equation*}
 i.e ($\mod 2\pi k$)
  \begin{equation*}
     Q(y,t)= \begin{cases}
         \arcsin \hat u_y, \quad\text{ if } \quad \hat\alpha^2-\hat\beta^2>0,\\
         \pi-\arcsin \hat u_y, \quad\text{ if } \quad \hat\alpha^2-\hat\beta^2<0.
     \end{cases}
 \end{equation*}
  Differentiating the relation $\sin Q=\hat u_y$ with respect to $t$ yields
\[\hat u_{yt}=Q_t \cos Q.\]
On the other hand, noticing that in both cases, 
\[\hat u_{yt}=\hat u \cos Q\]
 (due to the choice of cosine), we conclude that
\[
Q_t=\hat u,
\]
and differentiating the last equation  with respect to $y$ yields 
\[
Q_{yt}=\sin Q.
\]

Introducing the change of variables $x(y,t)$ by 
\[
x_y=-\frac{1}{\hat q}, \quad x_t=\frac{\hat u^2}{2},
\]
(these equations are compatible due to $ \hat q_t=\hat u \hat q\hat w \Leftrightarrow \left(\frac{1}{\hat q}\right)_t=-\frac{1}{2}(\hat u^2)_y$)
and $u(x,t):=\hat u(y(x,t),t)$, $w(x,t):=\hat w(y(x,t),t)$ and $q(x,t):=\hat q(y(x,t),t)$, this equation in $(x,t)$-variables reads as
\[
q_t(x,t)=\red{-}\frac{1}{2}(q(x,t)u^2(x,t))_x,\quad q(x,t)=\red{-}\sqrt{1+w^2(x,t)},\quad w(x,t)=\red{-}u_x(x,t).
\]
\end{remark}

\begin{remark} First, notice that $\hat u(y,t)$ and $\hat u_y(y,t)$ remain bounded for all $(y,t)$.

In case $\hat\alpha^2(y^*,t^*)-\hat\beta^2(y^*,t^*)=0$, we have $x_y(y^*,t^*)=0$, and $\hat q(y,t)$ and $\hat w(y,t)$ become unbounded at $(y^*,t^*)$. 

In particular, the change of variables $x\leftrightarrow y$ becomes non-invertible at this point, and $u_x(x^*,t^*)$ with $x^*=x(y^*,t^*)$ becomes unbounded.

\end{remark}

Now, still working in the case $\hat\alpha^2-\hat\beta^2>0$, introduce the reverse change of variable $(y,t)\mapsto (x,t)$ by
\begin{equation}\label{x_y_}  x(y,t):=\int_{\eta(t)}^y(\hat \alpha^2-\hat \beta^2)(\xi,t)\dd \xi=\int_{\eta(t)}^y\frac{1}{\hat q}(\xi,t)\dd \xi.
\end{equation}
It is clear that $x(\eta(t),t)=0$ and $x_y(y,t)=\frac{1}{\hat q(y,t)}$. Furthermore, differentiating \eqref{x_y_} with respect to 
$t$ and using \eqref{SP_in_y_alt}, we obtain \eqref{x_t}.

Then $u(x,t)=\hat u(y(x,t),t)$, $w(x,t)=\hat w(y(x,t),t)$ and $q(x,t)=\hat q(y(x,t),t)$ satisfy \eqref{SP} by Proposition \ref{prop:hatu_u}.

\begin{remark}
    Exploring the compatibility condition of \eqref{Lax_y_Psi} and \eqref{Lax_t_Psi}, we obtain
    the following equations:
\begin{subequations}\label{rel}
\begin{align}\label{sys_rec1_}
&2\hat \xi_t+\hat \alpha\hat \beta=0,\\\label{sys_rec3_}
&\hat \alpha\hat \alpha_y-\hat \beta\hat \beta_y-4\hat \xi\hat \beta \hat \alpha=0,\\\label{sys_rec4_}
&\hat \alpha_y\hat \beta+\hat \alpha\hat \beta_y+2\hat \xi(\hat \alpha^2-\hat \beta^2)=0.
\end{align}
\end{subequations}
In particular, we have
\begin{equation}
    \label{hmbbgg_2}
    \xi=-\frac{(\alpha\beta)_y}{2(\alpha^2-\beta^2)}=-\frac{\hat u_{yy}\hat q}{4}.
\end{equation}
\end{remark}

\begin{remark}
    Notice that $\hat f_2=-\frac{2\hat\beta_t}{\hat \alpha}$. This means that $\hat u$, $\hat w$, $\hat q$ and $x(y,t)$ can be written in terms of $\hat \alpha$ and $\hat \beta$ only.
\end{remark}

\begin{remark}

Integrating \eqref{SP_in_y_alt} first with respect to $\xi$ from $0$ to $\eta(\tau)$, then integrating the resulting expression with respect to $\tau$ from $0$ to $t$, we obtain
\[
\int_0^t \int_0^{\eta(\tau)}\left(\frac{1}{\hat q}\right)_t (\xi,t)\dd\xi=-\frac{1}{2}\int_0^t \hat u^2(\eta(\tau),\tau)d\tau+\frac{1}{2}\int_0^t \hat u^2(0,\tau)d\tau.
\]
Now, applying Fubini’s theorem to the left-hand side yields
\[
\int_0^t \int_0^{\eta(\tau)}\left(\frac{1}{\hat q}\right)_t (\xi,\tau)\dd\xi\dd\tau=\int_0^{\eta(t)}\int_{\eta^{-1}(\xi)}^t \left(\frac{1}{\hat q}\right)_t (\xi,\tau)\dd\tau\dd\xi=\int_0^{\eta(t)} \frac{1}{\hat q}(\xi,t)\dd\xi-\int_0^{\eta(t)} \frac{1}{\hat q}(\xi,\eta^{-1}(\xi))\dd\xi,
\]
and making the substitution $\xi=\eta(\tau)$ in the second integral, we get
\[
\int_0^t \int_0^{\eta(\tau)}\left(\frac{1}{\hat q}\right)_t (\xi,\tau)\dd\xi\dd\tau=\int_0^{\eta(t)} \frac{1}{\hat q}(\xi,t)\dd\xi-\frac{1}{2}\int_0^t \hat u^2(\eta(\tau),\tau)d\tau.
\]
This implies
\[
\int_0^{\eta(t)} \frac{1}{\hat q}(\xi,t)\dd\xi=\frac{1}{2}\int_0^t \hat u^2(0,\tau)d\tau,
\]
and, thus, assuming that the solution of the RH problem $\hat M(y,t,k)$ is defined for all $y\geq 0$, $t\geq 0$, we can rewrite \eqref{x_y_} as 
\begin{equation}\label{change_l}
x(y,t)=\int_0^{y} \frac{1}{\hat q}(\xi,t)\dd\xi-\frac{1}{2}\int_0^t \hat u^2(0,\tau)d\tau.
\end{equation}

This formula can also be  obtained as the inversion of $y(x,t)$ characterized by \eqref{x-yt} with condition $x(0,0)=0$. Moreover, if one determines $\eta$ as solution of $x(\eta(t),t)=0$, then \eqref{eta_diff} holds.
\end{remark}

\begin{remark}
    Notice that the ``conservation law'' (see Remark \ref{rem:cons_law}) implies that (cf. \eqref{p-alt})
    \[
    x(y,t)=y-\int_y^\infty \left(\frac{1}{\hat q}-1\right)(\xi,t)\dd\xi + \int_0^\infty \left(\frac{1}{\hat q}-1\right)(\xi,0)\dd\xi.
    \]
\end{remark}

Now let us discuss sufficient conditions on the contour $\Gamma$,  jump matrix $\hat J_0(k)$, and residue conditions that provide the symmetries \eqref{sym_1} and \eqref{sym_2}.

\begin{proposition}\label{prop:sym}
    Let the contour $\Gamma$ be invariant w.r.t. $k\mapsto-\bar k$ and $k\mapsto\bar k$, the jump matrix $\hat J_0$ satisfy the symmetries
    \begin{equation}
 \label{sym_J_1}
 \hat J_0(k)=\overline{\hat J_0(-\bar k)}
 \end{equation}
and 
 \begin{equation}
 \label{sym_J_2}
J_0^{{-1}}(k)=\sigma_2 \overline{J_0(\bar k)}\sigma_2,
 \end{equation} 
 and the parameters of the residue conditions satisfy \eqref{res-sym}.

    Then the solution $\hat M$ of the RH problem  \eqref{jump_cond}--\eqref{res_hatM__loc} satisfies symmetries \eqref{sym_1} and \eqref{sym_2}
\end{proposition}

\subsection{Differentiability of the solution of RH problem}\label{sec:dif}

First, we observe that the residue conditions can be eliminated by introducing small circular contours and reformulating them as jump conditions. To fix ideas, we therefore restrict to the case in which no residue conditions are present.

Furthermore, assume that the contour $\Gamma$ is invariant with respect to $k\mapsto-\bar k$ and $k\mapsto\bar k$, the jump matrix $\hat J_0$ satisfy the symmetries
    \begin{equation}
 \label{sym_J__}
 \hat J_0(k)=\overline{\hat J_0(-\bar k)},\quad J_0^{-1}(k)=\sigma_2 \overline{J_0(\bar k)}\sigma_2,
 \end{equation}
and the only unbounded part of $\Gamma$ is a real line. 

In order to obtain a solution of the SP equation \eqref{SP_in_y} in the $(y,t)$ scale, we need to show that $\hat M(y,t,k)$ is once differentiable with respect to $y$ and once differentiable with respect to $t$.

\begin{proposition}\label{prop:diff}
    Assume that $\hat J(y,t,k)$ admits a factorization $\hat J=(v^-)^{-1}v^+$ and $w^+:=v^+-I,~w^-:=I-v^-$ that satisfy the following conditions
\begin{enumerate}[(i)]

\item $v^{\pm}$ are continuous on $\Gamma$;

\item $w^\pm$ are nilpotent,

    \item $w^+,~w^-=O\left(\frac{1}{k^2}\right)$ as $k\to\infty$, $k\in\mathbb{R}$,

    \item $w^+,~w^-=O\left(k^2\right)$ as $k\to 0$, $k\in\mathbb{R}$,

    \item $w^+,~w^-\in L^2(\Gamma)$.   
\end{enumerate}
Then the function $\hat M(y,t,k)$ is $C^1$ in $y$ and $C^1$ in $t$. 
    \end{proposition}

\begin{proof} First, notice that  we can apply Lemma \ref{Lem 5.3} which, together with Zhou's vanishing lemma, ensures the solvability of the
RH problem \eqref{jump_cond}--\eqref{norm-m-hat_loc}. In particular,
\begin{equation}\label{mu_sol}
 \hat\mu(y,t,k)=I+(\mathds{1}-C_w)^{-1}(y,t)C_w(y,t)I
\end{equation}
and the solution $\hat M(y,t,k)$ of the Riemann--Hilbert problem \eqref{jump_cond}--\eqref{norm-m-hat_loc} can be written in terms of $\hat\mu(y,t,k)$ as 
\begin{equation}\label{M_via_mu}
    \hat M(y,t,k)=I+\frac{1}{2\pi\ii}\int_\Sigma \frac{\hat\mu(y,t,k)(w^+ + w^-)(y,t,z)}{z-k}dz.
\end{equation}
Therefore, the product rule together with \textit{(iii)} implies that the statement boils down to proving that the map $\hat \mu(y,t,k)$ defined by \eqref{mu_sol} is $C^1$ in $y$ and $C^1$ in $t$.

For this purpose, notice that  $\hat \mu(y,t,k)$ is a composition of the following maps:
\begin{align*}
&(y,t)\mapsto C_{ w(y,t)}I: (\mathbb{R}^+,[0,T))\to L^2(\Gamma);\\
&(y,t)\mapsto (\mathds{1}-C_{ w(y,t)})^{-1}:(\mathbb{R}^+,[0,T))\to \mathcal{B}(L^2(\Gamma)).
\end{align*}
In turn, $(y,t)\mapsto C_{w(y,t)}I: (\mathbb{R},[0,\infty)\to L^2(\Gamma)$ is a composition of 
\begin{align*}
    &(y,t)\mapsto  w^\pm(y,t):(\mathbb{R}^+,[0,T))\to L^2(\Gamma);\\
    &Q:( w^+, w^-)\mapsto C_{ w }I:L^2(\Gamma)\times L^2(\Gamma) \to L^2(\Gamma)
\end{align*}
and $(y,t)\mapsto (\mathds{1}-C_{ w(y,t)})^{-1}:(\mathbb{R}^+,[0,T))\to \mathcal{B}(L^2(\Gamma))$ is a composition of
\begin{align*}
    &(y,t)\mapsto w^\pm(y,t):(\mathbb{R}^+,[0,T))\to L^\infty(\Gamma);\\
    &F:(w^+, w^-)\mapsto  (\mathds{1}-C_{ w})^{-1}:L^\infty(\Gamma)\times L^\infty(\Gamma)\to \mathcal{B}(L^2(\Gamma)).
\end{align*}

Differentiating $\eul^{2\ii k\hat p(y,t,k)}$, we obtain
\[
\frac{d}{dx} (\eul^{2 \ii k\hat p(y,t,k)})=2\ii k, \quad \frac{d}{dt} (\eul^{2 \ii k\hat p(y,t,k)})=\frac{1}{2\ii k}. 
\]
Thus, taking into account properties \textit{(ii)--(iv)}, we conclude that the map $(y,t)\mapsto  w^\pm(y,t):(\mathbb{R}^+,[0,T))\to L^\infty(\Gamma)$ is differentiable in both $y$ and $t$. Furthermore, combining these properties with the dominated convergence theorem yields the differentiability of $(y,t)\mapsto \hat w^\pm(y,t):(\mathbb{R}^+,[0,T))\to L^2(\Gamma)$.

The maps $F$ and $Q$ are linear, and estimates \eqref{C_w_ess_1} and \eqref{C_w_ess_2} imply their differentiability. Moreover, $dF:(u^+,u^-)\mapsto -C_u:L^\infty(\Gamma)\times L^\infty(\Gamma)\to \mathcal{B}(L^2(\Gamma))$ and $dQ:(u^+,u^-)\mapsto C_u I:L^2(\Gamma)\times L^2(\Gamma)\to L^2(\Gamma)$. Thus, by the chain rule (see \cites{TFA}), the map $(y,t)\mapsto (\mathds{1}-C_{ w(y,t)})^{-1}:(\mathbb{R}^+,[0,T))\to \mathcal{B}(L^2(\Gamma))$ is differentiable together with its inverse. The chain rule again yields the differentiability of the map $(y,t)\mapsto C_{w(y,t)}I: (\mathbb{R}^+,[0,T))\to L^2(\Gamma)$. Therefore, $\hat\mu(y,t,k)$ is differentiable as a composition of differentiable maps.
\end{proof}

\begin{corollary}\label{cor:diff}
\begin{enumerate}[(i)] Assume that the assumptions of Proposition \ref{prop:diff} are satisfied. Then
    \item    The functions $\hat u(y,t)$, $\hat w(y,t)$, $\hat q(y,t)$ defined by \eqref{hmbbgg} are $C^1$ in $y$ and $C^1$ in $t$. 

    \item Assume also that $\hat q(y,t)$ defined by \eqref{hmbbgg} satisfies $\hat q(y,t)\geq c>0$, and $x(y,t)$ is given by\eqref{change_l}.
 Then $u(x,t):=\hat u(y(x,t),t)$, $w(x,t):=\hat w(y(x,t),t)$ and $q(x,t):=\hat q(y(x,t),t)$ are $C^1$ in $x$ and $C^1$ in $t$. 
\end{enumerate}
\end{corollary}

\begin{proof}
    \noindent (i) The claim follows directly from Proposition \ref{prop:diff}.

    \noindent (ii) The condition $\hat q(y,t)\geq c>0$, together with $\hat q(\cdot,t)\in C^1$, implies that $x(y,t)$ defines a $C^1$ bijection.
\end{proof}

\begin{remark}
   To apply Proposition \ref{prop:diff} to the master Riemann--Hilbert problems \textbf{RH$^{(xt)}$}, we first transform them into equivalent problems with continuous jump matrices.
\end{remark}

%-------------------%

\section{The existence result for Problem I}\label{sec:8}

In Section~\ref{sec:3}, we assumed the existence of a solution to Problem~I for the SP equation with prescribed behavior as $x \to +\infty$, and constructed the associated Riemann--Hilbert problem from the corresponding Jost solutions. In this section, we instead consider the Riemann--Hilbert problem generated by initial data from a suitable class and show that the constructed RH problem yields a solution to Problem~I.

To fix ideas,
we assume that $a(k)$ has no zeros in $\overline{\mathbb{C}^+}$.  Zeros in the open upper half-plane can be removed by introducing small circular contours and replacing the corresponding residue conditions with jump conditions. Similarly, zeros on the real axis can be handled by introducing small semicircular contours.

\begin{proposition}
    Assume that $u_0(x)\in H^{3,2}((0,+\infty))$. Assume also that the functions 
$\hat\alpha(y,t)$, $\hat\beta(y,t)$ and $\hat f_2 (y,t)$ evaluated from
    the solution $\hat M^{(xt)}(y,t,k)$ of the RH problem \eqref{jump_hatM_(xt)}--\eqref{res-hatM-_(xt)}, see \eqref{hatM_xt__i}, satisfies $\hat\alpha^2(y,t)-\hat\beta^2(y,t)\geq c>0$ for all $0\leq t<T$ and $y\geq 0$.
Then the problem  
\begin{subequations}\label{SP_IBVP}
\begin{align}\label{SP-2_IBVP}
&q_t = \frac{1}{2}(u^2q)_x,\quad \quad x \geq 0,\quad
     0\leq t < T\leq\infty\\
\label{tq_IBVP}
&q=\sqrt{1+w^2},\\
\label{tw_IBVP}
&w= u_x;\\
     \label{ic_IBVP}
     &u(x,0) = u_0(x), \quad x \geq 0;\\
     \label{as_IBVP}
     &u(x,t)\in \tilde H^{2,2}((0,+\infty)), \quad
     0< t < T<\infty.
\end{align}
\end{subequations}
has a unique solution, $u(x,t)$, which can be represented in terms of the solution of the associated RH problem \eqref{jump_hatM_(xt)}--\eqref{res-hatM-_(xt)} in the parametric form \eqref{hat_u_via_RH_xt}--\eqref{v_j_via_RH_xt}.

\end{proposition} 

\begin{proof} First, notice that the existence of a solution of the RH problem \eqref{jump_hatM_(xt)}--\eqref{res-hatM-_(xt)} follows from Zhou's vanishing lemma (see Proposition \ref{prop:exist}). Next, we construct $u(x,t)$ by following Steps~2--5 of Proposition \ref{uniq_I}. As shown in Section
 \ref{sec:7}, the function $u(x,t)$ constructed in this way indeed satisfies the SP equation \eqref{SP-2_IBVP}--\eqref{tw_IBVP}.

Therefore, it remains to prove that:
\begin{enumerate}[1.]
 \item $u(x,0) = u_0(x)$ for $x \geq 0$.
 
    \item $u(x,t)$, $w(x,t)$ and $q(x,t)$ are $C^1$ in $x$ and $C^1$ in $t$.
    
    \item $u(\cdot,t)\in \tilde H^{2,2}((0,+\infty))$. 
\end{enumerate}

\textbf{1.}
Notice that $\hat{ M}^{(xt)}(y,0,k)$ satisfies the RH problem \textbf{RH$^{(x)}$} (see Remark \ref{rem:jump_coins}). Thus
 by uniqueness of solution of RH problem \textbf{RH$^{(x)}$}, we conclude that $\hat{ M}^{(xt)}(y,0,k)=\hat{ M}^{(x)}(y,k)$. In particular, \eqref{m_0_via_RH_x_par} and \eqref{hat_w_via_RH_xt}--\eqref{v_j_via_RH_xt} imply that
\[
\hat {w}(y,0)=\hat{w}_{0}(y),
\]
and
\[
x(y)=x(y,0),
\]
where $x(y)$ on the left-hand side is the change of the variable associated with the $x$-problem for $w_{0}$, whereas the right-hand side is the change of the variable associated with the master RH problem at $t=0$. 

\textbf{2.} In order to apply the results of Section \ref{sec:dif}, we first modify the RH problem \eqref{jump_hatM_(xt)}--\eqref{res-hatM-_(xt)} so that its jump matrix satisfies the conditions of Proposition \ref{prop:diff}. 
For this purpose, we construct two functions $h$ and $\tilde h$ to improve the behavior at $k=\infty$ and $k=0$ (c.f. \cite{L12}). 

Notice that we can construct these functions in such a way that we achieve the matching near $k=\epsilon$. Indeed,
Let
\[
h_0(k)=\frac{r_1}{\ii k}, \qquad \tilde h_0(k)=\ii \tilde r_1 k-\tilde r_2k^2,
\]
where $r_1$ and $\tilde r_1$ are determined by the expansions $r(k)=\frac{r_1}{\ii k}+O\left(\frac{1}{k^2}\right)$ and $\tilde r(k)\eul^{-2\ii k\nu(0)}=\ii \tilde r_1 k-\tilde r_2k^2+O(k^3)$. 
Define
\[
F(k)=-\frac{\kappa_2^0\eul^{-\ii k \nu(0)}}{a(k)\tilde a(k)}+h_0(k)-\tilde h_0(k),
\]
which is analytic in \(\mathbb{C}^+\cap\{|k-\epsilon|<\epsilon/2\}\). For $n\ge1$, set
\[
G(k)=(k+\ii)^{n+2}F(k), \qquad 
T_nG(k)=\sum_{j=0}^n \frac{G^{(j)}(\epsilon)}{j!}(k-\epsilon)^j,
\]
and
\[
H_0(k)=\frac{T_nG(k)}{(k+\ii)^{n+2}}, \qquad 
H(k)=\tfrac12\big(H_0(k)+\overline{H_0(-\bar k)}\big).
\]
Then
\[
h(k)=\frac{r_1}{\ii k}+H(k), \qquad 
\tilde h(k)=\ii \tilde r_1 k-\tilde r_2k^2
\]
satisfy $f(k)=\overline{f(-\bar k)}$ and are continuous on $|k|=\epsilon$. Moreover,
\begin{subequations}\label{prop_h_til_h}
    \begin{equation}\label{beh_inf_0}
        h(k)=O\left(\frac{1}{k}\right),\ k\to\infty, \qquad 
\tilde h(k)=O(k),\ k\to 0,
    \end{equation}
       \begin{equation}\label{beh_inf_0_2}
       h(k)-r(k)=O\left(\frac{1}{k^2}\right),\ k\to\infty, \qquad 
\tilde h(k)-\tilde r(k)\eul^{-2\ii k\nu(0)}=O(k^3),\ k\to 0,
    \end{equation}
and
\begin{equation}\label{beh_eps}
    h(k)-\tilde h(k)-\frac{\kappa_2^0\eul^{-\ii k \nu(0)}}{a(k)\tilde a(k)}
=O\big((k-\epsilon)^n\big), \quad k\to\epsilon.
\end{equation}
Furthermore, since $r(k)-\tilde r(k)\eul^{-2\ii k\nu(0)}-\frac{\kappa_2^0\eul^{-\ii k \nu(0)}}{a(k)\tilde a(k)}=0$ for $k\in\mathbb{R}$, we also have
\begin{equation}\label{beh_eps_2}
    h(k)-r(k)=-\tilde r(k)\eul^{-2\ii k\nu(0)} + \tilde h(k)+O\big((k-\epsilon)^n\big), \quad k\to\epsilon.
\end{equation}
\end{subequations}

Notice that
\begin{equation}\label{prop_h_1}
 \overline{h(-\bar k)}=h(k),\quad h(k)\in C(\{|k|=\epsilon\}), \quad h(k)=O\left(\frac{1}{k}\right),\quad k\to\infty.   
\end{equation}
\begin{equation}\label{prop_tilh_1}
 \overline{\tilde h(-\bar k)}=\tilde h(k),\quad \tilde h(k)\in C(\{|k|=\epsilon\}), \quad \tilde h(k)=O\left(k\right),\quad k\to 0.   
\end{equation}

Consider a matrix-valued function $\tilde{\hat M}(y,t,k)$ related to
$\hat M(y,t,k)$ by
\begin{equation}\label{hat_M}
  \tilde { \hat M}(y,t,k)=\begin{cases}
       \hat  M(y,t,k)\eul^{-\ii k \hat p(y,t,k)\sigma_3}\begin{pmatrix}
            1&0\\
        -h(k)&1
   \end{pmatrix}\eul^{\ii k \hat p(y,t,k)\sigma_3},\quad \{|k|>\epsilon\}\cap\mathbb{C}_+,\\
        \hat  M(y,t,k)\eul^{-\ii k \hat p(y,t,k)\sigma_3}\begin{pmatrix}
            1&h^*(k)\\
        0&1
\end{pmatrix}\eul^{\ii k \hat p(y,t,k)\sigma_3},\quad \{|k|>\epsilon\}\cap\mathbb{C}_-,\\
            \hat  M(y,t,k)\eul^{-\ii k \hat p(y,t,k)\sigma_3}\begin{pmatrix}
            1&0\\
        -\tilde h(k)&1
   \end{pmatrix}\eul^{\ii k \hat p(y,t,k)\sigma_3},\quad \{|k|<\epsilon\}\cap\mathbb{C}_+,\\
       \hat   M(y,t,k)\eul^{-\ii k \hat p(y,t,k)\sigma_3}\begin{pmatrix}
            1&\tilde h^*(k)\\
        0&1
\end{pmatrix}\eul^{\ii k \hat p(y,t,k)\sigma_3},\quad \{|k|<\epsilon\}\cap\mathbb{C}_-,
    \end{cases}  
    \end{equation}
which satisfies the following RH problem:
\begin{enumerate}
    \item Jump condition across $\mathbb{R}\cup\{|k|=\epsilon\}$
    \begin{subequations}
        \label{jump_M_(xt)_hat}
        \begin{equation}
          \tilde { \hat M}_-^{(xt)}(y,t,k)=\tilde { \hat M}_+^{(xt)}(y,t,k)\tilde { \hat J}^{(xt)}(y,t,k),\quad k\in \mathbb{R}\cup\{|k|=\epsilon\} 
        \end{equation}
        where
         \begin{equation} \label{jump_M_(xt)_matr_hat}
        \tilde { \hat J}^{(xt)}(y,t,k)=\eul^{-\ii k\hat p(y,t,k)\sigma_3} \tilde{\hat J}^{(xt)}_0(k)\eul^{\ii k\hat p(y,t,k)\sigma_3}
        \end{equation} 
        with $\hat p(y,t,k)=y-\frac{t}{4 k^2}$ and
\begin{equation}\label{jump_M_(xt)_matr_0_hat}
\tilde{\hat J}^{(xt)}_0(k)=\begin{cases}
      \begin{pmatrix}
       1&0\\h(k)-r(k)&1
   \end{pmatrix} \begin{pmatrix}
       1&h^*(k)-r^*(k)\\0&1
   \end{pmatrix},\quad k\in\mathbb{R}\cap\{|k| >\epsilon\},\\
   \begin{pmatrix}
       1&\tilde r^*(k)\eul^{2\ii k\nu(0)}-\tilde h^*(k)\\0&1
\end{pmatrix}\begin{pmatrix}
       1&0\\\tilde r(k)\eul^{-2\ii k\nu(0)}-\tilde h(k)&1
   \end{pmatrix},\quad k\in\mathbb{R}\cap\{|k| <\epsilon\},\\
   \begin{pmatrix}
       1&0\\
       h(k)-\tilde h(k)-\frac{\kappa_2^0\eul^{-\ii k\nu(0)}}{a(k)\tilde a(k)}&1
   \end{pmatrix},\quad k\in\mathbb{C}^+\cap\{|k| =\epsilon\},\\
   \begin{pmatrix}
       1&h^*(k)-\tilde h^*(k)-\frac{\kappa_2^0\eul^{\ii k\nu(0)}}{a^*(k)\tilde a^*(k)}\\
       0&1
\end{pmatrix},\quad k\in\mathbb{C}^-\cap\{|k| =\epsilon\}
   \end{cases} 
\end{equation}
and $r(k)=\frac{b^*(k)}{a(k)}$, $\tilde r(k)=\frac{\tilde b^*(k)}{\tilde a(k)}$.
    \end{subequations}

 \item Behavior at $\infty$:
\begin{equation}\label{inf_hatM_(xt)_hat}
     \tilde{\hat M}(y,t,k)=I+O(\frac{1}{ k}),\quad k\to\infty.
\end{equation}

\item Residue conditions:

\begin{align}\label{res-tilhatM+_(xt)}
\Res_{ k_j} \tilde{\hat M}^{(xt)(1)}(y,t, k)&=\frac{e^{2\ii k_j\left(y-\frac{t}{4k_j^2}\right)}}{\dot a( k_j)b( k_j)} \tilde{\hat M}^{(xt)(2)}(y, t,k_j),\\
\label{res-tilhatM-_(xt)}
\Res_{\bar k_j} \tilde{\hat M}^{(xt)(2)}(y,t, k)&=-\frac{e^{-2\ii \bar k_j\left(y-\frac{t}{4\bar k_j^2}\right)}}{\dot a^*(\bar k_j)b^*(\bar k_j)} \tilde{\hat M}^{(xt)(1)}(y,t,\bar k_j).
\end{align}   
    
\end{enumerate}
Now, we can apply the results of Section \ref{sec:dif} to conclude that $u(x,t)$, $w(x,t)$ and $q(x,t)$ are $C^1$ in $x$ and $C^1$ in $t$.

\textbf{3.} First, notice that the following analog  of Lemmas 2.3 and 2.9 from 
\cite{Zhou98} holds:

\begin{lemma}\label{lem_est} Assume that $\Gamma$ is the contour shown in Figure \ref{fig:contour_RH_x}.
 Assume also that $\hat J(y,t,k)$ admits a factorization $\hat J=(v^-)^{-1}v^+$ and $w^+:=v^+-I,~w^-:=I-v^-$ that satisfy the following conditions
\begin{enumerate}[(i)]

\item $w^\pm$ are nilpotent;

        \item $w^\pm(k)\in H^n(\mathbb{R})$, $w^\pm(k)\in H^n(\Gamma\setminus\mathbb{R})$;

        \item $w^\pm(k)$ vanishes to order at least $2n$ at $k=0$;

        \item $w^\pm(k)$  decays exponentially fast as $y\to +\infty$ ($y\to -\infty$) on $\Gamma\setminus\mathbb{R}$ away from the real line;

        \item $w^\pm(k)=O((k\pm\epsilon)^{n-1})$.

    \end{enumerate}
    Then 
\begin{equation}
    \label{est_1}
    \|C_- w^+\|_{L^2}\leq \frac{\|w^+\|_{H^n}}{(1+y^2)^{\frac{n}{2}}}, \quad \|C_+ w^-\|_{L^2}\leq \frac{\|w^-\|_{H^n}}{(1+y^2)^{\frac{n}{2}}}.
\end{equation}
\end{lemma}

The proof follows the same arguments as in \cite{Zhou98} and is therefore omitted. We note, however, that one must additionally assume that $w^\pm(k)$ vanishes to order at least $2n$ at $k=0$ to ensure that
\[
w^\pm(k)e^{-\frac{i t}{2k}} \in H^n(\mathbb{R}).
\]

Notice that
\begin{equation}\label{hatM_via_mu}
     \tilde{\hat M}(y,t,k)=I+\frac{1}{2\pi\ii}\int_\Sigma \frac{\tilde \mu(y,t,z)(\tilde w^+ + \tilde w^-)(y,t,z)}{z-k}dz,
\end{equation}
where
$w^\pm$ are given by
\begin{equation}\label{w+}
\tilde w^-(y,t,k)=\begin{cases}
       \begin{pmatrix}
       0&(h^*(k)-r^*(k))\eul^{-2\ii k\hat p(y,t,k)}\\0&0
   \end{pmatrix},\quad k\in\mathbb{R}\cap\{|k| >\epsilon\},\\
   \begin{pmatrix}
       0&0\\(\tilde r(k)\eul^{-2\ii k\nu(0)}-\tilde h(k))\eul^{2\ii k\hat p(y,t,k)}&0
   \end{pmatrix},\quad k\in\mathbb{R}\cap\{|k| <\epsilon\},\\
   0,\quad k\in\mathbb{C}^+\cap\{|k| =\epsilon\},\\
   \begin{pmatrix}
       0&\left(h^*(k)-\tilde h^*(k)-\frac{\kappa_2^0\eul^{\ii k\nu(0)}}{a^*(k)\tilde a^*(k)}\right)\eul^{-2\ii k\hat p(y,t,k)}\\
       0&0
   \end{pmatrix},\quad k\in\mathbb{C}^-\cap\{|k| =\epsilon\}
   \end{cases} 
\end{equation}
and
\begin{equation}\label{w-}
\tilde w^+(y,t,k)=\begin{cases}
       \begin{pmatrix}
       0&0\\(r(k)-h(k))\eul^{2\ii k\hat p(y,t,k)}&0
   \end{pmatrix},\quad k\in\mathbb{R}\cap\{|k| >\epsilon\},\\
   \begin{pmatrix}
       0&(\tilde h^*(k)-\tilde r^*(k)\eul^{2\ii k\nu(0)})\eul^{2\ii k\hat p(y,t,k)}\\0&0
   \end{pmatrix},\quad k\in\mathbb{R}\cap\{|k| <\epsilon\},\\
   \begin{pmatrix}
0&0\\-\left(h(k)-\tilde h(k)-\frac{\kappa_2^0\eul^{-\ii k\nu(0)}}{a(k)\tilde a(k)}\right)\eul^{2\ii k\hat p(y,t,k)}&0
   \end{pmatrix},\quad k\in\mathbb{C}^+\cap\{|k| =\epsilon\},\\
   0,\quad k\in\mathbb{C}^-\cap\{|k| =\epsilon\}.   \end{cases} 
\end{equation}
In particular, the properties of $h$ and $\tilde h$ together with the dominated convergence theorem (limits are not tangential) imply that
\begin{equation}\label{xi_via_mu_}
\begin{aligned}
\hat \xi(y,t)=&\frac{2}{\pi}\int_{\mathbb{R}\setminus(-\epsilon,\epsilon)} \tilde \mu_{11}(y,t,z) \left(r^*(z)-h^*(z)\right)\eul^{-2\ii z y+\frac{\ii t}{2z}} \dd z-\\
 &\frac{2}{\pi}\int_{(-\epsilon,\epsilon)} \tilde \mu_{11}(y,t,z) \left(\tilde r^*(z)\eul^{2\ii k\nu(0)}-\tilde h^*(z)\right)\eul^{-2\ii z y+\frac{\ii t}{2z}} \dd z -\\
  &\frac{2}{\pi}\int_{|z|=\epsilon, z\in \mathbb{C}^-} \tilde \mu_{11}(y,t,z) \left(h^*(k)-\tilde h^*(k)-\frac{\kappa_2^0\eul^{\ii k\nu(0)}}{a^*(k)\tilde a^*(k)}\right)\eul^{-2\ii z y+\frac{\ii t}{2z}}\dd z.
\end{aligned}   
\end{equation}
\begin{equation}\label{alpha_via_mu_}
\begin{aligned}
\hat  \alpha(y,t)=&1+\frac{1}{2\pi\ii}\int_{\mathbb{R}\setminus(-\epsilon,\epsilon)} \tilde \mu_{21}(y,t,z) \frac{(r^*(z)-h^*(z))}{z}\eul^{-2\ii z y+\frac{\ii t}{2z}} \dd z-\\
 &\frac{2}{\pi}\int_{(-\epsilon,\epsilon)}\tilde \mu_{21}(y,t,z) \frac{(\tilde r^*(z)\eul^{2\ii k\nu(0)}-\tilde h^*(z))}{z}\eul^{-2\ii z y+\frac{\ii t}{2z}} \dd z -\\
  &\frac{2}{\pi}\int_{|z|=\epsilon, z\in \mathbb{C}^+}\tilde  \mu_{21}(y,t,z) \frac{\left(h^*(k)-\tilde h^*(k)-\frac{\kappa_2^0\eul^{\ii k\nu(0)}}{a^*(k)\tilde a^*(k)}\right)}{z}\eul^{-2\ii z y+\frac{\ii t}{2z}} \dd z.
\end{aligned}   
\end{equation}
\begin{equation}\label{beta_via_mu_}
\begin{aligned}
\hat \beta(y,t)=&\frac{1}{2\pi\ii}\int_{\mathbb{R}\setminus(-\epsilon,\epsilon)} \tilde \mu_{11}(y,t,z) \frac{(r^*(z)-h^*(z))}{z}\eul^{-2\ii z y+\frac{\ii t}{2z}} \dd z-\\
 &\frac{2}{\pi}\int_{(-\epsilon,\epsilon)}\tilde \mu_{11}(y,t,z) \frac{(\tilde r^*(z)\eul^{2\ii k\nu(0)}-\tilde h^*(z))}{z}\eul^{-2\ii z y+\frac{\ii t}{2z}} \dd z -\\
  &\frac{2}{\pi}\int_{|z|=\epsilon, z\in \mathbb{C}^+}\tilde  \mu_{11}(y,t,z) \frac{\left(h^*(k)-\tilde h^*(k)-\frac{\kappa_2^0\eul^{\ii k\nu(0)}}{a^*(k)\tilde a^*(k)}\right)}{z}\eul^{-2\ii z y+\frac{\ii t}{2z}} \dd z.
\end{aligned}   
\end{equation}

Again, following the ideas of 
\cite{Zhou98} (see Theorem 2.1 ), we conclude that the following holds:

\begin{lemma}\label{lem_est_2}
    Assume that $r(k)\in H^n(\mathbb{R})$ and $r(k)$ vanishes to order at least $2n$ at $k=0$. Then, for all fixed $t$
\begin{equation}
    \label{est_2}
  \hat  \xi(y,t)\in L^2((1+y^{2n})\dd y,(0,+\infty)).
\end{equation}
Assume additionally that $r(k)$ vanishes to order at least $2n+1$ at $k=0$.
then
\begin{equation}
    \label{est_3}
    \hat \alpha(y,t)-1,~\hat \beta(y,t)\in L^2((1+y^{2n})\dd y,(0,+\infty)).
\end{equation}
\end{lemma}

Under the assumption that $u_{0x}(x) \in H^{2,2}((0,+\infty))$, one can proceed as in \cite{L12} to show that $\frac{d}{dk} r(k) \in L^2(\mathbb{R})$. It then follows from Lemmas~\ref{lem_est} and \ref{lem_est_2} (with $n=1$) that
\begin{equation}\label{beh_xi_alp_beta_2}
        \hat \xi(y,t),~\hat \alpha(y,t)-1,~\hat \beta(y,t)\in L^2((1+y^{2})\dd y,(0,+\infty)).
    \end{equation}
Therefore, $\hat\alpha^2(y,t)-\hat\beta^2(y,t)\geq c>0$ together with \eqref{hmbbgg} and \eqref{hmbbgg_2} yields
\begin{equation}\label{beh_w_q_u_yy}
        \hat w(y,t),~\frac{1}{\hat q(y,t)}-1,~\hat u_{yy}(y,t)\in L^2((1+y^{2})\dd y,(0,+\infty)).
    \end{equation}
Moreover, since $1\leq \hat q(y,t)=\frac{1}{\hat\alpha^2(y,t)-\hat\beta^2(y,t)}\leq \frac{1}{c}$, we can conclude that  $x(y,t)\sim y$ as $y\to+\infty$. Consequently, 
\begin{equation}\label{beh_w_q_u_yy__}
         u_x(x,t),~ u_{xx}(x,t)\in L^2((1+x^{2})\dd x,(0,+\infty)).
    \end{equation}
\end{proof}

\begin{remark}
    It is the positivity assumption $\hat\alpha^2(y,t)-\hat\beta^2(y,t)\geq c>0$
which can be considered as a specific feature in the implementation of the RHP
approach to the characterization of initial boundary problems for the SP equation.
\end{remark}

\section{The existence result for Problem II}\label{sec:9}

\begin{proposition}
    Let $\{u_0(x)\in H^{3,2}((-\infty,0)); ~v_0(t)\in  H^{2}(0, T), ~v_1(t)\in  H^{2}(0, T)\}$ be a set of functions satisfying internal compatibility  conditions
        \begin{enumerate}[(i)]
    \item $\partial_x^j u_0(0)=v_j(0)$, $j=0,1$
    \item if $v_0(t)=0$ on $(T_1,T_2)$, then $v_{1t}(t)=0$ on $(T_1,T_2)$. 
\end{enumerate}
Assume also that the functions 
$\hat\alpha(y,t)$, $\hat\beta(y,t)$ and $\hat f_2 (y,t)$ evaluated from
    the solution $\hat M^{(xt)}(y,t,k)$ of the RH problem \eqref{jump_hatM_(xt)__}--\eqref{res-hatM-_(xt)__}, see \eqref{hatM_xt__i}, satisfies $\hat\alpha^2(y,t)-\hat\beta^2(y,t)\geq c>0$ for all $0\leq t<T$ and $y\leq \eta(t)$. Then the problem
    \begin{subequations}\label{SP_IBVP__}
\begin{align}\label{SP-2_IBVP__}
&q_t = \frac{1}{2}(u^2q)_x,\\
\label{tq_IBVP__}
&q=\sqrt{1+w^2},\\
\label{tw_IBVP__}
&w= u_x;\\
\label{boundary_IBVP__}
    &u(0,t) = v_0(t), \quad u_x(0,t) = v_1(t), \quad
     0\leq t < T\leq\infty;\\
     \label{ic_IBVP__}
     &u(x,0) = u_0(x), \quad x \leq 0;\\
     \label{as_IBVP__}
     &u(x,t)\in \tilde H^{2,2}((-\infty,0)), \quad
     0\leq t < T<\infty.
\end{align}
\end{subequations}
has a unique solution, $u(x,t)$, which can be represented in terms of the solution of assotiated RH problem \eqref{jump_hatM_(xt)__}--\eqref{res-hatM-_(xt)__} in the parametric form \eqref{hat_u_via_RH_xt__}--\eqref{v_j_via_RH_xt__}.

\end{proposition}

\begin{proof}

The existence of a solution to the RH problem \eqref{jump_hatM_(xt)__}--\eqref{res-hatM-_(xt)__}, the fact that $u(x,t)$ satisfies the SP equation, as well as its differentiability and decay as $x\to -\infty$, can all be established by arguments similar to those used for Problem~I.

Therefore, it remains to prove that:
\begin{enumerate}[1.]

    \item  $u(x,0)=u_0(x)$.

    \item $ u(0,t) = v_0(t)$ and $u_x(0,t) = v_1(t)$.
\end{enumerate}

\textbf{1.} Following the approach used for the right half-line problem,
we establish a relation between the
 solution $\hat{ M}^{(xt)}(y,t,k)$ of the \textbf{RH$^{(xt)}$} problem evaluated at $t = 0$ and the solution $\hat M^{(x)}(y,k)$ of the \textbf{RH$^{(x)}$} in the form of the multiplication by an matrix factor:
 \begin{subequations}
     \begin{align}
         &P^{(x)}(y,k)=\eul^{-\ii k\hat p(y,0,k)\sigma_3}P_0^{(x)}(k)\eul^{\ii k\hat p(y,0,k)\sigma_3},\\
         & P_0^{(x)}(k)=\begin{cases}
        \begin{pmatrix}
           1&-\frac{B(k)}{a^*(k)d^*(k)}\\
           0&1
        \end{pmatrix},\quad k\in\mathbb{C}^+\setminus\{|k|<\epsilon\}\\
\begin{pmatrix}
           1&-\frac{\tilde B(k)\eul^{2\ii k\nu(0)}}{\tilde a^*(k)\tilde d^*(k)}\\
           0&1
        \end{pmatrix},\quad k\in\mathbb{C}^+\cap\{|k|<\epsilon\}\\
        \begin{pmatrix}
           1&0\\
           \frac{B^*(k)}{a(k)d(k)}&1
        \end{pmatrix},\quad k\in\mathbb{C}^-\setminus\{|k|<\epsilon\},\\
      \begin{pmatrix}
           1&0\\
           \frac{\tilde B^*(k)\eul^{-2\ii k\nu(0)}}{\tilde a(k)\tilde d(k)}&1
        \end{pmatrix},\quad k\in\mathbb{C}^-\cap\{|k|<\epsilon\}.
     \end{cases}
     \end{align}
 \end{subequations}

Thus, the problem reduces to show that  $\hat{ N}^{(xt)}(y,k)$ defined by
 \begin{equation}
   \hat{ N}^{(xt)}(y,k)=\hat{ M}^{(xt)}(y,0,k)P^{(x)}(y,k)
 \end{equation}
 satisfies the Riemann--Hilbert problem \textbf{RH$^{(x)}$}.

 \begin{itemize}
     \item The jump conditions match by construction.

     \item Noticing that $-\frac{B(k)}{a^*(k)d^*(k)}\eul^{-2\ii k y}=O(\frac{1}{k})$ for $y\leq0$ (notice that $\eta(0)=0$) in $\mathbb{C}^+$, the normalization condition is satisfied.

     \item Let's check that the residue conditions match.

    For $k\in \mathbb{C}^+$, we have

    \begin{align}
     \hat{ N}^{(xt)(1)}(y,k)&=\hat{M}^{(xt)(1)}(y,0,k),\\
     \hat{ N}^{(xt)(2)}(y,k)&=-\frac{B(k)\eul^{-2\ii k y}}{a^*(k)d^*(k)}\hat{M}^{(xt)(1)}(y,0,k)+\hat{M}^{(xt)(2)}(y,0,k).
 \end{align}
Taking into account the residue condition \eqref{res-hatM-_(xt)__}, the singularities of the second column at $k=\bar\eta_j$ cancel.

On the other hand, at $k=\bar k_j$, the second column is singular due to the singularity of $\frac{1}{a^*(k)}$, and, using $d^*(\bar k_j)=b^*(\bar k_j)B(\bar k_j)$, the corresponding residue condition takes the form \eqref{res-hatM-_(x)__}. 

The singularities in $\mathbb{C}^-$  can be treated in a similar way.

 \end{itemize}

 Thus, by uniqueness of solution of RH problem \textbf{RH$^{(x)}$}, we conclude that $\hat{ M}^{(xt)}(y,0,k)P^{(x)}(y,k)=\hat{ M}^{(x)}(y,k)$.

 Recall that $w_0(x)$ is expressed in terms of the first coefficient in the expansion of $\hat{ M}^{(x)}(y,k)$ near $k=0$, and thus it is important to control $P^{(t)}(y,k)$ near this point. In this regard, notice that 
\[
 -\frac{\tilde B(k)}{\tilde a^*(k)\tilde d^*(k)}\eul^{2\ii k\nu(0)}=O(k)
 \]
and thus
\[
P^{(x)}(y,k)=I+O(k),
\]
which implies that
\[
\hat {w}(y,0)=\hat{w}_{0}(y), \quad x(y)=x(y,0).
\]

\textbf{2.}
Introduce
\[
     P^{(t)}_0(k)=\begin{cases}
         \begin{pmatrix}
             \frac{A(k)}{d^*(k)}&0\\
             b^*(k)&\frac{d^*(k)}{A(k)}
         \end{pmatrix},\quad k\in\mathbb{C}^+\cap\{|k|>\epsilon\},\\
                \begin{pmatrix}
             \frac{\tilde A(k)\eul^{\ii k\nu(0)}}{\tilde d^*(k)}&0\\
             \tilde b^*(k)\eul^{-\ii k\nu(0)}&\frac{\tilde d^*(k)\eul^{-\ii k\nu(0)}}{\tilde A(k)}
         \end{pmatrix},\quad k\in\mathbb{C}^+\cap\{|k|<\epsilon\},\\
                \begin{pmatrix}
             \frac{d(k)}{A^*(k)}&-b(k)\\
            0&\frac{A^*(k)}{d(k)}
         \end{pmatrix},\quad k\in\mathbb{C}^-\cap\{|k|>\epsilon\},\\
                  \begin{pmatrix}
             \frac{\tilde d(k)\eul^{\ii k\nu(0)}}{\tilde A^*(k)}&-\tilde b(k)\eul^{\ii k\nu(0)}\\
            0&\frac{\tilde A^*(k)\eul^{-\ii k\nu(0)}}{\tilde d(k)}
         \end{pmatrix},\quad k\in\mathbb{C}^-\cap\{|k|<\epsilon\},
     \end{cases}
 \]
and define  
\[
   \hat{\hat{ N}}^{(xt)}(t,k):=\hat{ M}^{(xt)}(\eta(t),t,k)P^{(t)}(t,k)
 \]
with $P^{(t)}(t,k)=\eul^{-\ii k\hat p(\eta(t),t,k)\sigma_3}P^{(t)}_0(k)\eul^{\ii k\hat p(\eta(t),t,k)\sigma_3}$, where $\hat p(\eta(t),t,k) = \eta(t)-\frac{t}{4k^2}$.

 \begin{enumerate}
     \item The jump condition 
     for $\hat{\hat{ N}}^{(xt)}(t,k)$ coincides with that for $\hat{ M}^{(t)}(\eta(t),t,k)$
    by construction.

    \item Noticing that $\frac{A(k)}{d^*(k)}=1+O(\frac{1}{k})$ and $b^*(k)\eul^{2\ii k\left(\eta(t)-\frac{t}{4 k^2}\right)}=O(\frac{1}{k})$ in $\mathbb{C}^+$ for all $t\geq0$, the normalization condition is satisfied.

    \item Let's check that the residue conditions match.

    For $k\in \mathbb{C}^+$, we have

      \begin{align}
     \hat{\hat{ N}}^{(xt)(1)}(t,k)&=\frac{A(k)}{d^*(k)}\hat{M}^{(xt)(1)}(\eta(t),t,k)+b^*(k)\eul^{2\ii k\left(\eta(t)-\frac{t}{4 k^2}\right)}\hat{M}^{(xt)(2)}(\eta(t),t,k),\\
     \hat{\hat{ N}}^{(xt)(2)}(t,k)&=\frac{d^*(k)}{A(k)}\hat{M}^{(xt)(2)}(\eta(t),t,k).
 \end{align}
Taking into account the residue condition \eqref{res-hatM-_(xt)__} together with $A(\bar\eta_j)=-\frac{b^*(\bar\eta_j) B(\bar\eta_j)}{a^*(\bar\eta_j)}$, the singularities of the first column at $k=\bar\eta_j$ cancel.

On the other hand, at $k= \mu_j$, the second column is singular due to the singularity of $\frac{1}{A(k)}$, and, using $d^*(\mu_j)=b^*(\mu_j)B(\mu_j)$, the corresponding residue condition takes the form \eqref{res-hatM+_(t)}. 

The singularities in $\mathbb{C}^-$  can be treated in a similar way.
    \end{enumerate}

    By the uniqueness of solution of the RH problem \textbf{RH$^{(t)}$}, it follows  that
 $\hat{\hat{ N}}^{(xt)}(t,k)=\hat{ M}^{(t)}(\eta(t),t,k)$.

Since $v_0(t)$ and $v_1(t)$ are expressed in terms of the first two coefficients in the expansion of $\hat{ M}^{(t)}(\eta(t),t,k)$ near $k=0$, it is important to control $P^{(t)}(t,k)$ near this point. 
In this regard, notice that in $\mathbb{C}^+$ as $k\to 0$
\[
\frac{\tilde A(k)\eul^{\ii k\nu(0)}}{\tilde d^*(k)}=1+O(k) 
\]
and
\[
\tilde b^*(k)\eul^{-\ii k\nu(0)}\eul^{2\ii k\left(\eta(t)-\frac{t}{4 k^2}\right)}=O(k)\eul^{- \frac{\ii t}{2k}}.
\]
In particular, $\tilde b^*(k)\eul^{-\ii k\nu(0)}\eul^{2\ii k\left(\eta(t)-\frac{t}{4 k^2}\right)}$ decays exponentially fast as $k\to 0$ in $\mathbb{C}^+$ for all $t>0$. Therefore,
 \[
 P^{(t)}(t,k)= I+\begin{pmatrix}
         O(k)& 0\\
         O(k)\eul^{- \frac{\ii t}{2k}} &O(k)
     \end{pmatrix}, \quad k\to 0,
 \]
 and the expressions for $\hat u(\eta(t),t)$, $\hat q(\eta(t),t)$ and $\hat w(\eta(t),t)$ given in \eqref{sol_via_RH_xt__}
 remain unchanged.   
\end{proof}

\begin{remark}
    Remarkably, we do not need to impose additional conditions (``global relations'') to prove that the boundary conditions hold.

\end{remark}

%\section{Concluding Remarks}

\appendix

\section{Basic facts about $L^2$-RH problems.}

For our purposes, we will need several results concerning the equivalence between
$L^2$ Riemann--Hilbert problems (whose jump matrices may be discontinuous at
intersection points of the contour) and singular integral equations; see, for
example, \cites{L18,ZH89}. We collect the necessary definitions and results
below.

Throughout this section $\Gamma$ denotes an oriented contour consisting of a
finite union of piecewise $C^1$ arcs. The arcs may intersect only at their
endpoints, and all intersections are assumed to be transversal. In particular,
this class of contours includes configurations such as the union of the real
line with finitely many circles.

The orientation of $\Gamma$ is chosen so that
\[
\Gamma=\partial D_+=-\partial D_-,
\]
where $\mathbb{C}\setminus\Gamma=D_+\cup D_-$.

For a bounded component $D$ of $\mathbb{C}\setminus\Gamma$, denote by $E^2(D)$
the class of analytic functions $f$ in $D$ for which there exists a sequence of
rectifiable Jordan curves $\{C_n\}\subset D$ such that for every compact set
$D_c\subset D$ there exists $N$ so that $C_n$ surrounds $D_c$ for all $n>N$ and
\[
\sup_{n\ge1}\int_{C_n}|f(z)|^2|dz|<\infty .
\]

If $D$ is unbounded, define $f\in E^2(D)$ if
$f\circ\phi^{-1}\in E^2(\phi(D))$, where
\[
\phi(z)=\frac{1}{z-z_0}
\]
for some $z_0\in\mathbb{C}\setminus\overline{D}$. We also write
\[
\dot E^2(D)=\{f\in E^2(D):zf(z)\in E^2(D)\},
\]
that is, $f(z)\to0$ as $z\to\infty$.

Denote by $\mathcal{B}(L^2(\Gamma))$ the space of bounded linear operators on
$L^2(\Gamma)$ and by $\mathcal{F}(L^2(\Gamma))$ the set of Fredholm operators
on $L^2(\Gamma)$. The index map
\[
\mathrm{Ind}:\mathcal{F}(L^2(\Gamma))\to\mathbb{Z}
\]
is constant on the connected components of $\mathcal{F}(L^2(\Gamma))$.

For $h\in L^2(\Gamma)$ the Cauchy transform is defined by
\begin{equation}\label{Cauchy}
(Ch)(\lambda)=\frac{1}{2\pi i}\int_\Gamma \frac{h(z)}{z-\lambda}\,dz .
\end{equation}
The non-tangential limits of $Ch$ as $\lambda$ approaches $\Gamma$ from the
left and right exist and will be denoted by $C_+h$ and $C_-h$, respectively.
The operators $C_\pm$ are bounded on $L^2(\Gamma)$ and satisfy
\[
C_+-C_-=\mathds{1}, \quad  (C_+)^2=C_+,\quad  (C_-)^2=-C_-, \quad C_+C_-=C_-C_+=0.
\]

For functions $w^\pm\in L^2(\Gamma)\cap L^\infty(\Gamma)$ define the operator
$C_w:L^2(\Gamma)+L^\infty(\Gamma)\to L^2(\Gamma)$ by
\begin{equation}\label{C-w}
C_w(f)=C_+(fw^-)+C_-(fw^+).
\end{equation}
Moreover,
\begin{subequations}
\begin{alignat}{4}\label{C_w_ess_1}
||C_w||_{\mathcal{B}(L^2(\Gamma))}&\leq C \max\{||w^+||_{L^\infty}(\Gamma),||w^-||_{L^\infty}(\Gamma)\},\\\label{C_w_ess_2}
||C_w h||_{L^2(\Gamma)}&\leq C ||h||_{L^\infty(\Gamma)} \max\{||w^+||_{L^2(\Gamma)},||w^-||_{L^2(\Gamma)}\}
\end{alignat}
with $C=2\max \{||C_+||_{\mathcal{B}(L^2(\Gamma))},||C_-||_{\mathcal{B}(L^2(\Gamma))}\}$.
\end{subequations}
where
\[
C=2\max\{\|C_+\|_{\mathcal{B}(L^2(\Gamma))},\|C_-\|_{\mathcal{B}(L^2(\Gamma))}\}.
\]

The following lemma (\cites{L18}) establishes the equivalence between the
$L^2$ Riemann--Hilbert problem determined by $(\Gamma,v)$ and the singular
integral equation for $\mu\in I+L^2(\Gamma)$
\begin{equation}\label{mu_int_eq}
\mu-I=C_w(\mu).
\end{equation}

\begin{lemma}\label{Lem 5.2}
Given $v^\pm:\Gamma\to GL(n,\mathbb{C})$. Let $v=(v^-)^{-1}v^+$, $w^+=v^+-I$, $w^-=I-v^-$. Suppose $v^\pm$, $(v^\pm)^{-1}\in I + L^2(\Gamma)\cap L^\infty(\Gamma)$. If $m\in I + \dot E^2(D)$ satisfies $L^2-$RH problem determined by $(\Gamma,v(\lambda))$, then $\mu = m_+(v^+)^{-1}= m_-(v^-)^{-1}\in I+L^2(\Gamma)$ satisfies \eqref{mu_int_eq}. Conversely, if $\mu \in I+L^2(\Gamma)$ satisfies \eqref{mu_int_eq}, then $m=I+C(\mu(w^+ + w^-))\in I + \dot E^2(D)$ satisfies $L^2-$RH problem determined by $(\Gamma,v(\lambda))$.
\end{lemma}

The following lemma (\cites{L18}) shows that operator $\mathds{1}-C_w$ is Fredholm.

\begin{lemma}\label{Lem 5.3}
Given $v^\pm:\Gamma\to GL(n,\mathbb{C})$. Let $v=(v^-)^{-1}v^+$, $w^+=v^+-I$, $w^-=I-v^-$. Suppose $v^\pm$, $(v^\pm)^{-1}\in I + L^2(\Gamma)\cap L^\infty(\Gamma)$, and $v^\pm$ is continuous on $\Gamma$. Then
\begin{enumerate}[(1)]
    \item The operator $\mathds{1}-C_w:L^2(\Gamma)\to L^2(\Gamma)$ is Fredholm.
    \item If $w^\pm$ are nilpotent, then $\mathds{1}-C_w$ has Fredholm index $0$; in this case TFAE:
    \begin{enumerate}[(a)]
        \item  The map $\mathds{1}-C_w:L^2(\Gamma)\to L^2(\Gamma)$ is bijective.
        \item The $L^2-$RH problem determined by $(\Gamma,v(\lambda))$ has a unique solution.
        \item The homogeneous $L^2-$RH problem determined by $(\Gamma,v(\lambda))$ has only the zero solution.
        \item The map $\mathds{1}-C_w:L^2(\Gamma)\to L^2(\Gamma)$ is injective.
    \end{enumerate}
\end{enumerate}
\end{lemma}

\bibliographystyle{RS}
\bibliography{shepelsky_etal}

@article{IST01,
  title = {{Method for solving the Korteweg-de Vries equation}},
  author = {Gardner, Clifford S. and Greene, John M. and Kruskal, Martin D. and Miura, Robert M.},
  journal = {Physical Review Letters},
  volume = {19},
  issue = {19},
  pages = {1095--1097},
  numpages = {0},
  year = {1967},
  month = {Nov},
  publisher = {American Physical Society},
}

@article{IST02,
  title={Exact theory of two-dimensional self-focusing and one-dimensional self-modulation of waves in nonlinear media},
  author={V. E. Zakharov and A. B. Shabat},
  journal={Journal of Experimental and Theoretical Physics},
  year={1972},
  volume={34},
  pages={62--69}
}

@book{novikov1984theory,
  title={Theory of solitons: the inverse scattering method},
  author={Novikov, S and Manakov, Sergei V and Pitaevskii, Lev Petrovich and Zakharov, Vladimir Evgenevi{\v{c}}},
  year={1984},
  publisher={Springer Science \& Business Media}
}

@book{ablowitz1991solitons,
  title={Solitons, nonlinear evolution equations and inverse scattering},
  author={Ablowitz, Mark J and Ablowitz, MA and Clarkson, Peter A},
  volume={149},
  year={1991},
  publisher={Cambridge University Press}
}

@article{ablowitz1974inverse,
  title={{The inverse scattering transform-Fourier analysis for nonlinear problems}},
  author={Ablowitz, Mark J and Kaup, David J and Newell, Alan C and Segur, Harvey},
  journal={Studies in Applied Mathematics},
  volume={53},
  number={4},
  pages={249--315},
  year={1974},
  publisher={Wiley Online Library}
}

@book{trogdon2015riemann,
author = {Trogdon, Thomas and Olver, Sheehan},
title = {{Riemann--Hilbert problems, their numerical solution, and the computation of nonlinear special functions}},
publisher = {Society for Industrial and Applied Mathematics},
year = {2015},
address = {Philadelphia, PA}
}

@book{fokas2008unified,
author = {Fokas, Athanassios S.},
title = {{A unified approach to boundary value problems}},
publisher = {Society for Industrial and Applied Mathematics},
year = {2008},
}

@article{BSL17,
 author = {A.~{Boutet de Monvel} and D. Shepelsky  and L.  Zielinski},
 title = {{The short pulse equation by a Riemann–Hilbert approach}},
 journal = {Letters in Mathematical Physics},
 year = {2017},
 volume = {107},
 pages = {1345--1373},
}

@article{BKS20,
    author = {A.~{Boutet de Monvel} and I. Karpenko and D. Shepelsky},  
    title = {{A Riemann--Hilbert approach to the modified Camassa--Holm equation
   with nonzero boundary conditions}},
    journal = {J. Math. Phys.},
    year = 2020,
    volume = {61, 3},
    pages = {031504}
}

@article{fokas1997,
  title={{A unified transform method for solving linear and certain nonlinear PDEs}},
  author={Fokas, A. S.},
journal={Proceedings of the Royal Society of London. Series A: Mathematical, Physical and Engineering Sciences},
  volume={453},
  pages={61411--1443},
  year={1997}
}

@ARTICLE{BFS06,
  author={Boutet de Monvel, Anne and Fokas, Athanassios and Shepelsky, Dmitry},
  journal={Commun. Math. Phys.}, 
  title={Integrable nonlinear evolution equations on a finite interval}, 
  year={2006},
  volume={263},
  pages={133-172}
  }

@ARTICLE{FIS,
  author={Fokas, A. and Its, A. R. and {L-Y Sung} },
  journal={Nonlinearity}, 
  title={{The nonlinear Schr\"odinger equation on the half-line}}, 
  year={2005},
  volume={18},
  number={4},
  pages={1771}
  }

@ARTICLE{BS08,
  author={Boutet de Monvel, A. and Shepelsky, D.},
  title={{The Camassa--Holm Equation on the Half-Line: a Riemann--Hilbert Approach}},
  journal={J. Geom. Anal.},
  volume={18},
  year={2008},
  pages={285--323}
}

@ARTICLE{Z89,
  author={Zhou, X.},
  title={{The Riemann-Hilbert problem and inverse scattering}},
  journal={SIAM J. Math. Anal.},
  volume={20},
  year={1989}
}

@ARTICLE{SW04,
  author={Schäfer, T. and Wayne, C. E.},
  title={{Propagation of ultra-short optical pulses in nonlinear media}},
  journal={Physica D},
  volume={196},
  year={2004},
  pages={90--105}
}

@ARTICLE{CJWS05,
  author={Chung, Y. and Jones, C. K. R. T. and Schäfer, T. and Wayne, C. E.},
  title={{Ultra-short pulses in linear and nonlinear media}},
  journal={Nonlinearity},
  volume={18},
  number={3},
  year={2005},
  pages={1351--1374}
}

@ARTICLE{SS05,
  author={Sakovich, A. and Sakovich, S.},
  title={{The short pulse equation is integrable}},
  journal={J. Phys. Soc. Japan},
  volume={74},
  year={2005},
  pages={239--241}
}

@ARTICLE{KS25,
  author={Karpenko, I. and Shepelsky, D.},
  title={{The modified Camassa-Holm equation on the half line: a Riemann--Hilbert approach}},
  journal={arXiv:2512.05274},
year={2025}
}

@ARTICLE{Gambino2015,
  author={Gambino, G. and Tanriver, U. and Guha, P. and Choudhury, A. and Choudhury, S.},
  title={{Regular and singular pulse and front solutions and possible isochronous behavior in the short-pulse equation: phase-plane, multi-infinite series and variational approaches}},
  journal={Commun. Nonlinear Sci. Numer. Simul.},
  volume={20},
  number={2},
  year={2015},
  pages={375--388}
}

@ARTICLE{Matsuno2007,
  author={Matsuno, Y.},
  title={{Multiloop soliton and multibreather solutions of the short pulse model equation}},
  journal={J. Phys. Soc. Japan},
  volume={76},
  number={8},
  year={2007},
  pages={084003}
}

@ARTICLE{Sakovich2006,
  author={Sakovich, A. and Sakovich, S.},
  title={{Solitary wave solutions of the short pulse equation}},
  journal={J. Phys. A: Math. Gen.},
  volume={39},
  year={2006},
  pages={L361--L367}
}

@ARTICLE{Matsuno2008,
  author={Matsuno, Y.},
  title={{Periodic solutions of the short pulse model equation}},
  journal={J. Math. Phys.},
  volume={49},
  number={7},
  year={2008},
  pages={073508}
}

@ARTICLE{Coclite2015,
  author={Coclite, G. M. and di Ruvo, L.},
  title={{Well-posedness results for the short pulse equation}},
  journal={Z. Angew. Math. Phys.},
  volume={66},
  year={2015},
  pages={1529--1557}
}

@ARTICLE{Pelinovsky2010,
  author={Pelinovsky, D. and Sakovich, A.},
  title={{Global well-posedness of the short-pulse and sine-Gordon equations in energy space}},
  journal={Comm. Partial Differential Equations},
  volume={35},
  number={4},
  year={2010},
  pages={613--629}
}

@ARTICLE{Liu2009,
  author={Liu, Y. and Pelinovsky, D. and Sakovich, A.},
  title={{Wave breaking in the short-pulse equation}},
  journal={Dyn. PDE},
  volume={6},
  year={2009},
  pages={291--310}
}

@article{ZH89,
  author  = {X. Zhou},
  title   = {{The Riemann--Hilbert problem and inverse scattering}},
  journal = {SIAM Journal on Mathematical Analysis},
  volume  = {20},
  pages   = {966--986},
  year    = {1989}
}

@article{L18,
  author  = {Jonatan Lenells},
  title   = {{Matrix Riemann--Hilbert problems with jumps across Carleson contours}},
  journal = {Monatshefte für Mathematik},
  volume  = {186},
  number  = {1},
  pages   = {111--152},
  year    = {2018}
}

@book{TFA,
  author    = {Gerald Teschl},
  title     = {{Topics in Linear and Nonlinear Functional Analysis}},
  publisher = {American Mathematical Society},
  address   = {Providence, RI},
  year      = {to appear}
}

@article{L12,
  author  = {Lenells, Jonatan},
  title   = {{Nonlinear Fourier Transforms and the mKdV Equation in the Quarter Plane}},
  journal = {Studies in Applied Mathematics},
  volume  = {129},
  number  = {4},
  pages   = {347--371},
  year    = {2012},
  doi     = {10.1111/j.1467-9590.2012.00560.x}
}

@article{Zhou98,
  author    = {Xin Zhou},
  title     = {{$L^2$-Sobolev space bijectivity of the scattering and inverse scattering transforms}},
  journal   = {Communications on Pure and Applied Mathematics},
  volume    = {51},
  number    = {7},
  pages     = {697--855},
  year      = {1998},
  publisher = {John Wiley \& Sons},
}

\end{document}